\documentclass[12pt,reqno,twoside]{amsart}
\usepackage{amsmath}
\usepackage{amssymb}
  
  \newcommand{\grad}{\nabla}

\newcommand{\vre}{\varepsilon}

\newcommand\vare{\varepsilon}

\newcommand\CC{\mathbb C}
\newcommand\HH{\mathbb H}
\newcommand\NN{\mathbb N}
\newcommand\RR{\mathbb R}
\newcommand\ZZ{\mathbb Z}
\newcommand\QQ{\mathbb Q}
\newcommand\TT{\mathbb T}

\newcommand\FF{\mathbb F}

\newcommand\cH{\mathcal{H}}

\newcommand\cJ{\mathcal{J}}

\newcommand\cU{\mathcal{U}}

\renewcommand\Re{\operatorname{Re}}
\renewcommand\Im{\operatorname{Im}}
\newcommand\supp{\operatorname{supp}}

\newcommand\Ind{\operatorname{Ind}}

\newcommand\vol{\operatorname{vol}}
\newcommand\tr{\operatorname{tr}}

\newcommand\norm[1]{\left\|#1\right\|}

\newcommand\abs[1]{\left|#1\right|}

\newcommand\inn[1]{\left\langle #1 \right\rangle}
\newcommand\set[1]{\left\{{#1}\right\}}

\newtheorem*{claim}{Claim}
\newtheorem{step1}{Step}

\newtheorem{theorem}{Theorem}[section]
\newtheorem{lemma}[theorem]{Lemma}
\newtheorem{proposition}[theorem]{Proposition}

\theoremstyle{definition}
\newtheorem{definition}[theorem]{Definition}

\theoremstyle{remark}
\newtheorem{remark}[theorem]{Remark}

\theoremstyle{corollary}
\newtheorem{corollary}[theorem]{Corollary}

\numberwithin{equation}{section}

\begin{document}

\title{The ergodic theory of lattice subgroups}

\author{Alexander Gorodnik}
\address{University of Bristol, UK}
\email{a.gorodnik@bristol.ac.uk}
\thanks{The first author was supported in part by NSF Grant}

\author{Amos Nevo}
\address{Department of Mathematics, Technion IIT}
\email{anevo@tx.technion.ac.il}
\thanks{The second author was supported in part by the Institute for Advanced Study and 
an ISF grant.}


\date{Final version, September  2007}



\begin{abstract}
We prove mean and pointwise ergodic theorems for general families 
of averages on a semisimple algebraic (or $S$-algebraic) 
group $G$, together with 
an explicit rate of convergence when the 
action has a spectral gap. 
Given any lattice $\Gamma$ in $G$, we use the 
ergodic theorems for $G$ to 
solve the lattice point 
counting problem 
for general domains in $G$,  
and prove mean and pointwise ergodic theorems for arbitrary measure-preserving 
actions of the lattice, together with 
explicit rates of convergence when a spectral gap is present. 
 We also prove an equidistribution theorem in 
arbitrary isometric actions of the lattice.

For the proof we develop a general method 
to derive ergodic theorems 
for actions of a locally compact group $G$, and 
of a lattice subgroup $\Gamma$, 
provided certain natural spectral, 
geometric and regularity 
conditions are satisfied by the group 
$G$, the lattice $\Gamma$, and the domains where 
the averages are supported. In particular, we establish the 
general principle that under these conditions a quantitative 
mean ergodic theorem in $L^2(G/\Gamma)$ 
for a family of averages gives rise to a quantitative solution 
of the lattice point counting problem in their supports.
 We demonstrate the new explicit error terms that
 we obtain by a variety of examples.  
\end{abstract}

\maketitle
{\small \tableofcontents}

\section{Main results : Semisimple Lie groups case}

\subsection{Admissible sets}

Let $G$ be a locally compact second countable (lcsc) group, 
and $\Gamma\subset G$ a lattice subgroup.
 Consider the following four fundamental problems in ergodic theory 
that present themselves in this context, namely : 
\begin{enumerate}
\item  Prove 
ergodic theorems for general families of averages on $G$, 
 \item  Solve the lattice point counting 
problem (with explicit error term) for any lattice subgroup 
$\Gamma$ and for general domains on $G$, 
\item Prove 
ergodic theorems for {\it arbitrary} actions of a lattice subgroup  
$\Gamma$,
\item  Establish equidistribution results 
for isometric actions of the lattice $\Gamma$.
\end{enumerate}

Our purpose in the present paper is to give a complete  
solution to these problems for non-compact 
semisimple algebraic groups over arbitrary local fields,  
and any of their lattices. Our results apply also to lattices in 
products of such groups, and thus also to $S$-algebraic groups 
and their lattices.  
In fact, many 
of our arguments hold in greater generality 
still, and we will elaborate on that further in our discussion below. 
However, for simplicity 
of exposition we will begin by describing the 
main results, as well as some of their applications,
 in the case of connected semisimple Lie groups.


We start by introducing 
the following definition, which describes the families $\beta_t$ 
that will be the subject of our analysis.

Fix any 
left-invariant Riemannian metric on $G$, and let 
$$
\mathcal{O}_\vre=\{g\in G:\, d(g,e)<\vre\}.
$$
Let $m_G$ denote a fixed left Haar measure on $G$. 

\begin{definition}\label{admissible-Lie}
 {\rm
An increasing family of bounded Borel subsets $G_t$, $t>0$, of $G$ will be called 
{\it admissible} if 
there exists $c>0$ such that for all $t$ sufficiently large and $\vre$ sufficiently small 
\begin{align}
\mathcal{O}_\vre\cdot G_t\cdot \mathcal{O}_\vre &\subset G_{t+c\vre},\label{eq:int1}\\
m_G(G_{t+\vre})&\le (1+c\vre)\cdot  m_G(G_{t}).\label{eq:int2}
\end{align}
}
\end{definition}

Let us briefly note the following facts (see Prop. \ref{Riemannian}
 and Prop. \ref{stability} below, as well as the Appendix for the proof).
\begin{enumerate}
\item Admissibility is independent of the Riemannian metric chosen to define it. 
\item  Many of the natural 
families of sets in $G$ are admissible. In particular 
the radial sets $B_t$ projecting to the Cartan-Killing Riemannian 
balls on the symmetric space are admissible. Furthermore, the sets 
$\set{g\,;\, \log\norm{\tau(g)}< t}$ 
where $\tau$ is faithful linear representation 
are also admissible, for any choice of linear norm $\norm{\cdot}$.

\item Admissibility is 
invariant under translations, namely if $G_t$ is admissible, 
so is $gG_t h$, for any fixed $g,h\in G$.      
\end{enumerate}

It is natural to define also the corresponding H\"older conditions. As we shall see below, 
whenever a spectral gap is present, the assumption of admissibility can be weakened 
to H\"older admissibility.

\subsection{Ergodic theorems on 
semisimple Lie groups}

We define $\beta_t$ to be the probability
 measures on $G$ obtained as the restriction 
of Haar measure to $G_t$, normalized by $m_G(G_t)$. 

The averaging operators associated to $\beta_t$ when $G$ acts by measure-preserving 
transformations of a probability space $(X,\mu)$ are given by 
$$\pi(\beta_t)f(x)= \frac{1}{m_G(G_t)}\int_{G_t} 
f(g^{-1}x)dm_G(g) \,\,.$$

Assume $G$ is connected semisimple with finite center and no compact factors. Then 
\begin{enumerate}
\item
The family $\beta_t$ (and $G_t$) will be called 
(left-) radial if it is invariant under 
(left-) multiplication by some fixed 
maximal compact subgroup $K$, for all sufficiently large $t$. 
 Standard radial averages are those defined in Definition \ref{standard}.   
\item The action is called irreducible 
if every non-compact simple factor acts ergodically. 
\item The action is said to have a strong spectral gap if 
each simple factor has a spectral gap, namely 
admits no asymptotically invariant sequence of unit vectors (see 
\S 3.6 for a full discussion). 
\item The sets $G_t$ (and the averages $\beta_t$) will be called balanced if for every simple 
factor $H$ and every compact subset $Q$ of its complement, $ \beta_t(Q H)\to 0$.  
$G_t$ will be called well-balanced if the convergence 
 is at a specific rate (see \S 3.5 for a full discussion). 
  \end{enumerate}

Our first main result is the following pointwise ergodic theorem 
for admissible averages 
on semisimple Lie groups.

\begin{theorem}\label{th:Lie groups}
{\bf Pointwise ergodic theorems for admissible averages.}
 Let $G$ be a connected 
semisimple Lie group with finite center and no 
non-trivial compact factors. Let  
$(X,\mu)$ be a standard Borel space with a 
probability-measure-preserving ergodic action of $G$. Assume  
that $G_t$ is an admissible family. 
\begin{enumerate} \item 
Assume that $\beta_t$ is left-radial. 
 If the action is irreducible, 
then $\beta_t$ satisfies the pointwise ergodic theorem 
in $L^p(X)$, $ 1 < p < \infty$, namely for every 
$f\in L^p(X)$, and for almost every $x\in X$  :
$$\lim_{t\to\infty}\pi(\beta_t)f(x)=\int_X fd\mu\,. $$ 
The conclusion holds also in reducible actions of $G$, 
provided the averages are standard radial, well-balanced and boundary-regular 
 (see \S\S  3.4, 3.5 for the definitions).

\item If the action has a strong spectral gap, 
then $\beta_t$ converges  
to the ergodic mean almost surely exponentially fast, 
namely for every 
$f\in L^p(X)$, $1< p \le \infty$, 
and almost all $x\in X$ 
$$\abs{\pi(\beta_t)f(x) -\int_X fd\mu}\le 
C_p(f,x)e^{-\theta_p t}\,,$$
where $\theta_p > 0$ depends explicitly on the spectral gap 
(and the family $G_t)$.

The conclusion holds also in  actions of $G$ with a spectral gap,  
provided the averages satisfy the additional 
necessary condition of being well-balanced 
(see \S\S 3.5, 3.7 for the definitions).  

\end{enumerate}

\end{theorem}

Regarding Theorem \ref{th:Lie groups}(1), 
we remark that the
 proof of pointwise convergence in the case of reducible actions without a 
spectral gap is quite involved, and we have thus assumed 
in that case that the averages are standard radial, well-balanced and boundary-regular to make the analysis tractable. However, the reducible case will be absolutely indispensable for us below,  since we will induce actions 
of a lattice subgroup to actions of $G$, and these may be reducible. 

Regarding Theorem \ref{th:Lie groups}(2), 
we note that $\theta_p$ depends explicitly
 on the spectral gap of the action, and on natural 
geometric parameters of $G_t$, and we refer to \S 7.1 
for a full discussion including a formula for a lower bound. Furthermore, 
H\"older admissibility is sufficient for this part, as we will see below.

Let us now formulate
 the following invariance principle for ergodic 
actions of $G$, which will play an important role below, 
in the derivation of pointwise ergodic theorems for lattices.

\begin{theorem}\label{th:invariance}
{\bf Invariance principle.}
Let $G$, $(X,\mu)$ be as in Theorem 
\ref{th:Lie groups}, and let $G_t$ be an admissible family.  
Then for any given function $f\in L^p(X)$ the set
 where pointwise convergence to the ergodic mean holds, namely 
$$\set{x\in X\,;\, \lim_{t\to\infty} 
\frac{1}{m_G(G_t)}\int_{G_t} f(g^{-1}x)dm_G(g)=\int_X fd\mu}$$  
contains a $G$-invariant set of full measure.
\end{theorem}

We note that $G$ is a non-amenable group, and the sets $G_t$ are 
not asymptotically invariant under translations (namely do not 
have the F\o lner property). Thus the conclusion of Theorem 
\ref{th:invariance} is not obvious, even in the case where 
$X$ is a 
homogeneous $G$-action. The special case where $G=SO^0(n,1)$ and 
$\beta_t$ are the bi-$K$-invariant averages lifted from ball averages on  
hyperbolic space $\HH^n$ was considered earlier by \cite{BR}.

One of our applications of ergodic theorems on $G$ is to
the lattice point counting problem in $G_t$. The solution of 
the latter actually depends only 
on the {\it mean } ergodic theorem 
for $\beta_t$, which holds under more general conditions than 
the pointwise theorem.
Because of its later significance, 
we therefore formulate separately the following 
\begin{theorem}\label{th: Mean Lie groups}
{\bf Mean ergodic theorems for admissible averages.}
Let $G$ and $(X,\mu)$ be as in Theorem \ref{th:Lie groups}, and 
let $G_t$ be an admissible family. 
\begin{enumerate}
\item If the action is irreducible or $G_t$ are balanced, then 
$$\lim_{t\to \infty} \norm{\pi(\beta_t)f-\int_X fd\mu}_{L^p(X)}
=0\,\,\,\,, 1\le p < \infty\,\,.$$
\item If the action has a strong spectral gap, or a spectral gap and the averages 
are well balanced, then 
$$ \norm{\pi(\beta_t)f-\int_X fd\mu}_{L^p(X)}
\le B_p e^{-\theta_p t}\,\,\,\,, 1 < p < \infty\,\,$$
for the same $\theta_p > 0$ as in Theorem \ref{th:Lie groups}(2).  
\end{enumerate}  
\end{theorem}

\subsection{The lattice point counting problem in 
admissible domains}

Let now 
$\Gamma\subset G$ be any lattice subgroup; the lattice point 
counting problem is to determine the 
number of lattice points in the domains $G_t$. Its 
ideal solution calls for
 evaluating the main term in the asymptotic expansion, 
 establishing the existence of the limit,
and estimating explicitly the error term. Our second main result 
gives a complete solution to this problem for all lattices and 
all families of admissible domains. The proof we give below will 
establish the general principle asserting that a mean ergodic 
in $L^2(G/\Gamma)$ for the averages $\beta_t$ 
 (with explicit rate of convergence) implies  
a solution to the $\Gamma$-lattice point counting problem 
in the admissible domains $G_t$ (with an explicit 
estimate of the error term). We will show 
below that under certain natural assumptions 
this principle can be established in great generality 
for lattices in general lcsc groups, but 
will state it first for connected semisimple Lie groups. 

We note that in this case, the main term in the lattice count (namely part (1) of the following theorem)
 was established  \cite{Ba} (for uniform lattices),
\cite{DRS} (for balls w.r.t. a norm) and \cite{EM} (in general). 
Error term were considered for rotation-invariant norms  
in \cite{DRS} and for more general norms very recently in \cite{Ma}.  
For a comparison of part (2) of the following theorem 
with these results see \S 2.

\begin{theorem}\label{th:Lattice points}
{\bf Counting lattice points in admissible domains.}
Let $G$ be a 
connected semisimple Lie group with finite center, and no 
non-trivial compact factors. 
Let $G_t$ be an admissible family of sets, and let 
$\Gamma$ be any lattice subgroup. Normalize Haar measure $m_G$ 
to assign measure one to a fundamental domain of $\Gamma$ in $G$. 
\begin{enumerate}\item If $\Gamma$ is an irreducible lattice, 
or the sets $G_t$ are balanced, then   
 $$\lim_{t\to \infty} \frac{\abs{\Gamma\cap G_t}}{m_G(G_t)}=1\,\,.$$

\item If $(G/\Gamma,m_{G/\Gamma})$  has a strong spectral gap, 
or the sets $G_t$ are well balanced, then, for all $\vre > 0$ 
$$ \frac{\abs{\Gamma\cap G_t}}{m_G(G_t)}
= 1+O_\vre\left( \exp\left(\frac{-t(\theta-\vre)
 }{\dim G +1}\right)\right)\,,$$
where $\theta >0$ depends on $G_t$ and the spectral gap 
in $G/\Gamma$, via 
$$\theta = \liminf_{t\to \infty}-\frac1t\log 
\norm{\pi_{G/\Gamma}(\beta_t)}_{L^2_0(G/\Gamma)}\,.$$ 

\end{enumerate}
\end{theorem}

\begin{remark}
\begin{enumerate}
\item Recall that the $G$-action on $(G/\Gamma, m_{G/\Gamma})$ is 
irreducible if 
and only if $\Gamma$ is an 
irreducible lattice in $G$, namely the projection of
 $\Gamma$ to every simple factor of $G$ is a 
dense subgroup.
\item The $G$-action on $G/\Gamma$ always has 
a spectral gap, but whether it has a strong spectral gap 
seems to be an open problem, in general (see \S 3.5 for more details). 
\item When the action has a strong spectral gap, the parameter 
$\theta$ can be given explicitly in terms of the rate of volume growth 
of the sets $G_t$ and the size of the gap - see Remark \ref{KS} and \S 7.1.   
\item Note that under the normalization of $m_G$ given in Theorem 
\ref{th:Lattice points}, if $\Delta\subset \Gamma$ 
is a subgroup of finite index, then 
 $$\lim_{t\to \infty} \frac{\abs{\Delta\cap G_t}}{m_G(G_t)}=
\frac{1}{[\Gamma:\Delta]}\,.$$
\end{enumerate}
\end{remark}

Finally, we remark that the condition of 
admissibility is absolutely crucial in obtaining pointwise ergodic theorems for $G$, and thus also for $\Gamma$. This is true when the action does not have a spectral gap, but also when it does (although here H\"older-admissibility  is sufficient). However, lattice point counting results, quantitative or not, hold in significantly greater generality. Namely, it holds for families that satisfy the weaker condition $m_G(\mathcal{O}_\vre G_t \mathcal{O}_\vre) \le (1+c\vre) m_G( G_t)$, which amounts to a quantitative version of the well-roundedness condition of \cite{DRS} and \cite{EM}. This generalization is discussed systematically in \cite{GN}, where several applications, including to quantitative counting of lattice points in sectors, on symmetric varieties and on Adele groups are given.

\subsection{Ergodic theorems for lattice subgroups}

We now turn to our third main result, namely 
to the solution of the problem of 
establishing ergodic theorems for a 
general action of a lattice 
subgroup on a probability space $(X,\mu)$.
 This result also uses Theorem \ref{th:Lie groups} 
as a basic tool; 
here it is applied to 
the action of $G$ induced by the action of $\Gamma$ on 
$(X,\mu)$. This argument generalizes the
 one used in the proof of 
Theorem \ref{th:Lattice points}, where we considered 
the action of $G$ induced from the trivial action of 
$\Gamma$ on a point. However the increased generality 
requires a considerable number 
of additional further arguments. 

To formulate the result, consider the set of
 lattice points $\Gamma_t=\Gamma\cap G_t$. 
Let $\lambda_t$ denote the probability measure on $\Gamma$ 
uniformly distributed on $\Gamma_t$.

We begin with the following fundamental 
 mean ergodic theorem   
for arbitrary lattice actions.  

\begin{theorem}\label{mean}{\bf Mean ergodic theorem for 
lattice actions.}

Let $G$, $G_t$ and $\Gamma$, be as in  
Theorem \ref{th:Lattice points}. Let $(X,\mu)$ 
 be an ergodic measure-preserving action of $\Gamma$.   
\begin{enumerate}
\item Assume the action 
of $G$ induced from the $\Gamma$-action on 
$(X,\mu)$ is irreducible, or that $G_t$ are balanced. 
Then for every $f\in L^p(X)$, $ 1 \le p < \infty$, 
$$ \lim_{t\to\infty}\norm{\frac{1}{\abs{\Gamma_t}}\sum_{\gamma\in
 \Gamma_t} 
f(\gamma^{-1}x)-\int_X f d\mu}_{L^p(X)}=0\,\,.$$
\item  Assume that the action of $G$ induced from 
the $\Gamma$-action on 
$(X,\mu)$ has a strong spectral gap, 
or that it has a spectral 
gap and $G_t$ are well balanced. 
Then for every $f\in L^p(X)$, $1 < p < \infty$ 
$$\norm{\frac{1}{\abs{\Gamma_t}}\sum_{\gamma\in \Gamma_t} 
f(\gamma^{-1}x)-\int_X f d\mu}_{L^p(X)}
\le C e^{-\delta_p t}\norm{f}_{L^p(X)}\,\,,$$
where $\delta_p$ is determined explicitly by the  
spectral gap for the induced $G$-action (and depends also on 
the family $G_t$).
\end{enumerate}

\end{theorem}

One immediate application of Theorem \ref{mean} 
arises when we take 
$X$ to be a transitive action 
on a finite space, namely $X=\Gamma/\Delta$, $\Delta$ 
a finite index subgroup.

\begin{corollary}\label{th:finite}
{\bf Equidistribution in finite actions.}
Let $G$, $\Gamma$ and $G_t$ be as 
in Theorem \ref{th:Lattice points}. 
Let $\Delta\subset \Gamma $ 
be a subgroup of finite index, and $\gamma_0$ any element in $ \Gamma$. 
\begin{enumerate}
\item Under 
the assumptions of Theorem \ref{th:Lattice points}(1) 
$$
\lim_{t\to \infty}\frac{1}{\abs{\Gamma_t}}\cdot
\abs{\set{\gamma\in\Gamma\cap G_t:\,\gamma \cong \gamma_0 
\text{ mod } \Delta }}= 
\frac{1}{[\Gamma: \Delta]}\,.
$$
\item Under the assumptions of Theorem \ref{mean}(2) 
$$
\frac{1}{|\Gamma_t|}
\cdot|\{\gamma\in\Gamma\cap G_t:\,\gamma\cong
 \gamma_0 \text{ mod } \Delta\}= \frac{1}{[\Gamma: \Delta]}+
O(e^{-\delta t})
$$
where $\delta > 0$, and is     
determined explicitly by the spectral gap in $G/\Delta$.  
\end{enumerate}
\end{corollary}

We remark that the first conclusion 
in Corollary \ref{th:finite}, namely equidistribution 
of the lattice points in $\Gamma\cap G_t$ among the cosets of $\Delta$ in $ \Gamma$  
can also be obtained using the method of [GW], 
which employs Ratner's theory of
unipotent flow. It is also possible to derive 
this result from 
considerations related to the 
mixing property of flows on $G/\Gamma$. 

Another application of the mean ergodic theorem  
 is in the proof of an equidistribution theorem  
for the corresponding averages in isometric 
actions of the lattice. The result is as follows.

\begin{theorem}\label{th:equidistribution}
{\bf Equidistribution 
in isometric 
actions of lattices.} 
Let $G$, $G_t$, and $\Gamma$ be as 
in Theorem \ref{mean}. 
Let $(S,d)$ be 
a compact metric space on which $\Gamma$ 
acts by isometries, and assume the action is 
ergodic with respect to an invariant probability measure 
$\mu$ whose support coincides with $S$. 
Then under the assumptions of Theorem \ref{mean}(1), 
for every continuous function $f$ on $S$ and {\it every point}
 $s\in S$
$$ \lim_{t\to\infty}\frac{1}{\abs{\Gamma_t}}
\sum_{\gamma\in \Gamma_t} 
f(\gamma^{-1}s)=\int_S f d\mu\,\,.$$
and the convergence is uniform in $s\in S$ 
(i.e. in the supremum norm on $C(S)$). 
\end{theorem}

Let us now formulate pointwise ergodic theorems for 
general actions of lattices.

\begin{theorem}\label{th:Lattice subgroups}
{\bf Pointwise ergodic theorems for general lattice actions.}
Let $G$, $G_t$, $\Gamma$ and 
$(X,\mu)$ be as in Theorem \ref{mean}.

\begin{enumerate} \item 
Assume that the action induced to $G$ is irreducible, and $\beta_t$ are left-radial. 
Then the averages $\lambda_t$ satisfy the 
pointwise ergodic theorem in $L^p(X)$, 
$1 < p < \infty$, namely 
for $f\in L^p(X)$ and almost every $x\in X$ :
$$ \lim_{t\to\infty}\frac{1}{\abs{\Gamma_t}}\sum_{\gamma\in \Gamma_t} 
f(\gamma^{-1}x)=\int_X f d\mu\,\,$$
The same conclusion also holds when the induced action is reducible, provided $\beta_t$ are 
standard radial, well-balanced and boundary-regular. 
\item Retain the assumption of Theorem \ref{mean}(2). 
Then the convergence of $\lambda_t$ to the ergodic mean 
is almost surely exponentially fast, 
namely for $f\in L^p(X)$, $1 < p < \infty$  
and almost every $x\in X$   
$$\abs{\frac{1}{\abs{\Gamma_t}}\sum_{\gamma\in \Gamma_t} 
f(\gamma^{-1}x)-\int_X f d\mu}
\le C_p(x,f)e^{-\zeta_p t}$$
where $\zeta_p$ is determined explicitly by the  
spectral gaps for the induced $G$-action (and the family $G_t$).
\end{enumerate}
\end{theorem}

\begin{remark}\label{lattice-pointwise}

\begin{enumerate}
\item Note that if $G$ is simple, then of course any 
action of $G$ induced from an ergodic action of a lattice 
subgroup is irreducible. However, if $G$ is not simple, then 
the induced action can be reducible and then the assumption 
that the averages are balanced is necessary in Theorem 
\ref{th:Lattice subgroups}(1). We assume in fact that they are standard 
radial, well-balanced and boundary-regular, as we will apply
Theorem \ref{th:Lie groups}(1) to the induced action.  
\item Note further that if $G$ is simple and 
has property $T$, then the assumption 
of strong spectral gap stated in Theorem 
\ref{th:Lattice subgroups}(2) 
is satisfied for every ergodic action of every 
lattice subgroup. Furthermore, in that case  
$\zeta_p$ has an explicit positive lower  
 bound depending on $G$ and $G_t$ only and independent of 
$\Gamma$ and $X$. 
\item It may be the case that whenever
$G/\Gamma$ has a strong spectral gap, 
so does every action of $G$ induced 
from an ergodic action of the 
irreducible lattice $\Gamma$ which has a spectral gap, but 
this problem also seems to be open. 
\item As we shall see in \S 6.1, the possibility of 
utilizing the induced $G$-action 
to deduce information on {\it pointwise convergence} 
in the inducing $\Gamma$-action depends 
on the invariance 
principle stated in Theorem \ref{th:invariance}   
for admissible averages on $G$.

\end{enumerate}
\end{remark} 

{\it On the scope of the method}. 
In light of remarks (1) and (2) above, 
let us explain the reason we avoided the (considerable) temptation 
to retrict our attention to simple groups and their lattices. 
First, such a restriction 
rules out of course a solution to the lattice point counting 
problem even for such natural examples as $SL_2(\ZZ[\sqrt{2}])$, 
which is a lattice in $SL_2(\RR)\times SL_2(\RR)$.  
Second, even 
for certain lattices in the latter group, 
the existence of strong spectral gap in $G/\Gamma$ is unknown (see \S 3.7). 
Thus when considering lattice points on 
product groups, whether the spectral gap is strong and the averages balanced or well-balanced 
become necessary considerations. 

Third, 
we have formulated our ergodic theorems for $G$ also in the 
case of reducible actions, but this again is unavoidable. Indeed, 
the ergodic theorems for the lattices are proved by 
induction to $G$, and it is unknown when the resulting action is 
irreducible (it is when the $\Gamma$-action is mixing or isometric
 \cite{St}). Finally, in order to handle such obvious examples as 
$SL_2(\ZZ[\frac1p])$ (which is a lattice in 
$SL_2(\RR)\times SL_2(\QQ_p)$) it is necessary to extend 
the theory to include $S$-algebraic groups, a task we will take on below. 

Below we will give a complete analysis valid 
for $S$-algebraic groups and their lattices in all cases, but let us here 
demonstrate our results in a more concrete fashion, which shows, 
in particular, that sets $\Gamma_t$ satisfying all the assumptions required do exist. 
 Indeed, let $G$ be a connected semisimple Lie group 
with finite center and no compact factors. Let $G/K$ be its symmetric space and $d$ the 
Riemannian distance associated with the Cartan-Killing form, and  let 
$B_t=\set{g\in G\,;\, d(gK,K)\le t}$, $\beta_t$ be the Haar-uniform 
averages. Then $G_t$ are admissible and well-balanced, and it has been established in 
\cite{N1}\cite{N2}\cite{NS1}\cite{MNS} that 
in every ergodic probability measure preserving action 
of $G$, the family $\beta_t $ satisfies the pointwise ergodic 
theorem in $L^p$, $1 < p < \infty$.  Furthermore,  if the action 
has a spectral gap, then the convergence to the ergodic 
mean is exponentially fast, as in Theorem \ref{th:Lie groups}(2).


Now let $\Gamma\subset G$ be any lattice subgroup. Then the following result, announced in \cite[Thm. 14.4]{N5}, holds.

\begin{theorem}\label{killing}{\bf Ergodic theorems for lattice points in 
Riemannian balls.}
Let $G$, $B_t$ and $\Gamma$ be as in the preceding paragraph, 
 and $\lambda_t$ the uniform averages on $\Gamma\cap B_t$. Then in every probability 
measure-preserving action of $\Gamma$, $\lambda_t$ satisfy 
the mean ergodic theorem in $L^p$, $1 \le p < \infty$ and the 
pointwise ergodic theorem in $L^p$, $ 1 < p \le \infty$. 
If the $\Gamma$-action has a spectral gap, then $\lambda_t$ 
satisfy the exponentially fast mean and pointwise ergodic theorem 
as in Theorem \ref{mean}(2) and Theorem 
\ref{th:Lattice subgroups}(2). Finally, $\lambda_t$ satisfy 
the equidistribution theorem w.r.t. an ergodic invariant 
probability measure of full support 
in every isometric action of $\Gamma$. 
\end{theorem}

As is clear from the statements 
of the foregoing theorems, the distinction between 
actions with and without a 
spectral gap is fundamental in determining which 
ergodic theorems apply, and the 
two cases call for rather different methods of proof. Thus 
the results will be established
 according to the following scheme :
\begin{enumerate}
\item Ergodic 
theorems for general averages on semisimple $S$-algebraic groups 
in the presence of a spectral gap.
\item Ergodic 
theorems for general averages on semisimple $S$-algebraic 
 groups in the absence of a spectral gap. 
\item Stability of admissible averages on semisimple 
$S$-algebraic groups, and 
an invariance principle for their ergodic actions.
\item Mean, maximal and pointwise Ergodic theorems
 for lattice subgroups, in the absence 
of a spectral gap. 
\item Exponentially fast 
pointwise ergodic theorem for lattice actions 
in the presence of a spectral gap. 
\item Equidistribution for isometric lattice actions. 

\end{enumerate}

As we shall see below, this scheme applies in a much 
wider context than that of semisimple $S$-algebraic groups. 
We will formulate it in \S\S 5.2 and  6.2
 below as a general 
recipe to derive ergodic theorems 
for actions of an lcsc group, and 
of a lattice subgroup $\Gamma$, 
provided certain natural spectral, 
geometric and regularity 
conditions are satisfied by the group 
$G$, the lattice $\Gamma$, and the sets $G_t$.

\section{Examples and applications}

Let us now consider some concrete 
examples and applications of the results 
stated above, and compare our results 
to some precedents in the 
literature. 

\subsection{Hyperbolic lattice points problem}

We begin by applying Theorem 
\ref{th:Lattice points} to the 
classical lattice point counting problem in hyperbolic space. Let us call a lattice subgroup $\Gamma$ tempered if the spectrum of the representation of the isometry group in $L^2_0(G/\Gamma)$ is tempered.  

\begin{corollary}\label{hyp-space}

Let $\HH^n$ be hyperbolic $n$-space taken with constant curvature 
$-1$ and the resulting volume form. Let $B_t$ be  
the Riemannian balls centered at a given point, and let $\Gamma$ be any lattice. Then 
\begin{enumerate}
\item 
$$\frac{\abs{\Gamma\cap B_t}}{\vol(B_t)}=\frac{1}{\vol(\HH^n/\Gamma)}+O_\vre \left(
\exp -t\left(\frac{\theta}{n+1}-\vre\right)\right)$$
provided $\norm{\pi_{G/\Gamma}(\beta_t)}_{L^2_0(G/\Gamma)}
\le C^\prime e^{-t(\theta-\vre) }$.
\item  In particular, if $\Gamma$ is tempered, then 
$$\frac{\abs{\Gamma\cap B_t}}{\vol(B_t)}=\frac{1}{\vol(\HH^n/\Gamma)}+O_\vre \left(
\exp -t\left(\frac{n-1}{2(n+1)}-\vre\right)\right)\,.$$
\end{enumerate}
\end{corollary}

We remark that the bound stated above is actually 
better than that provided by Theorem \ref{th:Lattice points}, 
as the error term here is given by $\theta/(\dim G/K +1)$ rather than 
$\theta/(\dim G +1)$. This is a consequence of the fact that we have 
taken here bi-$K$-invariant averages on $G$, 
so that the arguments 
used in the proof of Theorem  \ref{th:Lattice points} 
can be applied on $G/K$ rather than $G$. 
The same bound holds for any choice of
bi-$K$-invariant admissible sets $G_t$. The spectral gap 
parameter is given by $\theta=\frac12(n-1)$ in the 
tempered case, 
since the convolution norm of $\beta_t$ on $L^2(G)$ 
is dominated 
by $\vol(B_t)^{-1/2+\vre}$ (see Remark \ref{KS}), 
which is asymptotic to
$ \exp -t\left(\frac12 (n-1)-\vre\right)$ (recall that $\vol B_t$ is asymptotic to $c_ne^{(n-1)t}$).  





 For comparison, the best existing bound for a tempered   
lattice in hyperbolic $n$-space ($n \ge 2$) is due to 
Selberg \cite{Se} and Lax and Phillips \cite{LP}, and is given by 
$$\frac{\abs{\Gamma\cap B_t}}{\vol B_t}= 
\frac{1}{\vol(G/\Gamma)}+
O_\vre\left(\exp -t\left(\frac{n-1}{n+1}-\vre\right)\right) \,\,.$$


The method developed 
in \cite{LP} uses detailed estimates on solutions to 
the wave equation, and in \cite{Se} the 
method uses refined properties of the spectral 
expansion associated with the Harish 
Chandra spherical transform. In particular these methods  
assume that $G_t$ are 
bi-$K$-invariant sets. 

On the other hand, the estimate of Theorem \ref{th:Lattice points}
holds for {\it any} family of admissible sets $G_t$.   
Thus the following sample corollary 
seems to be new, even in the classical case of 
$G=PSL_2(\RR)$ (or $G=PSL_2(\CC)$). Define
for $1 \le r< \infty$,  $\norm{A}_r=\left(\sum_{i,j=1}^2 
a_{i,j}^r\right)^{1/r}$, and $\norm{A}_\infty=\max \abs{a_{i,j}}$. 

\begin{corollary}\label{closed geod}  For 
any tempered finite-covolume Fuchsian group and for any $1 \le r\le \infty$, 
with the normalization $\vol(G/\Gamma)=1$  
$$\frac{\abs{\set{\gamma\in \Gamma\,;\, 
 \norm{\gamma}_r\le T}}}{\vol \set{g \in SL_2(\RR)\,;\, 
\norm{g}_r\le T}}=
1+O_{\vre,r}\left(T^{-1/4}\right)\,.$$
In particular, this holds for $\Gamma=PSL_2(\ZZ)$.  
\end{corollary} 

%




\subsection{Counting integral unimodular matrices}

Let $G=SL_n(\RR)$ $n \ge 2$ be the group of unimodular 
matrices, and 
$\Gamma=SL_n(\ZZ)$ the group of integral matrices.
A natural choice of balls here are those  
defined by taking the defining representation and 
the rotation-invariant linear norm on $M_n(\RR)$ given by 
$(\tr A^tA)^{1/2}$. Let $B^\prime_T$ denote the norm 
ball of radius $T$ 
intersected with $SL_n(\RR)$. Here the  
best result to date is due to \cite{DRS} and is given by
$$\frac{\abs{\Gamma\cap B^\prime_T}}{\vol B^\prime_T}
=1+O_\vre\left( T^{- \frac{1}{n+1}+\vre}\right)$$

Letting $t=\log T$, the family $B_t=B^\prime_{e^t}$ is 
admissible. 
For our estimate, we need to bound $\theta$, the rate 
of decay of 
$\norm{\pi_{G/\Gamma}(\beta_t)}$ in $L^2_0(G/\Gamma)$.  
For $n=2$, $\theta=1/2+\vre$ 
as noted above, since the representation 
is tempered. For 
$n\ge 3$, $SL_n(\RR)$ has property $T$, and we can simply 
use a bound valid for all of its represenetations simultaneously 
(provided only that 
they contain no invariant unit vectors). Note that  
in the case of $L^2_0(SL_n(\RR)/SL_n(\ZZ))$ this also happens 
to be the best possible estimate, since the spherical 
function with 
slowest decay does in fact occur in the spectrum. 
According to 
\cite{DRS}, every non-constant 
spherical function on $SL_n(\RR)$ is 
in $L^p$ for $p > 2(n-1)$. This implies (see Theorem 
\ref{Kfinite estimate} below) that the matrix 
coefficients of $\pi$ have an estimate in terms of 
$\Xi_G^{1/(n-1)+\vre}$, where $\Xi_G$ is the Harish 
Chandra function. Using the standard estimate 
for $\Xi_G$ (see. e.g. \cite{GV} and also Remark \ref{KS}) 
$$\norm{\pi(\beta_t)}\le 
\left(C_0\vol(B_t)^{-1/2+\vre_0}\right)^{1/(n-1)}\le 
C \exp\left(-t\left(\frac{n^2-n}{2(n-1)}-\vre\right)\right)$$ 
where the last estimate uses the fact that 
(see \cite{DRS})
$$\vol(B^\prime_T)=\vol\set{g\in 
SL_n(\RR)\,;\, \norm{g}_2 \le T}
\cong c_n T^{n^2-n}\,.$$ 

Therefore we have the estimate $\theta =n/2-\vre $, 
so that 
$\theta/(\dim G+1)=(1-\vre)/(2n)$ and 
$\theta/(\dim G/K +1)=(1-\vre)/(n+1)$. Thus 
we recapture the bound given by \cite{DRS},   
for the case of balls defined by the Euclidean norm 
$(\tr A^tA)^{1/2}$. This bound holds whenever the balls 
are bi-$K$-invariant. More generally, letting $n_e$ denote the
 least even integer greater than $n-1$, we have 

\begin{corollary}\label{admiss-SLn}
For any family of admissible 
sets $B_t\subset SL_n(\RR)$, and in particular those 
defined by any norm on $M_n(\RR)$, and for any lattice 
subgroup $\Gamma$, the following bound holds:
$$\frac{\abs{\Gamma\cap B_t}}{\vol B_t}=1+
O_\vre\left(vol(B_t)^{-1/(2 n^2 n_e)+\vre}\right)
$$
\end{corollary}

We note that the method of \cite{DRS} utilizes
the commutativity of the algebra of bi-$K$-invariant 
measures on $G$. Extending this method beyond  
the case of bi-$K$-invariant sets is in principle possible but 
would require further 
elaboration regarding the spectral analysis of $K$-finite functions.

Recently, F. Maucourant \cite{Ma} has obtained a bound 
for the lattice point counting problem for certain
simple groups and certain norms, subject to some 
constraints. Thus for the standard 
representation of $SL_n(\RR)$, when   
$n \ge 7$ the error estimate obtained in \cite{Ma} 
is $1/(6n)+\vre$ which is weaker than the estimate 
above. That is the case 
also for $ 3 \le n \le 6$. The case $n=2$ is not 
addressed in \cite{Ma}.

\subsection{Integral equivalence 
of  general $n$-forms}

\subsubsection{Binary forms} 
 Let us revisit the problem of integral 
equivalence of binary forms considered in \cite{DRS}.  
Let $W_n$ denote the vector space of binary forms of 
degree $n \ge 3$ 
$$W_n=\set{f(x,y)=a_0x^n+a_1x^{n-1}y+\cdots + 
a_ny^n}\,.$$
$SL_2(\RR)$ acts on $W_n(\RR)$ by acting linearly on  
the variables of the form, and when $n\ge 3$ the 
stability group of a generic form is finite. 
Two forms are in the same $SL_2(\RR)$-orbit iff 
they are equivalent under a linear substitution, 
and two forms are in the same $SL_2(\ZZ)$-orbit 
iff they are integrally equivalent. Fix {\it any} norm on $W_n(\RR)$, one example 
being the norm considered in \cite{DRS} 
$$\norm{f}^2=\norm{(a_0,\dots,a_0)}^2=\sum_{i=0}^n 
\binom{n}{i}^{-1} a_i^2$$
The orbits of $SL_2(\RR)$ are closed, and for each 
orbit we can consider the lattice point counting 
problem, or equivalently, 
the problem of counting forms integrally equivalent 
to a given form. Thus fix
 some $f_0$ with finite stabilizer and non-zero discriminant, 
denote   
$B^\prime_T=\set{f\,;\, f\cong_\RR f_0, \norm{f}\le T}$, 
and note that it has been 
established in \cite{DRS} that when the level sets of the form $f_0$ are compact, 
$\vol(B^\prime_T)\sim cT^{2/n}$. We further assume that the form 
satisfies $f_0(x,y) \neq 0$ for $(x,y)\neq (0,0)$. 
We then have the following corollary of Theorem 
\ref{th:Lattice points}. 
\begin{corollary}
Notation being as above, 
 the number of form integrally equivalent with $f_0$ of norm at most $T$ is estimated by 
$$\abs{\frac{\abs{\set{f\,;\,f\cong_\ZZ f_0\,,\, \norm{f}\le T}}}{\vol(B^\prime_T)
}-\frac{1}{\abs{St_{SL_2(\ZZ)}(f_0)}}}$$
$$
\le C(\vre, n, f_0)\vol(B^\prime_T)^{-1/8+\vre}
\le C^\prime(\vre,n,f_0)T^{-1/(4n)+\vre}\,.$$
\end{corollary}

Indeed, the problem under consideration 
is simply that of counting the points   
$\norm{\tau_n(\gamma)f_0}\le T $ where  $\gamma\in SL_2(\ZZ)$, 
for a particular choice of finite 
dimensional representation 
$\tau_n$ of $SL_2(\RR)$, 
and a particular choice of norm on 
the representation space. 
Thus the corollary is an immediate consequence of the
 fact that the sets $B_t=B^\prime_{e^t}$  are admissible (see the Appendix, \S 8.4)
together with Corollary \ref{admiss-SLn} and the fact that the representation of $SL_2(\RR)$ on $SL_2(\RR)/SL_2(\ZZ)$ is tempered.

We note that the existence of the limit was established 
in \cite[Thm. 1.9]{DRS}. The method of proof employed  
there can in principle also be made effective and produce some error estimate. 

\subsubsection{Integral equivalence of forms in many variables}
 Our consideration are not limited to binary 
forms, and we can consider the problem of  
integral equivalence, as well as simultaneous integral equivalence, of $n$-forms in any number of variables. 
Thus let $W_{n,k}$ be the real 
vector space of all 
degree $n$ forms in $k$ variables. $SL_k(\RR)$
 admits a representation $\sigma_{n,k}$ 
on $W_{n,k}$, by acting linearly on the variables.
Fix any norm on $W_{n,k}$. As before, $k_e$ denotes the least even integer greater than $k-1$. 

Consider a form $f_0$ with compact stability group.  
Let us assume that $f_0(x)\neq 0$ for $x\neq 0$, so that the projection of the vectors $u f_0$ onto the highest weight subspace of $W_{n,k}$ never vanishes, as $u$ ranges over a fixed maximal compact subgroup. Let $B^\prime_T$ 
 denote the set of forms integrally equivalent to $f_0$ and of norm at most $T$. Then 
 $B_t=B_{e^t}^\prime$ is an admissible family (see the Appendix \S 8.4, where a volume asymptotic is also established). Hence  Corollary \ref{admiss-SLn}
  applies and yields the following 
  \begin{corollary}{\bf Integral equivalence of forms in many variables}.
  Notation and Assumptions being as in the preceding paragraph, we have 
$$\abs{\frac{\abs{\set{f\,;\,f\cong_\ZZ f_0\,,\, \norm{f}\le T}}}{\vol(B^\prime_T)
}-\frac{1}{\abs{St_{SL_n(\ZZ)}(f_0)}}}$$
$$
\le C(\vre, n, k, f_0)\vol(B^\prime_T)^{-1/(2k^2k_e)+\vre}\,.$$
\end{corollary}

Now let  
$f_1, \dots ,f_N$ be a fixed (but arbitrary) ordered 
 basis of $W_{n,k}$ 
($N=\dim W_{n,k}$). We can consider 
ordered bases $f_1^\prime,\dots,f_N^\prime$ which are 
integrally equivalent to it, namely 
$f_i\cong_\ZZ f_i^\prime $, $1 \le i\le N$. 
Let $B^\prime_T=\set{g\in SL_k(\RR)\,;\, 
\norm{gf_i^\prime}\le T\,,\, 1\le i\le N}$. 
Then $B_t=B^\prime_{e^t}$ is in fact a family of norm balls and any such family is admissible (see the Appendix) so we have, again from Corollary  \ref{admiss-SLn}, the following  
\begin{corollary}{\bf Simultaneous integral equivalence.}
Notation and assumption being as in the preceding paragraph, we have 
$$\abs{
\frac{\abs{\set{(f_1^\prime,\dots,f_N^\prime)\,;\,
 f_i^\prime\cong_\ZZ f_i\,,\,
 \norm{f_i^\prime}\le T\,,\, 
1 \le i\le N }}}
{\vol B^\prime_T }
-1}$$
$$
\le C_\vre \vol(B^\prime_T)^{-1/(2k^2 k_e)+\vre}\,.$$
\end{corollary}


\subsection{Lattice points in $S$-algebraic groups}
All of our results will in fact be formulated 
and proved 
in the context of $S$-algebraic groups. Let us demonstrate 
them in the following simple case, as a motivation for the
developments below. 

Let $p$ be a prime, and consider 
$G_n=PSL_n(\RR)\times PSL_n(\QQ_p)$ 
and the $S$-arithmetic 
lattice $\Gamma_n=PSL_n(\ZZ[\frac1p]))$. Take the norm 
on $M_n(\RR)$ whose square is 
$\tr A^t A$, and its (well-defined) restriction to  
$PSL_n(\RR)$. For $A\in PSL_n(\QQ_p)$ 
let $\abs{A}_p=\max_{1\le i,j\le n} \abs{a_{i,j}}_p$, 
where $\abs{a}_p$ is the $p$-adic absolute value of $a\in\QQ_p$,
 normalized as usual by $\abs{p}_p=\frac1p$. 
If $A\in M_n(\ZZ)$, we write $(A,p)=1$ if $(a_{i,j},p)=1$ for 
some entry $a_{i,j}$. Define the height function on $G_n$ by $H(A,B)=\norm{A}\abs{B}_p$.

 Let $C_T$ be 
of integral matrices with Euclidean norm 
bounded by $T$, and with 
 $\det A$ a power of $p^n$ and $(A,p)=1$, namely 
 $$C_T=\set{A\in 
M_n(\ZZ)\,;\, \tr A^t A\le T^2,\det A\in p^{n\NN}, (A,p)=1}\,\,.$$ 
 
 \begin{proposition}
 The family $C_T$ satisfies 
$$\frac{\abs{C_T}}{\vol\set{g\in G\,;\, H(g)\le T}}
=1+O_{\vre,p,n}\left(T^{-\frac{1}{2n}+\vre}\right)$$
\end{proposition}
The proposition is a consequence 
of Corollary \ref{Lie-error} and the fact that 
the set $C_T$ 
in question is in one-to-one correspondence 
with set of lattice points in balls 
$B_t$ ($t=\log T$) in $G$ defined by the natural 
height function. Indeed,  for $y=(u,v)\in G$ the height is 
$H(y)=\sqrt{\tr u^t u }\cdot\abs{v}_p$. Clearly if $u\in 
PSL_n(\ZZ[\frac1p])$ and $\abs{u}_p=p^k$ (where $k\ge 0$) 
then $A=p^k u\in M_n(\ZZ)$, $(A,p)=1$, and 
$\det A=p^{kn}\det u\in p^{\NN}$. Also  
$\norm{A}=\norm{p^k u }=p^k \norm{u}=H(\gamma)$ where 
$\gamma=(u,u)\in \Gamma$, 
so that $C_T$ maps bijectively with 
$\set{\gamma\in \Gamma\,;\, H(\gamma) \le T}
=\Gamma\cap B_T^\prime$, 
where $B_T^\prime=\set{y\in G\,;\, H(y)\le T}$. 

Now consider the basis of open sets at the identity in $G$ given by 
$\mathcal{O}_\vre=\mathcal{U}_\vre \times\mathcal{K}_p $, 
the product of Riemannian balls $\mathcal{U_\vre}$ 
on $PSL_n(\RR)$ and the 
compact open neighbourhood
$$\mathcal{K}_p =
\set{v\in PSL_n(\QQ_p)\,;\, 
\abs{v-I}_p\le 1\, }.$$
Defining $B_t=B^\prime_{e^t}$, the family $B_t$ is 
admissible w.r.t. to $\mathcal{O}_\vre$. This follows from 
Theorem \ref{Appendix:vol}(4), since the height is defined by a product of two norms. 
The unitary representation of $G_n$ on 
$L^2_0(G_n/\Gamma_n)$ is strongly 
$L^{2(n-1)+\vre}$, and hence (since $\beta_t $ are radial) $\norm{\pi_0(\beta_t)}\le 
\vol(B_t)^{-1/(2(n-1))+\vre}$ (see Remark \ref{KS}). 
A direct calculation 
of the volume of $B_T$ shows that $\vol(B_T^\prime)\le C_\vre  T^{n^2-n+\vre}$, and this gives the error 
term above.

\subsection{Examples of ergodic theorems for lattice actions}

\subsubsection{Exponentially fast convergence on the $n$-torus}

Fix a norm on $M_n(\RR^n)$, 
and consider the corresponding 
norm-balls $G_t\subset SL_n(\RR)$, and the 
averages $\lambda_t$ on $SL_n(\ZZ)\cap G_t$ 

The following result is a direct corollary of Theorem 
\ref{th:Lattice subgroups}, and 
the well-known fact that the action 
of $SL_n(\ZZ)$ on $\TT^n$ admits a spectral gap. 
\begin{corollary}\label{torus action}
 Consider the action of $SL_n(\ZZ)$ of 
$(\TT^n,m)$, where $m$ is Lebesgue measure. 
The averages $\lambda_t$ 
satisfy for every $f\in L^p(X)$, $1 < p < \infty$ for almost every $x\in X$  
$$\abs{\lambda_t f(x)-\int_{\TT^n} fdm}\le 
C_p(f,x) e^{-\eta_n t}$$
where $\eta_n > 0$ is explicit.

\end{corollary}


\subsubsection{Exponentially
 fast convergence in the space of unimodular lattices}
 Let $\Gamma$ be a lattice in a 
simple group $H$ and $\tau : H\to SL_n(\RR)$ 
a rational representation with finite kernel. 
Then the averages 
$\lambda^H_t$ on $\tau(H)\cap B_t$ 
($B_t$ defined w.r.t. a norm on $M_n(\RR)$) converges 
exponentially fast to the ergodic mean, in any of the actions
of $\Gamma$ of $SL_n(\RR)/\Delta$, $\Delta$ a lattice subgroup. 
In particular letting $\Delta=SL_n(\ZZ)$, the homogeneous space 
$\mathcal{L}_n= 
SL_n(\RR)/SL_n(\ZZ)$ can be identified with the space of unimodular 
lattices in $\RR^n$. For such a lattice $L\in \mathcal{L}_n$ 
let $f(L)$ be the number of vectors in $L$ whose 
length (w.r.t. the standard Euclidean norm) is at most one. Then 
for $n\ge 2$, $f\in L^p(\mathcal{L}_n)$, $1 \le p < n$ and 
we let $\kappa_n=\int_{\mathcal{L}_n}f(L)dm(L)$ denote the average 
number of vectors of length at most one in a unimodular 
lattice $L$. Note that by Siegel's formula $\kappa_n$ equals 
the volume of the unit ball in $\RR^n$.

We can now appeal to Theorem
\ref{mean} and Theorem \ref{th:Lattice subgroups} and apply 
them to the averages $\lambda_t^H$. We conclude
\begin{corollary}Let $n \ge 2$ and $1 < p< n$. 
Then for almost every unimodular 
lattice $L\in \mathcal{L}_n$, we have 
$$\frac{\#\set{\gamma\in \Gamma\cap H_t\,;\,
 \abs{\#(\gamma L\cap B_1(0))-\kappa_n}\ge \delta}}
{\#\set{\gamma\in \Gamma\cap H_t}}\le C_p\delta^{-p}
\norm{f}^p_{L^p(\mathcal{L}_n)}e^{-\zeta_{p,n} t}
$$
where $\zeta_{p,n} > 0$ is explicit and depends on 
the spectral gap of the $H$-action on $L^2(H/\Gamma\times \mathcal{L}_n)$
and the admissible family $H_t$. 

\end{corollary}

\subsubsection{Equidistribution and exponentially fast convergence}

Let us consider now the case where 
the lattice $\Gamma$ acts isometrically 
on a compact metric space, preserving a ergodic probability measure 
of full support. Two important families of examples 
are given by 

1) The action of $\Gamma$ on any of its profinite completions, 
with the invariant probability measure being Haar measure on the 
compact group. In particular, this includes the congruence 
completion when $\Gamma$ is arithmetic. 

2) The action of $\Gamma$ on the unit sphere in $\CC^n$ 
or $\RR^n$, via a finite-dimensional unitary or orthogonal 
representation with a dense orbit on the unit sphere (when such 
exist). 

We note that combining 
Theorem \ref{th:equidistribution} and Theorem 
\ref{th:Lattice subgroups}, the following 
interesting phenomenon emerges. 
\begin{corollary}
Let $\Gamma$ be a lattice subgroup in a connected almost 
simple non-compact Lie group with proeprty $T$. Let $G_t$ be 
admissible and $\lambda_t$ the averages uniformly distributed 
on $G_t\cap \Gamma$. 
Then in every isometric action of $\Gamma$ on 
a compact metric space $S$, 
ergodic with respect to a probability 
measure $m$ of full support, the following holds. For every 
continuous function $f\in C(S)$, $\lambda_tf(s)$ converges 
 to $\int_S fdm$ for {\it every} $s\in S$, and converges 
exponentially fast to $\int_S fdm$ 
for {\it almost every} $s\in S$. 
The exponential rate of convergence depends only on $G_t$
 and $G$, and is independent of
 $S$ and $\Gamma$.
\end{corollary}

\subsubsection{Ergodic theorems for free groups} 

Let us note some further ergodic theorems which follow from 
Theorem \ref{th:Lattice subgroups}. 
\begin{enumerate}
\item The index $6$ principal level $2$ 
congruence group $\Gamma(2)$ of $SL_2(\ZZ)$ 
is a free group on two generators. 
Theorem \ref{th:Lattice subgroups} thus gives new ergodic  
theorems for {\it arbitrary actions} of free groups, 
where the averages are taken are uniformly distributed on say 
norm balls. If the free group action has a spectral gap, the 
convergence is exponentially fast. 
These averages are completely different than the averages 
w.r.t. a word metric on the free group discussed in 
\cite{N0}\cite{NS}. 
\item Note that for the averages just described, the phenomenon 
of periodicity (see \cite[\S 10.5]{N5}) associated with the 
existence of the sign character 
of the free group does not arise :
the limit is always the ergodic mean. 

Thus in particular Theorem 
\ref{th:finite} implies that for any norm on $M_2(\RR)$, 
norm balls become equidistributed 
among the cosets of any finite index subgroup of $
\Gamma(2)\cong \FF_2$, at an exponentially fast rate. 

\item Similar comments also apply for example 
to the lattice $\Gamma=PSL_2(\ZZ)
=\ZZ_2\ast \ZZ_3\subset PSL_2(\RR)$ itself, and again the averages in question 
are different from the word-metric ones discussed in \cite{N0}.  
Another family of examples are lattices in $PGL_3(\QQ_p)$, 
to which our results stated in \S 4 apply. In particular 
this includes the lattices acting simply 
transitively on the vertices of the Bruhat-Tits building,  
generalizing \cite[Thm. 11.10]{N5} for these lattices. 

\end{enumerate}

\section{Definitions, preliminaries, and basic tools}

\subsection{Maximal and exponential-maximal inequalities}

Let $G$ be a locally compact second countable (lcsc) group, with a 
left-invariant Haar measure $m_G$. Let $(X,\mathcal{B},\mu)$ be  
a standard Borel space with a Borel measurable $G$-action 
preserving the probability measure $\mu$.
There is a natural isometric representation $\pi_X$ of 
$G$ on the spaces $L^p(\mu)$, $1\le p\le \infty$, defined by
$$
(\pi_X(g)f)(x)=f(g^{-1}x),\quad g\in G,\; f\in L^p(\mu).
$$
To each finite Borel measure $\beta$ on $G$, we associate 
the bounded linear operator
$$
(\pi_X(\beta)f)(x)=\int_G f(g^{-1}x)\,d\beta(g)
$$
acting on $L^p(\mu)$. In particular, 
given an increasing sequence $G_t$, $t>0$, 
of Borel subsets of positive finite measure of $G$,
we consider the Borel probability  measures
\begin{equation}\label{eq:cesaro}
\beta_t=\frac{1}{m_G(G_t)}\int_{G_t} \delta_g\, dm_G(g),
\end{equation}
and the operators $\pi_X(\beta_t)$ are the Haar-uniform 
averages over the sets $G_t$.

\begin{definition}{\bf Maximal inequalities and ergodic theorems.}
{\rm
Let $\nu_t$, $t>0$ be a one-parameter family of absolutely continuous probability
 measures on $G$ such that the map
$t\mapsto \nu_t$ is continuous in the $L^1(G)$-norm.
The maximal function $\sup_{t >  t_0} \abs{\pi_X(\nu_t)f}$, $f\in L^\infty(X)$ is then measurable. We  
define : 
\begin{enumerate}
\item The family $\nu_t$ satisfies the 
{\it strong maximal inequality in 
$(L^p(\mu),L^r(\mu))$,} $p\ge r$, if there exist $t_0\ge 0$ and $C_{p,r}>0$
such that for every $f\in L^p(\mu)$,
$$
\left\|\sup_{t >  t_0} |\pi_X(\nu_t)f|
\right\|_{L^r(\mu)}\le C_{p,r}\|f\|_{L^p(\mu)}.
$$

\item The family $\nu_t$ satisfies the 
{\it mean ergodic theorem in $L^p(\mu)$} if for every $f\in L^p(\mu)$,
$$
\left\|\pi_X(\nu_t)f-\int_X f\,d\mu\right\|_{L^p(\mu)}\to 0\quad \hbox{as}\quad t\to\infty.
$$
\item The family $\nu_t$ satisfies the {\it pointwise ergodic theorem in $L^p(\mu)$} if for every $f\in L^p(\mu)$,
$$
\pi_X(\nu_t)f(x)\to\int_X f\,d\mu\quad \hbox{as}\quad t\to\infty
$$
for $\mu$-almost every $x\in X$.
\item The family $\nu_t$ satisfies the {\it exponentially fast
 mean ergodic theorem in $(L^p(\mu),L^r(\mu))$}, $p\ge r$, if 
there exist $C_{p,r}>0$ and 
$\theta_{p,r}>0$ such that for every $f\in L^p(\mu)$,
$$
\left\|\pi_X(\nu_t)f-\int_X f\,d\mu\right\|_{L^r(\mu)}\le C_{p,r} 
e^{-t\theta_{p,r} }\|f\|_{L^p(\mu)}.
$$
\item The family $\nu_t$ satisfies {\it exponential strong 
maximal inequality in $(L^p(\mu),L^r(\mu))$},
 $p\ge r$, if there exist $t_0\ge 0$, $C_{p,r}>0$, and 
$\theta_{p,r}>0$ such that for every $f\in L^p(\mu)$,
$$
\left\|\sup_{t\ge t_0} e^{t\theta_{p,r} }
 \left|\pi_X(\nu_t)f-\int_X f\,d\mu\right|\right\|_{L^r(\mu)}
\le C_{p,r}\|f\|_{L^p(\mu)}.
$$
\item The family $\nu_t$ satisfies {\it exponentially fast  
pointwise ergodic theorem in $(L^p(\mu),L^r(\mu))$}, 
$p\ge r$, if there exist $t_0\ge 0$, and 
$\theta_{p,r}>0$ such that for every $f\in L^p(\mu)$,
$$
\left|\pi_X(\nu_t)f(x)-\int_X f\,d\mu\right|\le B_{p,r}(x,f)
e^{-t\theta_{p,r} }
\quad\hbox{ for $\mu$-a.-e. $x\in X$}
$$
with the estimator $B_{p,r}(x,f)$ satisfying the norm estimate 
$$
\|B_{p,r}(\cdot,f)\|_{L^r(\mu)}\le C_{p,r}\|f\|_{L^p(\mu)}.
$$

\end{enumerate} 
}
\end{definition}

\begin{remark}{\rm
The main motivation to consider the 
exponential strong maximal inequality in 
$(L^p(\mu),L^r(\mu))$ is that it implies
the exponentially fast pointwise ergodic theorem 
in $(L^p(\mu),L^r(\mu))$, together with norm convergence to the ergodic mean, 
at an exponential rate. 

We recall, in comparison, that the ordinary 
strong maximal inequality only implies pointwise convergence almost 
surely provided that  
we establish also the existence of a dense subspace where almost sure 
pointwise convergence holds. In addition, convergence in norm requires 
a separate further argument.}
\end{remark}

\begin{remark}{\rm 
If the mean ergodic theorem holds in $L^p(\mu)$, 
using appoximation by bounded functions and H\"older inequality,
one can deduce the mean ergodic theorem in $L^{p'}(\mu)$ 
for $1\le p'\le p$. 
Similarly, the strong maximal inequality (resp. 
the exponentially fast mean ergodic theorem, 
the exponential strong maximal inequality) in 
$(L^p(\mu),L^r(\mu))$ implies
the strong maximal inequality 
(resp. exponentially fast mean ergodic theorem, 
the exponential strong maximal inequality)
 in $(L^{p'}(\mu),L^{r'}(\mu))$ 
for $p'\ge p$ and $1\le r'\le r$.}

\end{remark}

\subsection{$S$-algebraic groups and upper local dimension}

We now define the class of $S$-algebraic 
groups which will be our main focus.

\begin{definition}\label{alg-gps}
{\bf $S$-algebraic groups.}
\begin{enumerate}

\item 
Let $F$ be a locally compact non-discrete field, 
and let $G$ be the group  
of $F$-points of a semisimple linear algebraic group defined 
over $F$, 
with positive $F$-rank (namely containing an $F$-split torus of 
positive dimension over $F$). We assume in addition 
that $G$ is algebraically connected, and does not have non-trivial 
anisotropic (i.e. compact) algebraic 
factor groups defined over $F$. We will also assume, for simplicity, that $G^+$ is of finite index in $G$ 
(see Remark \ref{G+}). 

\item By an $S$-algebraic group we mean any finite 
product of the groups 
described in (1). 
 \end{enumerate}
\end{definition}

The unitary representation theory of $S$-algebraic groups 
has a number of useful features which we will use 
extensively below. 
Another property of $S$-algebraic groups 
which is crucial for handling their lattice points is the finiteness of 
their upper local dimension, as defined by natural choices 
of neighborhood bases. Let us introduce the following 

\begin{definition}\label{ULD}
 For a family of neighborhoods 
$\{\mathcal{O}_\vre\}_{0<\vre<1}$ of $e$ in an lcsc group $G$ such that 
$\mathcal{O}_\vre$'s are 
symmetric, bounded, and increasing with $\vre$, we let 
\begin{equation}\label{eq:B_epoint3}
\varrho_0\stackrel{def}{=}
\limsup_{\vre\to 0^+}
 \frac{\log m_G(\mathcal{O}_\vre)}{\log\vre}<\infty.
\end{equation}
\end{definition}
\begin{remark}
\begin{enumerate}
\item When $M$ is a Riemannian manifold and 
$\mathcal{O}_\vre$ are the balls w.r.t. the Riemannian metric, 
the condition $m_M(\mathcal{O}_\vre)\ge C_\rho\vre^{\rho}$, 
$\vre > 0$  is 
equivalent to $\dim (M)\le \rho$.
\item When $G$ is an $S$-algebraic group, we will always take 
$\mathcal{O}_\vre$ to be the sets $\mathcal{U}_\vre\times K_0$, where $\mathcal{U}_\vre$ 
is the family of Riemannian balls in the Archimedean component of $G$ (if it exists), and $K_0$ a fixed compact open subgroup of the totally disconnected component of $G$. Thus the local dimension  
of $\mathcal{O}_\vre$ is the dimension of the Archimedean component. 
\end{enumerate}
\end{remark}

\subsection{Admissible and coarsely admissible sets}

We begin our discussion of admissibility 
by introducing a coarse version of it, which will be useful in what follows.   

\begin{definition}{\bf Coarse admissibility.} {\rm
Let $G$ be an lcsc group with left Haar measure $m_G$. 
An increasing family of bounded Borel subsets $G_t$ ($t\in \RR_+$ or $t\in \NN_+$) 
of $G$ will be called
{\it coarsely admissible} if 
\begin{itemize}
\item For every 
bounded $B\subset G$, there exists $c=c_B>0$ such that for all 
sufficiently large $t$,
\begin{equation}\label{eq:B_max}
B\cdot G_t\cdot B \subset G_{t+c}.
\end{equation}

\item For every $c>0$, there exists $d>0$ such that for all 
sufficiently large $t$,
\begin{equation}\label{eq:B_max2}
m_G(G_{t+c})\le d\cdot  m_G(G_{t}).
\end{equation}
\end{itemize}
}
\end{definition}

It will be important in our considerations later on 
 that coarse admissibility 
implies at least a certain minimal 
amount of volume growth for our family $G_t$, provided that the group is compactly generated. 
This property will play a role  
in the spectral estimates that will arise in the proofs of Theorem
\ref{th:semisimple-1} and Theorem \ref{th:semisimple-T}. 
Thus let us note the following. 
\begin{proposition}\label{growth}
{\bf Coarse admissibility implies growth.}
When $G$ is compactly generated, coarse admissiblity for 
an increasing family of 
bounded Borel subset $G_t$, $t>0$, of $G$ implies that 
 for any bounded  symmetric 
generating set $S$ of $G$, there exist 
$a=a(S) > 0$, $b=b(S)\ge 0$ such that $S^n\subset G_{an+b}$. 
\end{proposition}
\begin{proof}
Let $S$ be a compact symmetric generating set. Taking $B$ to be a 
bounded open set containing the identity together with 
$G_{t_0}\cup G_{t_0}^{-1}$, 
and applying condition (\ref{eq:B_max}) we 
conclude that $G_{t_0+c}$ contains an open 
neighborhood of the identity. 
Then, assuming without loss of generality 
 that $e\in S$ we have 
$S\subset SG_{t_0+c}S \subset G_{t_1}$. Applying condition 
(\ref{eq:B_max}) repeatedly, we conclude that  
$S^n \subset G_{t_1+nc_1}$ and the required property follows. 
\end{proof}

\begin{remark}{\bf Sequences in totally disconnnected groups} 
If $G$ is totally disconnected, 
and $K\subset G$ is a compact open subgroup,
 then $G/K$ is a discrete countable 
metric space. If $G_t\subset G$, 
$t\in \RR_+$ is an increasing family 
 of bounded sets, then their projections to 
$G/K$ will yield only a sequence of distinct sets. 
Since it is the large scale behaviour of the sets 
that we are mostly interested in, it is natural 
to assume that in the totally 
disconnected case the family $G_t$ is in fact 
countable, and we then parametrize it by 
$G_t$, $t\in \NN_+$. 
This convention will greatly simplify our notation below. 
\end{remark} 

We now consider the following 
abstract notion of admissible families, which (as we shall see) 
generalizes the one 
introduced in \S 1.  
\begin{definition}\label{admissible}{\bf Admissible families.} {\rm
\begin{enumerate}
\item {\it Admissible $1$-parameter families}. 
Let $G$ be an lcsc group, 
fix a family of neighborhoods 
$\{\mathcal{O}_\vre\}_{0<\vre<1}$ of $e$ in $G$ such that 
$\mathcal{O}_\vre$'s are symmetric, 
bounded, and decreasing with $\vre$.

 An increasing $1$-parameter 
family of bounded Borel subset $G_t$, $t\in \RR_+$, 
on an lcsc group 
 $G$ will be called 
{\it admissible} (w.r.t. to the family $\mathcal{O}_\vre$)
 if it is coarsely admissible and   
there exist $c>0$, $t_0 > 0$ and $\vre_0 > 0$ such that for 
 $t\ge t_0$ and  $0 < \vre\le \vre_0$
\begin{align}
\mathcal{O}_\vre
\cdot G_t\cdot \mathcal{O}_\vre &\subset G_{t+c\vre},\label{eq:1}\\
m_G(G_{t+\vre})&\le (1+c\vre)\cdot  m_G(G_{t}),\label{eq:2}
\end{align}

\item {\it Admissible sequences}. An increasing sequence 
 bounded Borel subset $G_t$, $t\in \NN_+$, 
on an lcsc totally disconnected group 
 $G$ will be called 
{\it admissible}  if it is coarsely admissible, 
and there exists $t_0 > 0$ and 
a compact open subgroup $K_0$ such that for $t\ge t_0$ 
\begin{align} K_0 G_t K_0 =G_t\,. \label{eq:4}
\end{align} 
\end{enumerate}
}

\end{definition}

Let us note the following regarding admissibility. 
\begin{remark}\label{loc.Lip}
\begin{enumerate}

\item 
When $G$ is connected and 
$\mathcal{O}_\vre$ are Riemannian balls 
every $\mathcal{O}_\vre$ generates $G$, 
and so it 
is clear that admissibility of the $1$-parameter family 
$G_t$ implies 
coarse admissibility (and thus also the minimal growth 
condition). 
However this argument fails for $S$-algebraic groups which have a totally disconnected simple component, and 
so we have required coarse admissibility explicitly in the definition.

\item Condition (\ref{eq:2})
 is of course equivalent to the function 
$\log m_G(G_t)$ being {\it uniformly} locally Lipschitz continuous, 
for sufficiently large $t$. Furthermore, note that 
$$\norm{\beta_{t+\vre}-\beta_t}_{L^1(G)}=\int_G 
\frac{\abs{m_G(G_t)
\chi_{G_{t+\vre}}-m_G(G_{t+\vre})\chi_{G_t}}}{m_G(G_t)\cdot 
m_G(G_{t+\vre})}dm_G =$$
$$=\frac{2(m_G(G_{t+\vre})-m_G(G_t))}{m_G(G_{t+\vre})}\,.$$
It follows that admissibility implies that the map $t\mapsto \beta_t$ is uniformly 
locally Lipschitz continuous as a map 
from $[t_0,\infty)$ to the Banach space $L^1(G)$. 
The converse also holds, provided we assume in addition 
that the ratio of 
$m_G(G_{t+\vre})$ and 
$m_G(G_t)$ is uniformly bounded for $t\ge t_0$ and 
$\vre \le \vre_0$. 
\item We note that we can relax the Lipschitz conditions in the 
definition of admissibility to the corresponding 
H\"older conditions. 
Such averages will be called H\"older-admissible, and will be discussed further below.
%

\end{enumerate}

\end{remark}

\subsection{Absolute 
continuity, and examples of admissible averages}

Admissible $1$-parameter families posses a regularity property  
which will be crucial in the proof of ergodic theorems in the 
absence of a spectral gap, and thus in 
Theorem \ref{th:semisimple-1} and Theorem 
\ref{th:semisimple_lattice-1}.

To define the property, let us 
first note that a $1$-parameter family $G_t$ gives rise to
 the gauge $\abs{\cdot} : G\to \RR_+$ defined by 
$\abs{g}=\inf\set{s> 0\,;\, g\in G_s}$. If the family $G_t$ 
satisfies the condition 
$\cap_{r> t}G_r=G_t$ for every $t\ge t_0$, then conversely   
the family $G_t$ is determined by the 
gauge, namely 
$G_t=\set{g\in G\,;\, \abs{g} \le t}$ for $t\ge t_0$.
Note that clearly admissibility implies that $\cap_{r > t} G_r$ 
can only differ from $G_t$ by a set of measure zero. 
Clearly the resulting family is still admissible, and so 
we can and will assume from now on 
that $G_t$ is indeed determined 
by its gauge.

\begin{proposition}\label{prop:abs.cont}
{\bf Absolute continuity.}
An admissible $1$-parameter family $G_t$ (w.r.t. a basis 
$\mathcal{O}_\vre$, $0 < \vre < \vre_0$) on an lcsc group $G$ 
has the following property.  
The map $g\mapsto \abs{g}$ from $G$ to $\RR_+$ given 
by the associated gauge maps Haar 
measure on $G$ 
to a measure on $[t_0,\infty)$ 
which is absolutely continuous with respect to linear 
Lebesgue measure.
\end{proposition}

\begin{proof} 

 The measure $\eta$ induced on $\RR_+$ by the 
map $g\mapsto \abs{g}$ 
is by definition $\eta(J)=m_G(\set{g\in G\,;\, \abs{g}\in J})$, for 
any Borel set $J\subset \RR_+$. Assume that $J\subset [t_0,t_1)$ and that 
$\ell(J)=0$, namely $J$ has linear 
Lebesgue measure zero, and let us show that $\eta(J)=0$. 
Indeed, for any $\kappa > 0$ there exists a covering of 
$J$ by a sequence of intervals $I_i$, with 
$\sum_{i=1}^\infty \ell(I_i)< \kappa$. Subdividing the intervals 
if necessary, we can assume that $\ell(I_i)< \vare_0$. 
By (\ref{eq:2}), for all $\vare < \vre_0$ and $t\in [t_0,t_1)$
$$\eta((t,t+\vre])=m_G(\set{g\,;\, 
t < \abs{g}\le t+\vre})= m_G(G_{t+\vre})-m_G(G_t)
$$
$$\le  c\vre m_G(G_t)\le 
c m_G(G_{t_1})\ell((t,t+\vre])\,.$$
Denoting $cm_G(G_{t_1})$ by $C$, we see that 
$$\eta(J)\le \sum_{i=1}^\infty\eta(I_i)\le C\sum_{i=1}^\infty
 \ell(I_i)\le C\kappa$$
and since $\kappa$ is arbitrary it follows that $\eta(J)=0$ and thus 
$\eta$ is absolutely continuous w.r.t. $\ell$. 
\end{proof}

Let us now  
verify the first 
assertion made regarding admissibility 
in \S 1. The second assertion is discussed immediately below
and the third is proved in Lemma \ref{stability}.

\begin{proposition}\label{Riemannian}
When $G$ is a connected Lie group 
and $\mathcal{O}_\vre$ are the balls 
defined by a left-invariant 
Riemannian metric, admissibility is independent of the Riemannian 
metric chosen to define it 
(but the constant $c$ may change). 

\end{proposition}

\begin{proof} 
We will verify that (\ref{eq:1})
 is still satisfied, possibly with another 
constant $c$ (but keeping $\vare_0$ and $t_0$ the same) 
if we choose another 
Riemannian metric. 
Fix such a Riemannian metric and denote its balls by 
$\mathcal{O}_\vre^\prime$. 
First note that 
 it suffices to verify (\ref{eq:1}) for 
all $\vre < a$ where $a$ is {\it any} positive constant. Indeed then 
for $t\ge t_0$ and $\vre < a $ 
$$\mathcal{O}^\prime_{2\vre}G_t\mathcal{O}^\prime_{2\vre}= 
\mathcal{O}^\prime_\vre \mathcal{O}^\prime_\vre G_t 
\mathcal{O}^\prime_\vre 
\mathcal{O}^\prime_\vre \subset 
\mathcal{O}^\prime_\vre 
G_{t+c^\prime\vre}\mathcal{O}^\prime_\vre\subset G_{t+2c^\prime\vre}$$
Here we have used the property 
$\left(\mathcal{O}^\prime\right)_\vre^n=\mathcal{O}^\prime_{n\vre}$ 
which is valid for invariant Riemannian metrics. 
It follows that $\mathcal{O}^\prime_\vre G_t\mathcal{O}^\prime
_\vre\subset G_{t+c^\prime \vre}$ 
 holds for all $0 < \vre \le \vre_0$. 

Now note that it is possible to choose $a> 0$ small enough so 
that for $\vre < a$ 
there exists a fixed $m$ independent of $\vre$ such that 
$$\mathcal{O}_\vre ^\prime \subset \mathcal{O}^m_\vre
=\mathcal{O}_{m\vre}\,.$$ 
This fact follows by applying the 
exponential map in a sufficiently small 
ball in the Lie algebra of $G$, and using the fact that 
any two norms on the Lie algebra are equivalent. It then follows 
that for $\vre < a$, $t\ge t_0$ 
$$\mathcal{O}_\vre^\prime G_t \mathcal{O}_\vre^\prime \subset 
\mathcal{O}_\vre^m G_t \mathcal{O}_\vre^m \subset G_ {t+mc\vre}$$
as required, with $c^\prime=mc$.

\end{proof}

Admissible averages exist in abundance on $S$-algebraic groups. We refer 
to \S 8.1 in the Appendix for a  proof of the following result. 

\begin{theorem}\label{Appendix:vol}
For  an $S$-algebraic group $G=G(1)\cdots G(N)$ as in Definition \ref{alg-gps},
the following families of sets $G_t\subset G$ are admissible, where $a_i$ are any positive constants.
\begin{enumerate}
\item Let $S$ consist of infinite places, and let $G(i)$ be a closed
subgroup of the isometry group of a symmetric space $X_i$ of nonpositive curvature
equipped with the Cartan--Killing metric. For $u_i,v_i\in  X_i$, define
$$
G_t=\{(g_1,\ldots,g_N): \sum_i   a_i d_i(u_i, g_i\cdot v_i)<t\}.
$$
\item Let $S$ consist of infinite places, and let $\rho_i:G(i)\to \hbox{\rm GL}(V_i)$
be proper rational representations. For norms $\|\cdot \|_i$ on $\hbox{\rm End}(V_i)$, define
$$
G_t=\{(g_1,\ldots,g_N): \sum_i a_i \log \|\rho_i(g_i)\|_i<t\}.
$$
\item For infinite places, let $X_i$ be the symmetric space of $G(i)$ equipped with the Cartan-Killing
 distance $d_i$, and  for finite places, let $X_i$ be 
the Bruhat-Tits building of $G(i)$ equipped with the path metric $d_i$ on its $1$-skeleton. 
For $u_i\in  X_i$, define
$$
G_t=\{(g_1,\ldots,g_N): \sum_i a_i d_i(u_i, g_i\cdot u_i)<t\}.
$$
\item Let $\rho_i:G(i)\to \hbox{\rm GL}(V_i)$
be proper representations, rational over the fields of definition $F_i$. For infinite places, let
$\|\cdot \|_i$ be a Euclidean norm on $End(V_i)$, 
 and assume 
that $\rho_i(G(i))$ is self-adjoint : $\rho_i(G(i))^t=\rho_i(G(i))$.  
For finite places, let $\|\cdot \|_i$ be the $\max$-norm on $\hbox{\rm End}(V_i)$.
Define
$$
G_t=\{(g_1,\ldots,g_N): \sum_i  a_i \log \|\rho_i(g_i)\|_i<t\}.
$$
\end{enumerate}
\end{theorem}

An important class of families are those defined by height functions on $S$-algebraic groups.  
In Theorem \ref{holder heights} (Appendix, \S 8.5)  will establish the following  

\begin{theorem}\label{HHA}{\bf Heights are H\"older-admissible.} 
For  an $S$-algebraic group $G=G(1)\cdots G(N)$ as in Definition \ref{alg-gps}, 
let $\rho_i:G(i) \to \hbox{\rm GL}(V_i)$
be proper representations, rational over the fields of definition $F_i$. For infinite places, let
$\|\cdot \|_i$ be any Euclidean norm on $\hbox{\rm End}(V_i)$. 
For finite places, let $\|\cdot \|_i$ be the $\max$-norm on $\hbox{\rm End}(V_i)$.
Define (for any positive constants $a_i$) 
$$
G_t=\{(g_1,\ldots,g_N): \sum_i  a_i \log \|\rho_i(g_i)\|_i<t\}.
$$
Then $G_t$ are H\"older-admissible. 
\end{theorem}

\subsection{Balanced and well-balanced  
 families on product groups}
 In any discussion of ergodic theorems for averages $\nu_t$ 
on a product group $G=G_1\times G_2$, it is necessary to discuss the 
behaviour of the two projections $\nu_t^1$ and $\nu_t^2$ to the factor groups. 
Indeed, consider the case where one of these projections, 
say $\nu_t^1$, 
 assigns a fixed fraction of its measure to a bounded set for all  
$t$. Then choosing an ergodic action of $G_1$, we can view it 
as an ergodic action of the product in which $G_2$ acts trivially, 
and it is clear that the ergodic 
theorems will fail for $\nu_t$ in this action. 
Thus it is necessary to require one of
 the following two conditions. Either the projections of the 
averages $\nu_t$ 
to the non-compact factors do not assign a fixed fraction of 
their measure to 
a bounded set, or alternatively 
 that the action is irreducible,
 namely every non-compact factor acts ergodically. 
This unavoidable 
assumption is reflected in the following definitions.

\begin{definition}\label{balanced}{\bf balanced and 
well-balanced averages.}

Let $G$=$H_1 \cdots  H_N $ be an 
almost direct product of $N$ non-compact 
compactly generated 
subgroups.
For a set $I$ of indices $I\subset [1,N]$, 
let $J$ denote its complement, 
and $H_I=\prod_{i\in I}H_i$. 
Let $G_t$ be an increasing family of sets 
contained in $G$. 
\begin{enumerate}
\item
$G_t$ will be called balanced if for every $I$ satisfying $1 < \abs{I}< N$, 
and every compact set $Q$ contained in $H_I$ 

$$\lim_{t\to\infty}\frac{m_G(G_t\cap H_J Q)}{m_G(G_t)}
=0\,\,.$$
\item  An admissible family $G_t$ will be called 
well-balanced if there 
exists $a > 0$ and $\eta > 0$ such that for all $I$ satisfying 
$1< \abs{I}<  N$  
 $$\frac{m_G(G_{an}\cap H_J \cdot S_I^n )}{m_G(G_{an})}
\le Ce^{-\eta n}\,\,$$
where $S_{I}$ is a compact generating set 
of $H_I$, 
consisting of products of compact
 generating sets of its component groups.
\end{enumerate}
\end{definition}

The condition of being well-balanced is independent of the choices 
of compact generating sets in the component groups, but the various constants 
may change. An explicit sufficient condition for a family of sets 
$G_t$ defined by a norm on a semisimple Lie group 
to be well-balanced is given in \S 7.3. In addition, we note the following important natural examples 
of admissible well-balanced families of averages, and state an estimate on their boundary measures
 which will play an important role in the proof of Theorem \ref{th:semisimple-1}. 
 A complete proof of Theorem \ref{CAT} will be given in the Appendix. 

\begin{theorem}\label{CAT}
Let $G=G(1)\cdots G(N)$ be an $S$-algebraic group and
$\ell_i$ denote the standard $CAT(0)$-metric on either
the symmetric space $X_i$ or the Bruhat--Tits building $X_i$
associated to $G(i)$. For $p>1$ and $u_i\in X_i$, define
$$
G_t=\{(g_1,\ldots, g_N):\, \sum_i \ell_i(u_i, g_iu_i)^p<t^p\}.
$$
Let $m$ be a Haar measure $G$.

\begin{enumerate}
\item[(i)] 
There exist $\alpha,\beta>0$ such that for every nontrivial projection $\pi:G\to L$,
$$
m(G_t\cap \pi^{-1}(L_{\alpha t}))\ll e^{-\beta t}\cdot m_t(G_t)\,,
$$
namely the averages are well balanced.
\item[(ii)] If $G$ has at least one Archimedian factor, then the family $G_t$ is admissible, and 
writing $m=\int_0^\infty m_t\, dt$ where $m_t$
is a measure supported on $\partial G_t$, the following estimate holds : 

There exist $\alpha,\beta>0$ such that for every nontrivial projection $\pi:G\to L$,
$$
m_t(\partial G_t\cap \pi^{-1}(L_{\alpha t}))\ll e^{-\beta t}\cdot m_t(\partial G_t)\,,
$$
namely the averages are boundary-regular. 

\end{enumerate}
\end{theorem}

Let us introduce the following definition : 
\begin{definition}\label{standard}{\bf Standard radial averages.}
Let $G$ be an $S$-algebraic group as in Definition \ref{alg-gps}, and represent $G$ as a  
product $G=G_1\cdots G_N$ of its simple components. 
We will refer to any of the families defined in Theorem \ref{Appendix:vol}, Theorem \ref{holder heights} and Theorem \ref{CAT} as standard  radial averages. 

If the family satisfies in addition the estimate in Theorem \ref{CAT} (ii) it will be called boundary-regular. 
\end{definition}

\subsection{Roughly radial and quasi-uniform sets}
We now define several other stability properties for families  
of sets $G_t$ that 
 will be useful in the arguments below.

\begin{definition}{\bf Quasi-uniform families.} {\rm
An increasing $1$-parameter 
   family of bounded Borel subset $G_t$, 
$t>0$, of $G$ will be called
{\it quasi-uniform} if it satisfies the following two conditions.
\begin{itemize}
\item Quasi-uniform local stability. For every $\vre>0$, there exists a 
neighborhood $\mathcal{O}$ of $e$ in $G$ such that for all sufficiently large $t$,
\begin{equation}\label{eq:01}
\mathcal{O}\cdot G_t \subset G_{t+\vre}.
\end{equation}
\item Quasi-uniform continuity. 
For every $\delta>0$, there exist $\vre>0$ such that for all 
sufficiently large $t$,
\begin{equation}\label{eq:02}
m_G(G_{t+\vre})\le (1+\delta)\cdot  m_G(G_{t}),
\end{equation}
\end{itemize}
}
\end{definition}
Note that (\ref{eq:02}) is equivalent to
 the function $\log m_G(G_t)$ being
 quasi-uniformly continuous in $t$, and implies 
that $t\mapsto \beta_t$ is quasi-uniformly 
continuous in the $L^1(G)$-norm. 
If the ratio of 
$m_G(G_{t+\vre})$ and $m_G(G_t)$ is uniformly bounded 
for $ 0< \vre\le \vre_0$ the converse holds as well.    

An important ingredient in our analysis below will be the existence 
of a radial structure on the groups under consideration.
Thus let $G$ be an lcsc group, 
and $K$ a compact subgroup. Sets which are bi-invariant 
under translations by 
$K$ will be used 
in order to dominate
 sets which are not necessarily bi-$K$-invariant. 

In particular, we shall utilize special bi-$K$-invariant sets, 
called ample sets, which play a key role in the ergodic theorems 
proved in \cite{N4}, 
which we will use below. We recall the definitions.

\begin{definition}\label{radial sets}
{\bf Roughly radial sets and ample sets.}
{\rm Let $K\subset G$ be a fixed compact subgroup, $\mathcal{O}$ a fixed neighbourhood of $e\in G$, and $C$, $D$  positive constants. 
\begin{enumerate}
\item $B\subset G $ is called left radial (or more precisely $K$-radial) if it satisfies $KB=B$, where $K$ is of finite 
index in a maximal compact subgroup of $G$. 
\item \cite{N3} A measurable set $B\subset G$ of positive 
finite measure 
will be called roughly radial (or more precisely  
$(K,C)$-radial) 
provided that $m_G(KBK)\le C m_G(B)$.
\item (see \cite{N4}) A measurable 
set $B\subset G$ of positive finite measure is
 called ample (or more precisely 
$(\mathcal{O},D,K)$-ample) if it satisfies 
$m_G(K\mathcal{O}B K)\le D m_G(B)$. 
\end{enumerate}
}
\end{definition}

To illustrate the definition of ampleness, consider first 
the case where
 $G$ is a connected semisimple Lie group.
 We can fix a maximal compact subgroup 
$K$ of $G$, and consider the symmetric space $S=G/K$, 
with the distance $h$ derived from the 
Riemannian metric associated with the Killing form. 
Ampleness can be equivalently defined as follows. 
 For a $K$-invariant set $B\subset G/K$, consider 
 the $r$-neighborhood of $B$ in the symmetric space, given by 
$$U_r(B)=\set{gK\in G/K\,;\, h(gK, B)< r} \,.$$
Then $B$ is $(\mathcal{O}_r,D, K)$-ample 
iff $m_{G/K}(U_r(B))\le Dm_{G/K}(B)$, where 
$\mathcal{O}_r$ is the lift to $G$ of 
a ball of radius $r$ and center $K$ in $G/K$.

The following simple facts are obvious from the definition, but 
since they will be used below we record them for completeness.

\begin{proposition}\label{KC-radial}
 Any family of coarsely admissible sets on an lcsc 
group is $(K,C)$-radial for some finite $C$ and a (good) 
maximal compact subgroup $K$, as well as $(\mathcal{O},D,K)$-ample, 
for some neighborhood $\mathcal{O}$ and $D > 0$. 
\end{proposition}
\begin{proof} By definition of 
coarse admissibility, for the compact set $K$ there exists $c\ge 0$ with 
$K\cdot G_t\cdot K \subset G_{t+c}$. Therefore 
$m_G(KG_tK)\le d m_G(G_t)$, 
so that the family $G_t$ is $(K,d)$-radial.
The fact that $G_t$ are ample sets is proved in the same way. 
\end{proof}

Let us note that when $G$ is totally disconnected, 
and there exists a compact open subgroup $Q\subset 
\mathcal{O}_\vre$  satisfying 
$QG_t Q =G_{t}$ for all $t \ge t_0$, then $G_t$ are $(K,C)$-radial. 
Indeed $Q$ is of finite index in 
a good maximal compact subgroup $K$.  
Denoting the index by $N$, we have 
$$KG_t K=\cup_{i,j=1}^N k_iQ G_t  Qk_j \subset \cup_{i,j=1}^N 
k_i G_{t} k_j\,.$$ 
It follows that 
$m_G(KG_tK)\le N^2 m_G(G_t)$ and $G_t$ is $(K,N^2)$-radial.

\subsection{Spectral gap and strong spectral gap}

We recall the definition of spectral gaps, as follows.

\begin{definition}{\bf Spectral gaps}.
\begin{enumerate} 
\item A strongly continuous  
unitary representation $\pi$ of an lcsc group $G$ is said to have 
a spectral gap if  $\norm{\pi(\mu)}< 1$, for some (or equivalently,
all) absolutely 
continuous symmetric probability measure
$\mu$ whose support generates $G$ as a group. 
\item Equivalently, $\pi$ has a spectral gap if the Hilbert space
does not admit an asymptotically-$G$-invariant sequence
of unit vectors, namely a sequence satisfying 
$\lim_{n\to \infty}\norm{\pi(g)v_n-v_n}=0$
uniformly on compact sets in $G$.
\item A measure preserving action
of $G$ on a $\sigma$-finite measure space $(X,m)$ 
is said to have a spectral gap if the unitary representation $\pi_X^0$ of
$G$ in the space orthogonal to the space of $G$-invariant functions
has a spectral gap. Thus in the case of an ergodic
probability-preserving action, the representation in question is on the space 
$L^2_0(X)$ of function of zero integral. 
\item An lcsc group $G$ is said to have
 {\it property $T$} \cite{Ka}
 provided every strongly
continuous unitary representation which does not have $G$-invariant
unit vectors has a spectral gap. 
\end{enumerate}
\end{definition} 

If $G=G_1G_2$ is a (almost) direct product group, and there does not exists a 
sequence on unit vectors which is asymptotically invariant under every 
$g\in G$, it may still be the case that there exists such a sequence 
asymptotically invariant under the elements of a subgroup of $G$, for example 
$G_1$ or $G_2$. It is thus natural to introduce the following

\begin{definition}{\bf Strong spectral gaps.}
Let $G=G_1\cdots G_N $ be an almost direct product of $N$ 
lcsc subgroups. 
A strongly continuous unitary representation $\pi$ 
of $G$ has a strong spectral gap (w.r.t. the given decomposition) if the 
restriction of $\pi$ to every almost direct factor  
$G_i$ has a spectral gap.    
\end{definition}

\begin{remark}
\begin{enumerate}
\item 
Let $G_1=G_2=SL_2(\RR)$, $G=G_1\times G_2$ and 
$\pi=\pi_1\otimes\pi_2$, 
where $\pi_1$ has a spectral gap and 
$\pi_2$ does not (but has no invariant unit vectors). It is possible to 
construct two admissible families $G_t$ and $G_t^\prime$ on $G$, such that 
$\norm{\pi(\beta_t)}
\le Ce^{-\theta t}$, but $\norm{\pi(\beta_t^\prime)}
\ge b > 0$. For example, $G_t$ can be taken as the 
inverse images of the families of balls 
of radius $t$ in $\HH\times \HH$, w.r.t. the Cartan-Killing metric. 
 For a construction of $G_t^\prime$, one can use (in the obvious way) 
the non-balanced averages constructed in \S 7.2. 
\item 
In order to obtain conclusions which assert that the mean, maximal 
or pointwise theorem hold for $\beta_t$ with an exponential rate, 
it is of course necessary 
that in the representation $\pi_X^0$ in 
$L^2_0(X,\mu)$ 
$\beta_t$ have the exponential decay property, namely 
$\norm{\pi_X^0(\beta_t)}\le C e^{-\theta t}$, $\theta > 0$. 
In Theorem \ref{Exp decay} we will give sufficient conditions for 
the latter property to hold. 
\item Consider the special case $X=G/\Gamma$, where $G$ is a semisimple 
Lie group and $\Gamma$ a lattice subgroup. It is a standard 
corollary of the 
theory of elliptic operators on compact manifolds
 that if $\Gamma$ is co-compact, then the (positive) 
Laplacian $\Delta$ on $G/\Gamma$ has a spectral gap above zero, namely 
$\norm{\exp(-\Delta)}< 1$. It then follows 
that the $G$-action on $L^2_0(G/\Gamma)$ has 
a spectral gap. It was shown 
 that the same holds for any lattice, including 
non-uniform ones (see \cite{BG} and \cite[Lem. 3]{Be}).
\item When $\Gamma$ is a uniform lattice, 
it has been established in \cite[Thm 1.12]{KM} that $L_0^2(G/\Gamma)$ has a strong spectral gap. 
Nevertheless, whether $G/\Gamma$ always has this property 
is still an open problem, even for irreducible 
lattices in $SL_2(\RR)\times SL_2(\RR)$, which are arithmetic 
by Margulis's theorem.

This motivated the formulation of our results 
in a way which makes the dependence on the strong spectral gap 
explicit, namely this assumption is required only for averages 
which are not well-balanced. 
\end{enumerate}
\end{remark}


\section{Statement of results : general $S$-algebraic groups}
\subsection{Ergodic theorems for admissible sets}

Let us now formulate the two basic ergodic theorems for 
actions of the group
$G$ which we will prove in the following section. As usual, 
it is the maximal and exponential-maximal inequalities that 
will serve as our main technical tool in the proof 
of the ergodic theorems.  
The maximal inequalities will also 
be essential later on in establishing the connection between the 
averages on the group and those on the lattice. We will need the following 

\begin{definition}{\bf Weak mixing.}
A measure-preserving action of $G$ on $(X,\mu)$  is called 
\begin{enumerate}
\item weak-mixing if the unitary representation in 
$L_0^2(X)$ does not contain non-trivial finite-dimensional
 subrepresentations. 
 \item totally weak-mixing 
if for every non-trivial normal subgroup, 
the only finite dimensional
 subrepresentation it admits is the trivial one (possibly with 
multiplicity greater than one), or equivalently, if in the 
space orthogonal to its invariants no finite-dimensional 
subrepresentations occur. 
\item We will apply these notions below also to arbitrary 
unitary representations.
\end{enumerate}
\end{definition}

We will use below the notation and terminology established 
in \S 1 and \S 3. In the absence of a spectral gap, we have :

\begin{theorem}\label{th:semisimple-1}

Let $G$ be an $S$-algebraic group as 
in Definition \ref{alg-gps}. 
Let $(X,\mu)$ be a 
totally weak-mixing action of $G$, and let $\{G_t\}$
be a coarsely 
admissible $1$-parameter family, or sequence. 
Then the averages $\beta_t$ satisfy the strong maximal 
inequality in $(L^p,L^r)$ for $p\ge r \ge 1$,  $(p,r)\ne (1,1)$,
Furthermore, if the $G$-action is irreducible, $\beta_t$ 
satisfy, 
\begin{enumerate}
\item the mean ergodic theorem in $L^p$ for $p\ge 1$,
\item the pointwise ergodic theorem in $L^p$ for $p>1$, provided $G_t$ are admissible and left-radial. 
\end{enumerate}
 The conclusions still hold when the action is reducible,  
provided that $\beta_t$ are left-radial and balanced (for the mean theorem) 
or standard radial, well balanced and boundary-regular (for the pointwise theorem). 

\end{theorem}

We note that by Theorem \ref{CAT}, many natural radial averages do indeed satisfy all 
the conditions required in Theorem  \ref{th:semisimple-1}. 
In the presence of a spectral gap, we have

\begin{theorem}\label{th:semisimple-T}
Let $G$ be an $S$-algebraic group as in 
Definition \ref{alg-gps}. 
Let $(X,\mu)$ be a totally weak-mixing 
action of $G$ on a probability space. 
Let $G_t$ be a H\"older 
admissible $1$-parameter family, 
or an admissible sequence when $S$ consists of finite places.  
Assume that either the representation 
 of $G$ on $L^2_0(X)$  
has a strong spectral 
gap,  or that it has a spectral gap and the family $G_t$ 
is well-balanced. Then the averages $\beta_t$ 
satisfy 
\begin{enumerate}
\item the exponential mean ergodic theorem in $(L^p,L^r)$   
for $p\ge r \ge 1$,  $(p,r)\ne (1,1)$,  
\item the exponential strong maximal 
inequality in $(L^p,L^r)$ for $p>r\ge 1$.
\item the exponentially fast pointwise ergodic theorem in $(L^p,L^r)$ 
for $p>r\ge 1$.
\end{enumerate}
\end{theorem}

By Theorem \ref{Appendix:vol} andTheorem \ref{CAT} H\"older-admissible well-balanced families 
do exist in great abundance. 
  
Let us note the following regarding the necessity of the 
assumptions in Theorem \ref{th:semisimple-1} and Theorem 
\ref{th:semisimple-T}. 

\begin{remark}\label{assumption}{\bf 
On the exponential decay of operator norms. }
\begin{enumerate}
\item When $G$ is a product of simple groups but is 
not simple, the assumption that $G_t$ is balanced in Theorem 
\ref{th:semisimple-1} and well-balanced in Theorem 
\ref{th:semisimple-T} is obviously 
necessary in both cases. In the first case, we can 
simply 
take an ergodic action of $G=G_1\times G_2$ which is trivial  
on one factor. In the second, we can take an ergodic action 
with a spectral gap for $G$, but such that  
one of the factors admits an asymptotically 
invariant sequence of 
unit vectors of zero integral. 
\item  In general,  it will be seen below that the only property needed to prove Theorem 
\ref{th:semisimple-T}
 for a $1$-parameter H\"older-admissible 
family acting in $L^2_0(X)$, is the exponential decay of the 
operator norms: $ \norm{\pi_X^0(\beta_t)}\le C e^{-\theta t}$. 
It will be proved in Theorem \ref{Exp decay} below that 
this estimate holds for totally weakly mixing 
actions under the strong 
spectral gap assumption, or when the action has a
 spectral gap and the averages are well-balanced. 
\item We note that for averages given in explicit geometric form, 
 exponential decay of the operator norms  
$\norm{\pi_X^0(\beta_t)}$ can often 
be established directly, see e.g. \cite[Thm. 6]{N3}. The exponentially fast mean ergodic theorem 
for averages on $G$ holds in much greater generality and does not requite admissibility. This 
 fact is very useful in the solution of lattice point counting problems. We refer to  \cite{GN} for a full  discussion and further applications. 
\item Similarly, the radiality assumptions in Theorem 
\ref{th:semisimple-1} is made specifically in order to estimate 
certain spectral expressions that arise in the proof of the pointwise ergodic theorem,
 see remark \ref{radial sets}. At issue is 
the estimate of $\norm{\pi(\partial\beta_t)}$, where 
$\partial\beta_t$ is a singular 
probability measure supported on the boundary of $G_t$. We 
establish the required estimate for left-radial admissible averages in  irreducible actions, or for 
standard radial well-balanced averages in reducible actions. 
This accounts for the statement of Theorem 
\ref{th:semisimple-1}.

\end{enumerate}
\end{remark}

\begin{remark}\label{weak}{\bf Weak mixing. }
\begin{enumerate}
\item The assumption of 
weak-mixing of $G$ is necessary in Theorem 
\ref{th:semisimple-1} 
and in Theorem \ref{th:semisimple-T}, even for 
simple algebraic groups. 
Indeed, it suffices to consider $G=PGL_2(\QQ_p)$, 
and note that it admits a continuous character $\chi_2$ 
onto $\ZZ_2=\set{\pm 1}$. 
It is easily seen that for the natural radial 
averages $\beta_n$ on $G$  
(projecting onto the balls on the Bruhat-Tits tree), the sequence 
$\chi_2(\beta_n)$ {\it does not} converge at all. 
$\chi_2(\beta_{2n})$ does in fact converge, but not to 
the ergodic mean. 
Thus in general, the limiting 
value, if it exists, of $\tau(\beta_t)$ (or subsequences thereof) 
in  finite dimensional representations $\tau$ 
 must be incorporated explicitly 
into the formulation of the ergodic theorems for $G$. 
We refer to 
\cite[\S 10.5]{N5} for a fuller discussion.

Furthermore, note that obviously the ergodic 
action of $G$ on the two-point space $G/\ker \chi_2$ has a spectral 
gap, 
so weak mixing is essential also in Theorem \ref{th:semisimple-T} 
even for simple algebraic groups.

\item Alternatively, another formulation of 
Theorem \ref{th:semisimple-1} for simple algebraic groups  
 (in $L^2$, say) is that convergence to the ergodic 
mean (namely zero) holds 
for an arbitrary ergodic action, when we consider  
the functions in the orthogonal complement of the space spanned 
by all finite-dimensional subrepresentations. The complete picture 
requires evaluating the limits of $\tau(\beta_t)$ for finite-dimensional 
non-trivial representations $\tau$ (if they exist !). 

Similarly, in 
Theorem \ref{th:semisimple-T} if in fact 
 $\norm{\tau(\beta_t)}\le C\exp(-\theta t)$ for finite dimensional 
non-trivial $\tau$, we obtain the same conclusion.

\end{enumerate}
\end{remark}

\begin{remark}\label{G+}{\bf The group $G^+$. }
\begin{enumerate}
\item
Consider an algebraic group 
$G$ defined over a local field $F$ which 
is $F$-isotropic, almost simple and 
algebraically connected as in Definition \ref{alg-gps}. Then 
$G$ contains a canonical co-compact normal 
subgroup denoted $G^+$, which can be defined as the 
group generated by 
the unipotent radicals of a pair of opposite minimal parabolic $F$-subgroups of $G$ 
(see e.g. \cite[\S\S 1.5, 2.3]{M} for a discussion). 

We recall that if $F=\CC$, then 
$G^+=G$, and when $F=\RR$, $G^+$ is the connected component of 
the identity in the Hausdorff topology (which is of finite index in $G$). 
When the characteristic of $F$ is zero, $[G:G^+]< \infty$. 
We note that it is often the case that $G=G^+$
 even in the totally disconnected case. Thus when $G$ is simply connected 
and almost $F$-simple then $G^+=G$, and this includes for example
the groups $SL_n(F)$ and $Sp_{2n}(F)$
 (see e.g. \cite[\S\S 1.4, 2.3]{M}). 

\item A key property of $G^+$ is that it does not admit 
any proper finite index subgroup
(see e.g. \cite[Cor. 1.5.7]{M}). As a result, it follows  
that $G^+$ does not admit any non-trivial finite-dimensional unitary 
representations. Put otherwise, an irreducible unitary representation 
of $G$ is finite-dimensional if and only if it admit a 
$G^+$-invariant unit vector. In 
particular, every ergodic action of $G^+$ is 
weakly mixing, and if it is irreducible, each component is weak-mixing. 
\item  We note the following fact : when $[G:G^+] < \infty$ 
clearly every irreducible non-trivial unitary representation of $G^+$ appears as a subrepresentation  
of the representation of $G$ obtained from it by induction. The induced representation has no $G^+$-invariant unit vectors, and hence its matrix coefficients satisfy the estimates  
that irreducible infinite-dimensional unitary representations of $G$ without $G^+$-invariant unit vectors satisfy. In particular, the $K$-finite matrix coefficients are in $L^p(G^+)$ 
(see Theorem \ref{Lp representations}).

\end{enumerate}
\end{remark}

\subsection{Ergodic theorems for lattice subgroups}

Theorems \ref{th:semisimple-T} 
and \ref{th:semisimple-1} 
will be used to derive corresponding results 
 for {\it arbitrary} measure-preserving
 actions of lattice subgroups of $G$, provided that 
$G$ does not admit non-trivial finite-dimensional 
unitary representation $\tau$.
This condition is necessary, 
and without it the formulation of ergodic 
theorems for the lattice must take into account 
the possible limiting values of $\tau(\lambda_t)$,  
as noted in Remark \ref{weak}. Thus we will formulate our results 
for lattices $\Gamma$ contained in $G^+$, since $G^+$ does
 have the desired property. 
When $G$ is a Lie group, this 
amounts just to assuming the lattice is 
contained in the connnected component in the Hausdorff topology. 
In general, $\Gamma \cap G^+$ is a subgroup of finite index in $\Gamma$, since $\Gamma$ 
is finitely generated and every finitely generated subgroup of $G/G^+$ is finite. 
  
In the absence of a spectral gap, we will prove the following :

\begin{theorem}\label{th:semisimple_lattice-1}
Let $G$ be an $S$-algebraic group, as in Definition 
\ref{alg-gps}, 
and let $\Gamma$ be a lattice subgroup contained in $G^+$. 
Let $G_t$ be an 
admissible $1$-parameter 
family (or an admissible sequence) in $G^+$ 
 and $\Gamma_t=\Gamma\cap G_t$.
Let $(X,\mu)$ be an arbitrary ergodic probability 
measure-preserving action of $\Gamma$. 
Then the averages $\lambda_t$ satisfy the strong maximal 
inequality in $(L^p,L^r)$ for $p\ge r\ge 1$, $(p,r)\ne (1,1)$. 
If the action induced to $G^+$ is irreducible, $\lambda_t$ also 
satisfies 
\begin{enumerate}
\item the mean ergodic theorem in $L^p$ for $p\ge 1$,
\item the pointwise ergodic theorem in $L^p$ for $p>1$, 
assuming in addition that 
 $G_t$ are left-radial.
\end{enumerate}
The same conclusions hold when the induced 
action is reducible, provided the family  $G_t$ is left-radial and balanced (for the mean theorem) 
 or  standard radial, well-balanced and boundary-regular (for the pointwise theorem). 

\end{theorem}

In the presence of a 
spectral gap, we will prove the following :

\begin{theorem}\label{th:semisimple_lattice-T}
Let $G$ be an $S$-algebraic group as in Definition \ref{alg-gps} 
 and $\Gamma$ be a lattice contained in  $G^+$. Let 
$\{G_t\}_{t>0}$ be a
H\"older-admissible $1$-parameter family (or an admissible 
sequence) in $G^+$, 
and $\Gamma_t=\Gamma\cap G_t$. 
Let $(X,\mu)$ be an arbitrary 
probability measure-preserving 
action of $\Gamma$.  
Assume that either the representation of $G^+$ induced  
by the representation of $\Gamma$ on $L^2(X)$ 
has a strong spectral gap, or that it has a spectral 
gap and the family $G_t$ is 
well-balanced. Then 
 $\lambda_t$ satisfy 
\begin{enumerate}
\item the exponential mean ergodic theorem
 in $(L^p,L^r)$ for $p\ge r \ge 1$, $(p,r)\ne (1,1)$,
\item the strong exponential 
maximal inequality in $(L^p,L^r)$ for $p>r\ge 1$, 
\item the exponentially fast pointwise ergodic theorem :
 for every $f\in L^p(X)$, 
$1 < p < \infty$ and almost every $x\in X$ 
$$\abs{\pi_X(\lambda_t)f(x)
-\int_X fd\mu}\le B_r(f,x)e^{-\zeta_p t}$$
where $\zeta_p > 0 $ and $B_r(f,\cdot)\in L^r(X)$, $r< p$. 
\end{enumerate}
\end{theorem}

\begin{remark}
\begin{enumerate}
\item Regarding the assumptions of Theorem 
\ref{th:semisimple_lattice-1}, 
we note that the action 
of $G^+$ induced from the $\Gamma$-action is indeed often 
(but perhaps not always) irreducible. In addition to the obvious 
case where $G$ is simple, irreducibility holds 
(at least for groups 
over fields of zero characteristic) whenever   
the lattice is irreducible and the 
$\Gamma$-action is mixing \cite[Cor. 3.8]{St}. Another 
important case where it holds is 
when the lattice is irreducible and the 
$\Gamma$-action is via a  
dense emedding in a compact group
\cite[Thm. 2.1]{St}, 
or more generally when the $\Gamma$-action is 
isometric.

\item Regarding the assumptions of Theorem 
\ref{th:semisimple_lattice-T}, 
we note that the unitary representation 
of $G$ induced from the unitary representation 
of $\Gamma$ on $L^2(X)$  
always has a spectral gap provided the $\Gamma$-action on 
$(X,\mu)$ does, and it often 
(but perhaps not always) has a strong spectral gap. 
Indeed, by \cite[Ch. III, Prop. 1.11]{M} 
if the lcsc group $G$ has a spectral gap in $L^2_0(G/\Gamma)$ 
and the $\Gamma$-representation on $L^2_0(X)$ has a spectral 
gap,  then so does the representation 
induced to $G$. The existence 
of a spectral gap in $L^2_0(G/\Gamma)$ has long been established 
for all lattices in $S$-algebraic groups. Thus when the sets  
$G_t$ are well-balanced and left-radial, 
the conclusions of Theorem 
\ref{th:semisimple_lattice-T} 
hold provided only that 
the action of 
$\Gamma$ on $(X,\mu)$ has a spectral gap. If the induced action 
is irreducible and $G$ has property $T$, then the induced 
representation has a strong spectral gap and any admissible  
family $G_t$ will do. 

\end{enumerate}
\end{remark}

 Natural radial averages which satisfy all the required properties and thus also the ergodic theorems exist in abundance. To be concrete, let us concentrate on one family of examples, and generalize 
 Theorem \ref{killing} to the $S$-algebraic context. 
 
 Let $G$ be an $S$-algebraic group as in Definition \ref{alg-gps} and
$\ell$ denote the standard $CAT(0)$-metric on 
the symmetric space $X$ or the Bruhat--Tits building $X$
(or their product)  associated to $G$.
Let $\Gamma\subset G$ be a lattice subgroup, $\Gamma_t=G_t\cap \Gamma$, and $\lambda_t$ the uniform averages on $\Gamma_t$. 

\begin{theorem}\label{exist}
Notation being as in the preceding paragraph,  the averages $\lambda_t$  satisfy the mean, maximal and pointwise ergodic theorems in every ergodic action of $\Gamma$ (as in 
Theorem \ref{th:semisimple_lattice-1}). If the action has a spectral gap, then $\lambda_t$ satisfy the 
exponentially fast mean, maximal and pointwise ergodic theorems (as in 
Theorem \ref{th:semisimple_lattice-T}). 

In addition, in every isometric action of $\Gamma$ on a compact  metric space, preserving an ergodic 
probability measure of full support, the averages $\lambda_t$ become equidistributed (as in Theorem \ref{equi-char-0}).  
\end{theorem}

Theorem \ref{exist} is an immediate consequence of  Theorem \ref{CAT},  Theorem 
\ref{th:semisimple_lattice-1}, Theorem \ref{th:semisimple_lattice-T} and Theorem \ref{equi-char-0}. 

We note, however, that 
 the arguments in the Appendix employed to prove Theorem \ref{CAT} apply whenever certain growth and regularity conditions are met, and so Theorem \ref{exist} can in fact be extended to more general families of averages.

\section{Proof of ergodic theorems for $G$-actions}
Our purpose in this section is to prove the ergodic theorems 
for admissible averages on $G$ stated in 
Theorem  \ref{th:semisimple-T} and Theorem \ref{th:semisimple-1}.
Clearly we have to distinguish two cases, 
namely whether the action of $G$ on $X$ has a spectral gap or not.
The arguments that will be employed below 
in these two cases are quite different, but both 
use spectral theory in a material way. We will begin by recalling the 
relevant facts from spectral theory. 
Since we would like to consider all $S$-algebraic groups, 
we will work in the generality of groups admitting an 
Iwasawa decomposition, which we proceed 
to define. This set-up will have the added advantage 
that it incorporates a large class of subgroups 
of groups of automorphism of 
products of Bruhat-Tits buildings. This class contains 
more than all semisimple algebraic groups and $S$-algebraic groups 
and is of considerable interest.

\subsection{Iwasawa groups and spectral estimates}

Let us begin by defining the class of groups to be considered. 

\begin{definition}{\bf Groups with an Iwasawa decomposition.}
\begin{enumerate}\item 
An lcsc group $G$ has an Iwasawa decomposition 
if it has two closed
amenable subgroups $K$ and $P$, with $K$ compact and $G=KP$.

\item The Harish Chandra $\Xi$-function associated with the 
Iwasawa decomposition $G=KP$ of the unimodular group $G$ 
is given by 
$$\Xi(g)=\int_K \delta^{-1/2}(gk)dk$$
where $\delta$ is the left modular function of 
$P$, extended to a left-$K$-invariant
function on $G=KP$. (Thus if $m_P$ is left Haar measure on $P$, 
$\delta(p) m_P$ is right invariant, and 
$dm_G=dm_K \delta(p) dm_P$).
\end{enumerate}
\end{definition}


{\bf Convention}. The definition of an Iwawasa group 
involves a choice of a 
compact subgroup and an amenable 
subgroup. 
When $G$ is the $F$-rational points 
of a semisimple algebraic 
group defined over a locally compact non-discrete field $F$, 
$G$ admits an Iwasawa decomposition, and 
we can and will always choose below $K$ to be 
a good maximal compact subgroup, 
and $P$ a corresponding minimal 
$F$-parabolic group.
This choice 
will be naturally extended 
in the obvious way to $S$-algebraic groups. 

\noindent {\bf Spectral estimates}. 
Iwasawa groups possess a compact subgroup admitting an amenable 
complement, and so 
it is natural to consider the decomposition of 
a representation of $G$ to 
$K$-isotypic subspace. In general 
let $G$ be an lcsc group, $
K$ a compact subgroup, and $\pi : G\to \cU(\cH)$
be a strongly continuous unitary representation, 
where $\cU(\cH)$ is
the unitary group of the Hilbert space $\cH$. 

\begin{definition}{\bf $K$-finite vectors and 
strongly $L^p$-representations.}
\begin{enumerate}
\item A vector $v\in \cH$
is called $K$-finite (or $\pi(K)$-finite) if its orbit under $\pi(K)$
spans a finite dimensional space.  
\item The unitary representation $\sigma$ of $G$ is weakly
contained in the unitary representation $\pi$ if for every $F\in
L^1(G)$ the estimate $\norm{\sigma(F)}\le \norm{\pi(F)}$ holds. Clearly, if
$\sigma$ strongly contained in $\pi$ (namely equivalent to a
subrepresentation), then it is weakly contained in $\pi$. 
\item $\pi$ is called strongly $L^p$ is there exists a dense subspace 
$\cJ\subset \cH$, such that the matrix coefficients
$\inn{\pi(g)v,w}$ belong to $L^p(G)$, for $v,w\in \cJ$.  
\end{enumerate}
\end{definition}

We recall the following spectral estimates, which will play an 
important role below. 
\begin{theorem}\label{tensor}
{\bf Tensor powers and norm estimates.} 
\begin{enumerate}
\item \cite[Thm. 1]{CHH} 
If $\pi$ is strongly $L^{2+\vare}$ for all $\vare > 0$ then $\pi$ is
weakly contained in the regular representation $\lambda_G$.  
\item \cite{Co}\cite{H} (see \cite{HT} for a simple proof) 
If $\pi$ is strongly $L^p$, and $n$ is an integer satisfying 
$n \ge p/2$, then 
$\pi^{\otimes n}$ is strongly contained in $\infty\cdot \lambda_G$.  
\item \cite[Thm 1.1, Prop. 3.7]{N3}
 If $\pi$ is strongly $L^p$, 
and $n_e$ is an even integer satisfying 
$n_e \ge p/2$, then $\norm{\pi(\mu)}\le \norm{\lambda_G(\mu)}^{1/n_e}$
for every probability measure $\mu$ on $G$. If the probability measures $\mu$
and $\mu^\prime$ satisfy $\mu\le C\mu^\prime$ as measures
on $G$, then $\norm{\pi(\mu)}\le C^\prime \norm{\lambda_G(\mu^\prime)}^{1/n_e}$.
\end{enumerate}
\end{theorem}

We begin by stating the following basic spectral estimates for 
Iwasawa groups, which 
are straightforward generalizations of \cite{CHH}.

\begin{theorem}\label{Kfinite estimate}
Let $G=KP$ be a unimodular lcsc 
group with an Iwasawa decomposition,  
and $\pi$ a strongly continuous
unitary representation of $G$. Let 
$v$ and $w$ be two $K$-finite vectors, and denote the dimensions of their spans
under $K$ by $d_v$ and $d_w$. Then the following estimates hold, 
where $\Xi$ is the Harish Chandra $\Xi$-function.
\begin{enumerate}
\item 
 If $\pi$  
is weakly contained in the regular representation, then 
$$\abs{\inn{\pi(g)v,w}}\le \sqrt{d_v d_w}
\norm{v}\norm{w}\Xi(g)\,.$$

\item If $\pi$ is strongly $L^{2k+\vare}$ for 
all $\vare > 0$, then 
$$\abs{\inn{\pi(g)v,w}}\le
 \sqrt{d_v d_w}\norm{v}\norm{w}\Xi(g)^{\frac1k}\,.$$

\end{enumerate}
\end{theorem}
\begin{proof} Part (1) is stated in
 \cite[Thm. 2]{CHH} for semisimple
algebraic groups, but the same proof applies for any 
unimodular Iwasawa
group. 

Part (2) is stated in \cite{CHH} for irreducible
representations of semisimple algebraic groups, 
but the same proof
applies to an arbitrary representation of unimodular 
Iwasawa groups, since it reduces to (1) after taking 
a $k$-fold tensor product.
\end{proof}
\begin{remark} 

\begin{enumerate}
\item The quality of the estimate in  
Theorem \ref{Kfinite estimate} depends 
of course on the structure of $G$.
 For example, if $P$ is normal in $G$ (so that  
 $G$ is itself amenable) then $P$ is unimodular if $G$ is. Then 
$\delta(g)=1$ for $g\in G$ 
 and the estimate is trivial. 
\item In the other direction, Theorem \ref{Kfinite estimate} 
will be most useful when the Harish Chandra
function is indeed in some $L^p(G)$, $p < \infty$, 
so that Theorem 
\ref{tensor} applies.
\item For semisimple algebraic groups the $\Xi$-function 
is in fact in $L^{2+\vare}$ for all
$\vare > 0$, a well-known result due to Harish Chandra 
\cite{HC1}\cite{HC2}\cite{HC3}. 
\item When $\Xi$ is in some $L^q$, $q < \infty$, Theorem 
\ref{Kfinite estimate}(1) implies that any  
representation a tensor power of which is 
weakly contained in the regular representation is strongly 
$L^p$ for some $p$. This assertion uses of course also the 
density of $K$-finite vectors, which is a consequence of the 
 Peter-Weyl theorem. 
\end{enumerate}
\end{remark}

We remark that according to 
Theorem \ref{tensor}(3), it is possible to
bound the operator norm of a given measure by that of its
radialization. Thus in particular for a 
Haar-uniform probability measures on a 
$(K,C)$-radial set the norm is bounded in terms of the 
Haar-uniform probability measure on its radialization.



We now state the following result, 
which summarizes a number of results 
due to \cite{Co}\cite{HM}\cite{BW}. 
 in a form convenient for our purposes.

\begin{theorem}\label{Lp representations}
{\bf $L^p$-representations }\cite{Co}\cite{HM}\cite{BW}.  

Let $F$ be an locally compact 
non-discrete field. 
Let $G$ denote the $F$-rational 
point of an algebraically connected 
semisimple algebraic group which is almost $F$-simple.
Let $\pi$ be a unitary representation of $G$.
 without non-trivial finite-dimensional 
$G$-invariant subspaces (or equivalently 
without $G^+$-invariant unit vectors)
\begin{enumerate}
\item   
If the $F$-rank of $G$ is at least $2$ 
then $\pi$ 
is strongly $L^p$, 
for some fixed $p < \infty$ depending only on $G$. 
\item If the $F$-rank of $G$ is $1$, then 
any unitary representation $\pi$ admitting a spectral 
gap (equivalently, which does not contain an 
asymptotically invariant 
sequence of unit vectors) is strongly $L^p$ for 
some $p=p(\pi)< \infty$. In particular, every 
irreducible infinite-dimensional representation 
has this property.

\end{enumerate}
\end{theorem}

\begin{proof}

When $\pi$ is irreducible, both parts are stated in \cite[Thms. 2.4.2, 2.5.2]{Co} in the Archimedian case, and  in \cite[Thm. 5.6]{HM} in general (based on the theory of leading exponents in \cite{BW}). The passage to general unitary representations with a spectral gap via a direct integral argument presents no difficulty.   
\end{proof}

We remark that an explicit estimate of the relevant exponent $p$ is given by \cite{Li}\cite{LZ} in the Archimedian case, and by \cite{Oh}, in general.  The pointwise bounds for $K$-finite matrix coefficients developed in \cite{CHH}, \cite{H},\cite{HT} and in the general cae in  \cite[\S 5.7]{Oh},  imply the bound stated in Theorem \ref{Kfinite estimate}(2), 
showing that the matrix coefficients are indeed in   
$L^p(G)$, where $p$ depend 
only on $G$.

Note that any faithful unitary representation of an $S$-algebraic  group (as in Definition \ref{alg-gps})
with property $T$ 
is strongly $L^p$. 


\subsection{Ergodic theorems in the presence of a spectral gap}

When a spectral gap is present, 
a strong exponential maximal inequality holds for
general H\"older 
families of probability 
measures $\nu_t$ on $G$. For a proof we refer to 
\cite{MNS} and \cite{N3}, where the relation between the rate 
of exponential decay and the parameters $p$ and $r$ below is 
fully explicated.

\begin{theorem}\label{G-actions, spectral gap}
{\bf Exponential maximal inequality in the
presence of a spectral gap}\cite[Thm. 4]{N3}.
Let $G$ be an lcsc group, and assume that the family of 
probability 
measures $\nu_t$, 
is uniformly locally 
H\"older continuous in the total variation norm, namely 
$\norm{\nu_{t+\vre}-\nu_t}\le C\vre^a$, for all $t\ge t_0$ 
and $0 < \vre \le \vre_0$. Assume also that it is roughly monotone, 
namely $\nu_t\le C \nu_{[t]+1}$, where $C$ is fixed. 
\begin{enumerate} \item Assume that  
$\pi_X^0(\nu_t)$ 
have exponentially decaying 
norms in $L^2_0(X)$.
Then the strong exponential
 maximal inequality in $(L^p,L^r)$ holds in any 
  probability-measure-preserving action of $G$, and thus also 
the exponential pointwise ergodic theorem holds in $(L^p,L^r)$, 
for any $p> r > 1$.  
\item In particular, if the representation
$\pi^0_X$ on $L_2^0$ is strongly
$L^p$, and $\nu_t$ have exponentially decaying norms
 as convolution operators on $L^2(G)$, then the previous conclusion holds.
\end{enumerate}
\end{theorem}

\begin{remark} For future reference, let us recall 
the following simple observation, noted already 
in \cite[Thm. 2]{N3}. For {\it any sequence} of
 averages $\nu_n$, the  
exponential-maximal inequality in $(L^p,L^p)$, $1 < p < \infty$ 
is an immediate consequence of 
 the exponential decay condition on the norms
 $\norm{\pi_X^0(\nu_t)}$. So is the exponentially fast pointwise 
ergodic theorem, and both follow simply by considering 
the bounded operator $\sum_{n=0}^\infty e^{n\theta/2}\pi^0_X(\nu_n)$ 
on $L^2_0(X)$, and then using Riesz-Thorin interpolation.
\end{remark}

The next step is to establish the exponential decay conditions on 
the norms, when the averages are admissible. 
Note that according to 
Theorem \ref{tensor}(3), it is possible to
bound the operator norm of (say) a given 
 Haar-uniform probablity measures on a 
$(K,C)$-radial and in terms of its radialization. 
We formulate this fact as follows. 

\begin{proposition}\label{radial estimate}
{\bf Radialization estimate.}
Let $G=KP$ be an lcsc unimodular 
Iwasawa group and $\pi$ a strongly $L^p$-representation. 
Let $B$ be a set of positive finite measure and  
$\beta=\chi_B/\vol(B)$. 
\begin{enumerate}
\item  If $B$  
is bi-$K$-invariant, then  
$\norm{\lambda(\beta)}=\frac{1}{m(B)}\int_{B} \Xi_G(g)dm_G(g)$.
\item If $B$ is a $(K,C)$-radial set and $\tilde{B}=KBK$, then 
$$\norm{\pi(\beta)}\le C^\prime\left
(\frac{\int_{KBK}\Xi(g)dm_G(g)}{\vol(KBK)}\right)^{1/n_e}=C^\prime 
\norm{\lambda(\tilde{\beta})}^{1/n_e}$$
provided $n_e$ is even and $n_e > p/2$. 
\end{enumerate}
Both statments hold of course for every absolutely continuous 
$(K,C)$-radial 
probability measure (with the obvious definition of $(K,C)$-radial measures). 
\end{proposition}

\begin{proof}

Let us first recall 
the following estimate from \cite{CHH}, 
dual to C. Herz majorization principle. 
 The spectral norm of the convolution 
operator $\lambda_G(F)$ on $L^2(G)$ 
is estimated by : 
$$\norm{\lambda_G(F)}\le
\int_G\left(\int_K\int_K\abs{F(kgk^\prime)}^2dkdk^\prime
\right)^{1/2}
\Xi(g)dg\,\,$$
for any measurable function $F$ on $G$ for 
which the right hand side is
finite. 

Now consider a set $B$ which is  
$(K,C)$-radial. 
Clearly a comparison of the convolutions with the 
normalized probability measures 
immediately gives 
$$\norm{\lambda_G(\beta)}\le \frac{C}{\vol(KBK)}
 \norm{\lambda_G(\chi_{KBK})}\,.
$$
Utilizing the previous inequality for the bi-$K$-invariant 
measure $\tilde{\beta}$ uniformly 
distributed on $KBK$, 
we clearly obtain 
$$\norm{\lambda_G(\beta)}\le \frac{C}{\vol (KBK)} 
\int_{KBK} \Xi(g)dm_G(g)\,.$$
The equality in (1) for bi-$K$-invariant sets $\tilde{B}$ follows from 
the fact that since $P$ is amenable, The representation of $G$ 
on $G/P$ induced from the trivial representation of $P$ is weakly 
contained in the regular representation of $G$. $\Xi$ is a diagonal 
matrix coefficient of the latter representation, by definition,
so that 
$\frac{1}{m(\tilde{B})}
\int_{\tilde{B}} \Xi_G(g)dm_G(g)\le \norm{\lambda(\tilde{\beta})}$.

Now by definition 
of weak containment, the same inequality hold  
 also for unitary
representations $\pi$ which are weakly contained in the 
regular
representation. Taking tensor 
powers and using Theorem \ref{tensor}, 
we obtain a norm bound for $\norm{\pi(\beta)}$ in 
$L^p$-representations as well. 
\end{proof}

\begin{remark}\label{KS}{\bf Bounding spectral norms 
using the Kunze-Stein phenomenon.}

To complement the latter result let us recall, 
as noted in \cite[Thm. 3,4]{N3}, that
 it is also possible to estimate the operator norm 
of general (not necessarily $(K,C)$-radial)
 averages using the Kunze-Stein 
phenomenon, 
provided the representation 
is strongly $L^p$, $p < \infty$. 

1) Indeed, for the spectral norm, namely when 
$f\in C_c(G)$ 
and $\beta=\chi_B/m_G(B)$ ($B$ a bounded set), 
we have by the Kunze-Stein phenomenon, provided $1 < p < 2$ 
$$\norm{\beta\ast f}_{L^2(G)}\le 
C_p\norm{\beta}_{L^p(G)}\norm{f}_{L^2(G)}\,.$$
Taking $p=2-\vre$, and $f$ ranging over 
functions of unit 
$L^2(G)$-norm, we conclude that 
$$\norm{\lambda(\beta)}\le C^\prime_{\vre}
\norm{\beta}_{L^{2-\vre}(G)}= C^\prime_\vre 
m_G(B)^{-(1-\vre)/(2-\vre)}\le 
C''_\vre m_G(B)^{-1/2+\vre}\,\,.$$

2) If the representation $\pi$ satisfies 
$\pi^{\otimes n_e}\subset \infty\cdot \lambda_G$ (e.g. if it is
strongly $L^p$ and $n_e$ is even with $n_e\ge p/2$) then by part 
1) and Theorem \ref{tensor}(3) (see also \cite[Thm. 4]{N3})  
$\pi(\beta)\le C_\vre m_G(B)^{-1/(2n_e)+\vre}\,.$

3) It follows from (2) that if 
$\beta_t$ is a family of Haar-uniform averages on the sets 
$G_t$, and $\pi^{\otimes n_e}\subset \infty\cdot \lambda_G$, 
then the rate of volume growth of $G_t$ determines the largest 
parameter $\theta$ satisfying $\norm{\pi(\beta_t)}\le 
A_\vre e^{(-\theta+\vre) t}$ (for every $\vre > 0$)  
as follows :
$$\theta=\liminf_{t\to \infty}-\frac1t\log\norm{\pi(\beta_t)}
=\frac{1}{2n_e}\limsup_{t\to \infty} \frac1t \log m_G(B_t)\,\,.$$

\end{remark}

We now turn to establishing the necessary decay estimates 
for the norms of the operators $\pi(\beta_t)$, for 
admissible families $G_t$. 

\begin{theorem}\label{Exp decay}{\bf Exponential decay of
 operator norms.}
Let $G$ be $S$-algebraic as in 
Definition \ref{alg-gps}, and let $\sigma$ 
be totally weak-mixing unitary representation 
of $G$. Let 
$G_t$ be a coarsely 
admissible $1$-parameter family or sequence.
Assume that $\sigma$ has a strong 
spectral gap, or that it has a 
spectral gap and $G_t$ are well-balanced. 
Then for some 
$C$ and $\delta=\delta_\sigma > 0$ 
depending on $\sigma$ and $G_t$,  
$$\norm{\sigma(\beta_t)}\le Ce^{-\delta t}\,.$$
\end{theorem}

\begin{proof}
We have already shown in Proposition \ref{KC-radial} that
 coarsely admissible families are $(K,C)$-radial, and so let us
 use Proposition \ref{radial estimate}. 
If the representation $\sigma$ happens to be  
strongly $L^p$, then denoting 
the radializations by $\widetilde{G}_t=KG_tK$, 
the norm of $\sigma(\beta_t)$ is estimated by 
$$\norm{\sigma(\beta_t)}\le 
{C^\prime}\left(\frac{\int_{\tilde{G}_t}
 \Xi_G(g)dm_G(g)}{m_G(\tilde{G}_t)}\right)^{1/n_e}$$ 
where $n_e > p/2$ is even and 
$\Xi_G$ is the Harish Chandra $\Xi$-function. 

Recall that coarsely admissible sets satisfy the minimal growth 
condition $S^n\subset G_{an+b}$. It follows of course 
that their radializations satisfy $\tilde{S}^n\subset 
\tilde{G}_{a'n+b'}$, for a compact bi-$K$-invariant generating 
set $\tilde{S}$. The standard 
estimates of $\Xi_G(g)$ (see \cite{HC1}, \cite{HC2},\cite{HC3}) 
now imply that $\norm{\sigma(\beta_t)}$ 
decays exponentially.

 Now the assumption that the representation is strongly 
$L^p$ is satisfied 
when the representation is totally weak mixing 
and has a strong 
spectral gap. Indeed this follows immediately 
from Theorem
 \ref{Lp representations}, noting that (almost) every 
irreducible representation appearing in the direct integral 
decomposition of $\sigma$ w.r.t. a simple subgroup must infinite dimensional.

We note that the strong spectral gap assumption is 
necessary here. Indeed, consider the tensor product 
of an irreducible  principal series representation of $G=PSL_2(\QQ_p)$ and  
 a weak mixing representation of $G$ which admits an 
asyptotically invariant sequence of unit vectors. 
Then $\sigma$ has a spectral gap
 (as a representation of $G\times G$) and is totally weak mixing 
but is not strongly $L^p$ for any finite $p$.

To handle the case where $\sigma$ totally weak mixing but 
is not strongly $L^p$, 
let us recall that we now assume 
$G_t$ are well-balanced, 
and also coarsely admissible. It follows 
that the radialized sets $\tilde{G}_t$ are also well-balanced. 
Indeed,  since      
$G_t\subset KG_tK \subset G_{t+c}$ for every coarsely admissible 
family, clearly also (refering to Definition \ref{balanced})  
$G_t\cap H_IQ\subset KG_tK\cap H_I Q\subset G_{t+c}\cap H_I Q$ for every 
compact subset $Q$ of $H_J$. Taking $t=an$, $Q=S_I^n$ and using 
that $G_t$ are well-balanced, the claim follows. 
Now, if $\bar{\sigma}$ is the representation 
conjugate to $\sigma$, then for any probability measure 
$\nu$ on $G$ and every vector $u$ we have 
$$\norm{\sigma(\nu)u}^2=\inn{\sigma(\nu^\ast\ast\nu)u,u}\le 
\int_G\abs{\inn{\sigma(g)u,u}}d(\nu^\ast \ast \nu)(g)
$$
$$\le \left(\int_G
\left(\inn{\sigma(g)u,u}\right)^2 d(\nu^\ast \ast \nu)\right)^{1/2}
$$
$$
= \left(\inn{\sigma\otimes 
\bar{\sigma}(\nu^\ast\ast\nu)(u\otimes \bar{u}),
u\otimes \bar{u}}\right)^{1/2}
\le
\norm{\sigma\otimes\bar{\sigma}(\nu)(u\otimes\bar{u})}\,.
$$
Hence it suffices 
to prove that $\norm{\sigma\otimes \bar{\sigma}(\beta_t)}$ decays 
exponentially. But note 
that the diagonal matrix 
coefficients which we are now considering 
are all non-negative. Thus 
for each vector $u$, and every $(K,C)$-radial measure $\nu$, 
using Jensen's inequality and the previous argument   
$$\norm{\sigma(\nu)u}^4=\inn{\sigma(\nu^\ast\ast\nu)u,u}^2\le 
\abs{\inn{\sigma\otimes \bar{\sigma}(\nu^\ast\ast\nu)(u
\otimes \bar{u}),u\otimes \bar{u}}}
$$
$$
\le C^2 \inn{\sigma\otimes \bar{\sigma}
(\tilde{\nu}^\ast\ast\tilde{\nu})(u\otimes \bar{u}),
u\otimes \bar{u}}=C^2\norm{\sigma\otimes \bar{\sigma}(\tilde{\nu})
(u\otimes \bar{u})}^2\,.$$
We conclude that if the norm of 
$\sigma\otimes \bar{\sigma}(\tilde{\beta}_t)$ 
decays exponentially, so does the norm of 
$\sigma(\beta_t)$. 
Now if $\sigma$ has a spectral gap, and is totally weak-mixing, 
then $\sigma\otimes \bar{\sigma}$ has the same properties. 
This claim follows from 
Theorem \ref{Lp representations} and 
Theorem \ref{Kfinite estimate}.
 Indeed if an asymptotically 
invariant sequence of 
unit vectors exists in $\sigma\otimes \bar{\sigma}$, 
then there is also such a sequence 
which consists of $K$-invariant vectors, so 
that we can restrict attention to the spherical spectrum. 
But a sequence consisting of convex sums of products of  
normalized positive definite spherical functions
 cannot converge to $1$ uniformly on compact sets unless the 
individual spherical functions that occur in them have the same 
property, and this contradicts the spectral
 gap assumption on $\sigma$. Total weak-mixing follows from a 
direct integral decomposition and the fact that the 
tensor product of two irreducible infinite-dimensional 
unitary representations of an simple group do not have 
finite-dimensional subrepresentations.

Thus we are reduced to establishing 
the norm decay of a bi-$K$-invariant 
coarsely admissible family in a totally weak-mixing unitary 
representation with a spectral gap, 
which we continue to denote by $\sigma$.  

To estimate the norm of $\sigma(\tilde{\beta}_t)$  under these 
conditions write 
$G=G_1\times G_2$ where $G_1$ has property $T$ 
and $G_2$ is a product of groups of split rank one.  
Any irreducible 
unitary representation of $G$ is a tensor product of irreducible 
unitary representations of $G_1$ and $G_2$.  
Since $\tilde{\beta}_t$ are bi-$K$-invariant measure, we can 
clearly restrict our attention to infinite-dimensional 
spherical representations, namely 
those containing a $K$-invariant unit vector, and then estimate 
the matrix coefficient given by the spherical function.  
The non-constant 
spherical functions $\varphi_s(g)$ 
on $G$ are given by 
$\varphi_{s_1}(g_1)
\varphi_{s_2}(g_2)$, 
where  at least one of the 
factors is non-constant. If $\varphi_{s_1}(g_1)$ is non-constant, 
then again it is a multiple of spherical function on the 
simple components of $G_1$, one of which is non-constant. 
 $\varphi_{s_1}(g_1)$ is then bounded, according to Theorem 
\ref{Kfinite estimate} and Theorem \ref{Lp representations}, 
 by the function $ \Phi : 
g_1\mapsto \Xi_{G/H}(p_{G/H}(g_1))^{1/n}$ for some fixed $n$ 
depending only on $G_1$, where $G/H$ is one of the simple factors 
of $G_1$ and $p_{G/H}$ 
the projection onto it. Using 
the standard estimate of the $\Xi$-function, together with our 
assumption that the $G_t$ are well-balanced, 
we conclude that for $\eta_1> 0$ depending only 
on $G_1$ and $G_t$ 
$$\frac{1}{m_G(G_t)}\int_{G_t}
 \abs{\varphi_s(g)}dm_G(g)\le Ce^{-\eta_1 t}\,.$$

Now, if $\varphi_{s_1}(g_1)$ 
is constant, then the 
representation of $G$ in question is trivial on $G_1$ and 
factors to a non-trivial irreducible 
representation of $G_2$. The spherical functions of the 
complementary series on a 
split rank one group are the only ones we need to consider, 
and they have a simple parametrization as a subset of  
an interval. We can take to be say $[0,1]$, 
 where $0$ corresponds to the Harish Chandra function 
and $1$ to the constant function. 
The function  $\varphi_{s_2}(g_2)$ is a product of spherical 
functions on the real rank factors. 
The assumption that the original representation 
$\sigma= \pi^0_X$ has a spectral gap 
implies that for some $\delta > 0$, the parameter of at least 
one factor is outside $[1-\delta,1]$. The bound $\delta$ depends 
only on the original representation $\sigma$ and is uniform over 
the representations of $G_2$ that occur. It is then easily seen 
using the estimate of the $\Xi$-function in the split rank 
one case and the fact that $\tilde{G}_t$ are well-balanced, that 
$\varphi_s(g)$ also satisfies the foregoing estimate, for 
some $\eta_1^\prime > 0$. This 
gives a uniform bound over all irreducible unitary 
spherical representations of $G$ that 
are weakly contained in $\sigma$, and 
  it follows 
that $\norm{\pi_X^0(\tilde{\beta}_t)}\le Ce^{-\eta t}$, where 
$\eta=\eta(\sigma) > 0$. 
\end{proof}

{\it Proof of Theorem \ref{th:semisimple-T}}.

We now check that the assumption of 
Theorem \ref{G-actions, spectral gap} are satisfied for 
the $1$-parameter families of 
averages and the representations under consideration. Clearly, admissible 
families are also roughly monotone,  
and are uniformly locally Lipschitz continuous 
in the $L^1(G)$-norm, as noted already in \ref{loc.Lip}(2). 
Furthermore, coarsely admissible (and in particular, admissible) 
families satisfy the 
exponential decay condition in $L^2_0(X)$ in the
 representation under consideration, 
in view of Theorem \ref{Exp decay}. 
Thus the exponentially fast mean, pointwise and maximal ergodic 
theorems holds for $1$-parameter admissible families.


Finally, the case of sequence of admissible averages  
does not require an appeal to 
Theorem \ref{G-actions, spectral gap}, just to the remark 
following it, as well as to Theorem \ref{Exp decay}. 

This concludes the proof of 
all parts of Theorem \ref{th:semisimple-T}. \qed

\begin{remark}\label{holder}{\bf H\"older families.}
Clearly Theorem \ref{th:semisimple-T} 
remains valid with the same proof 
provided that $G_t$ satisy, for some $ 0 < a \le 1$ 
the following H\"older condition :
\begin{equation}\label{holder1}
m_G(G_{t+\vre})\le (1+c\vre^a)\cdot  m_G(G_{t})\,\,
\end{equation}
 rather than the Lipshitz 
condition above. 

Furthermore, we can also weaken 
condition (\ref{eq:1}) in the definition of admissibility by 
the H\"older condition 
\begin{equation}\label{holder2}
\mathcal{O}_\vre
\cdot G_t\cdot \mathcal{O}_\vre 
\subset G_{t+c\vre^a}\,\,.
\end{equation}

The H\"older assumptions are sufficient also for 
the proof of Theorem 
\ref{th:semisimple_lattice-T} 
as well as 
parts (4) and (5) of Theorem 
\ref{reduction} and part (4) 
of Theorem \ref{reduct-lattice-points}. 

\end{remark}

\subsection{Ergodic theorems in the absence of
a spectral gap, I}

We now turn to the proof of Theorem \ref{th:semisimple-1} and 
 consider 
actions which do 
not necessarily admit a spectral gap. We start by proving the 
strong maximal inequality for admissible 
families (and some more general ones),
and in the section that follows we 
establish pointwise convergence
 (to the ergodic mean) on a dense subspace. In the course 
of that discussion the mean ergodic theorem, namely 
norm convergence (to the ergodic mean) on a dense subspace 
will be apparent.
As is well-known these three ingredients suffice to prove 
Theorem \ref{th:semisimple-1} completely (see e.g. \cite{N5} for a full discussion).

\subsubsection{The maximal inequality}

Consider a family $\beta_t$ of probability measures on 
$G$, where $\beta_t$ is the Haar-uniform 
average on $G_t$, and $G_t$ are coarsely 
admissible. As noted in Proposition 
\ref{KC-radial} coarse admissibility implies that $G_t$ are 
a family of $(K,C)$-radial sets, 
with $C$ fixed and independent of $t$. As before, we denote by $\tilde{\beta}_t$ 
the Haar uniform averages of the sets $\widetilde{G}_t=KG_tK$, and again note that  
$\beta_t\le C \tilde{\beta}_t$ as measures on $G$. Hence for $f\ge 0$ we have almost surely  
$$f^\ast_\beta(x)=\sup_{t> t_0}\pi_X(\beta_t)f(x)\le C
 \sup_{t > t_0}\pi_X (\tilde{\beta}_t)f(x)= C
f^\ast_{\tilde{\beta}}(x)\,\,,$$
so that it suffices to
 prove the maximal inequalities for the averages $\tilde{\beta}_t$. 
Furthermore, 
if $G_t$ is coarsely admissible then 
clearly so are the sets $\widetilde{G}_t=KG_tK$.

Recall now that that coarse 
admissibility implies 
$(\mathcal{O}_r,D)$-ampleness for some constants $(r,D)$, as noted in 
Proposition \ref{KC-radial}.

Thus the maximal inequality for $\tilde{\beta}_t$ follows 
from the following  
\begin{theorem}\label{ample sets}
{\bf Maximal inequality for ample sets}(see \cite[Thm 3]{N4}). 
Let $G$ be an $S$-algebraic group as in Definition 
\ref{alg-gps}, and let $K$ be of finite index in a maximal compact subgroup. 
Let $(X,\mu)$ be a totally weak-mixing probability measure 
preserving action of $G$.
For set $E\subset G$ of positive finite measure, let $\nu_E$ 
denote the Haar-uniform average supported on $E$. Fix positive 
constants $r > 0$ and $D > 1$, and consider the 
maximal operator   
$$\mathcal{A}^\ast f(x)=\sup\set{\abs{\pi_X(\nu_E)f(x)}\,\,:\,\, 
E\subset G \text{  and $E$ is $(\mathcal{O}_r,D,K)$-ample }}\,.$$
Then $\mathcal{A}^\ast  f$ satisfies the maximal inequality ($1< p \le \infty$) 
$$\norm{\mathcal{A}^\ast f}_{L^p(X)} \le B_p(G,r,D,K)
\norm{f}_{L^p(X)}\,.$$
\end{theorem}
\begin{proof}
Theorem \ref{ample sets} is proved in full 
for connected semisimple 
Lie groups with finite center in \cite{N4}. 
The same proof applies to our present more general 
context without essential 
changes, as we now briefly note. First, if $G$ is a 
totally disconnected almost 
 simple algebraically 
connected non-compact algebraic group with property $T$, 
the analog of \cite[Thm. 2]{N4}, namely the exponential-maximal 
inequality for the cube averages defined there follows from  
the exponential 
decay of the norm of the {\it sequence} (in this case) of 
cube averages. The exponential decay in the space 
orthogonal to the invariants is 
assured by our assumption that the action is totally 
weak-mixing, which implies that Theorem 
\ref{Lp representations} and  Theorem \ref{Kfinite estimate} 
can be applied. Then
standard estimate of the $\Xi$-function yield the desired 
conclusion for the cube averages. 

Second, for groups of split rank one  
the maximal inequality for the sphere averages is established 
in \cite{NS} when the Bruhat-Tits tree has even valency, and
the same method gives the general case (using the description 
of the spherical function given e.g. in \cite{N0}). The fact 
that the maximal inequality for cube averages holds 
on product groups if it holds for the components 
is completely elementary, as in \cite[\S 2]{N4}. 
Finally, the fact 
that the maximal inequality for cube averages 
implies the maximal inequality for ample set in a given group 
depends only on analysis of the volume density associated with 
the Cartan decomposition and thus only on the root system, and 
the argument in \cite[\S 4]{N4} generalizes without difficulty.

This establishes the maximal inequality for every $S$-algebraic group as in 
Definition \ref{alg-gps}. 
\end{proof}

\subsection{Ergodic theorems in the absence of
a spectral gap, II}

\subsubsection{Pointwise convergence on a dense subspace}

In the present subsection we will make 
crucial use of the absolute continuity  
property established for admissible $1$-parameter of 
averages in Proposition \ref{prop:abs.cont}. 
This property implies 
 that $t\mapsto \beta_t$ is almost 
surely differentiable 
(in $t$, and w.r.t. the $L^1(G)$-norm) with globally bounded derivative. 

When $\beta_t$ are the Haar uniform averages on an 
increasing family of compact sets $G_t$, 
Almost sure differentiabilty is 
equivalent with the almost sure existence of the limit 
$$\lim_{\vre \to 0}\frac{1}{\vre}
\frac{m_G(G_{t+\vre})-m_G(G_t)}{m_G(G_t)}\,\,,$$
and for admissible families  the limit is 
uniformly bounded as a function of $t$.  
The almost sure differentiability will allows us to make use of 
a certain Sobolev-type argument developed originally in \cite{N1}.

The uniform local Lipshitz continuity for the averages is 
a somewhat stronger property than uniform local 
H\"older continuity, which was 
the underlying condition in  
 the case where the action has a spectral gap. However,  
the Lipshitz condition allows us to dispense with
the assumption of 
exponential decay of $\norm{\pi^0_X(\beta_t)}$. 

We note however that the (ordinary) strong $L^p$-maximal inequality 
holds for much more general averages, namely 
 under the sole conditions that $G_t$ 
are $(K,C)$-radial and their radializations are 
$(\mathcal{O}_r,D)$-ample averages on an 
$S$-algebraic group.
It is only pointwise convergence on a 
dense subspace that requires the additional 
regularity assumption of 
almost sure differentiability, which follows from the 
uniform local Lipshitz condition. 

Finally, we remark that the
 argument we give below 
is based solely on the spectral estimates
described in the previous sections. 
Thus it is does not require extensive 
considerations related to classification of unitary 
representations, and  
applies to all semisimple algebraic groups 
(and other Iwasawa groups).

Let $(G,m_G)$ denote an lcsc group $G$ with a left Haar measure. 
Let $N_t$, $t\in \RR_+$ be an admissible family. 
Then $N_t$ is an increasing 
family of bounded sets of positive 
measure, satisfying, without loss of generality  $N_t=\cap_{s > t}N_s$. 
Let $g\mapsto \abs{g}$ be the gauge 
$\abs{g}=\inf\set{s\,;\, g\in N_s}$. 
This condition implies that 
$N_t$ are determined by their gauge via 
$N_t=\set{g\in G\,;\, \abs{g}\le t}$. 
Thus the gauge is a measurable proper function with values in 
$\RR_+$. Define $\nu_t$, $t\in \RR_+$ to be the one-parameter
family of probability measures with compact supports on $G$,
absolutely continuous w.r.t. Haar measure, whose density is given by
the function $\frac{1}{m_G(N_t)}\chi_{N_t}(g)$. The map 
$t \mapsto \nu_t$ is a uniformly locally Lipshitz 
 function from $\RR_+$ to $L^1(G)$,
w.r.t. the norm topology, by assumption. 

We let $S_t=\set{g\,;\, \abs{g}=t}$, and clearly 
$N_t=\coprod_{0<  s\le t}S_s$ is a disjoint union. The map
$g\mapsto \abs{g}$ projects Haar measure on $G$ onto a measure 
on $\RR_+$, which is absolutely  
continuous measure w.r.t. Lebesgue measure on $\RR_+$, by Proposition 
\ref{prop:abs.cont}. 
The measure disintegration formula gives the representation  
$m_G=\int_0^\infty  m_r dr$, where $ m_r$ is a
 measure on $S_r$, defined for almost all $r$. Thus we can
write for any $F\in C_c(G)$ 
$$\nu_t(F)=\frac{\int_{N_t}F dm_G}{m_G(N_t)}=\frac{1}{m_G(N_t)}\int_0^t
\frac{m_r(F)}{m_r(S_r)} m_r(S_r)dr=$$
$$=\int_0^t \partial \nu_r(F)\psi_t(r)dr\,\,.$$
Here $\partial \nu_r=m_r/m_r(S_r)$ is a probability measure
on $S_r$ (for almost every $r$), and the
density $\psi_t(r)$ is given by $\psi_t(r)=m_r(S_r)/m_G(N_t)$. 
Here $\psi_t(r)$ is a measurable function, defined 
almost surely w.r.t. Lebesgue measure on $\RR_+$, and is almost surely
positive for $r\le t$. For any given
continuous function $F\in C_c(G)$, 
$\nu_t(F)$ is an absolutely continuous function on
$\RR_+$, given by integration against the $L^1$-density
$\partial\nu_r(F)\psi_t(r)$ (which is almost everywhere defined). 
In particular, $\nu_t(F)$ is differentiable almost everywhere, and 
its derivative is given, almost everywhere, as follows

\begin{proposition}\label{derivative}
Assume $N_t$ give rise to an absolutely continuous measure on
 $\RR_+$, as above. Then for almost all $t$ : 
$$\frac{d}{dt}(\nu_t(F))=
\frac{m_t(S_t)}{m_G(N_t)}\left(\partial \nu_t-\nu_t\right)(F)$$
\end{proposition}
\begin{proof}
We compute :

$$\frac{d}{dt}\nu_t
=\frac{d}{dt}\left( \frac{1}{m_G(N_t)}
\int_0^t m_r dr\right)=
$$
$$
\left(\frac{1}{m_G(N_t)}\right)^\prime \int_0^t m_r dr\,+\,
 \frac{1}{m_G(N_t)} \cdot m_t=$$
$$-\frac{m_G(N_t)^\prime}{m_G(N_t)^2}\int_0^t m_r dr+
\frac{m_t(S_t)}{m_G(N_t)}\partial \nu_t=$$
$$=\frac{m_t(S_t)}{m_G(N_t)}\left(\partial \nu_t-\nu_t\right)$$
We have used that $m_G(N_t)^\prime=m_t(S_t)$ for almost all $t$,  
which is a consequence of the Lebesgue differentiation theorem on 
the real line. 
\end{proof}

Let us note that when a uniform local Lipschitz condition 
is satisfied by $\log m_G(N_t)$, namely when  
$m_G(N_{t+\vare})\le (1+c\vare)m_G(N_t)$ for 
$0 < \vare\le 1/2$, and all $t \ge 1$, 
it results in a uniform estimate of the
ratio of the ``area of the sphere'' (i.e. $S_t$)
 to the ``volume of the ball'' (i.e. $N_t$).
Thus we have :
\begin{corollary}\label{unif bound}
Assume that $N_t$ is an admissible $1$-parameter family. Then 
$$\frac{m_G(N_{t+\vare})-m_G(N_t)}{m_G(N_t)}=
\frac{\int_t^{t+\vare} m_r(S_r) dr}{m_G(N_t)}\le c\vare $$
so that for almost all $t$ we have the uniform bound 
$$\frac{m_t(S_t)}{m_G(N_t)}=\frac{1}{m_G(N_t)}
 \lim_{\vare \to 0}\frac{1}{\vare}\int_t^{t+\vare}m_r (S_r)dr
 \le c\,\,.$$
\end{corollary}

The existence of the derivative almost everywhere 
of $\nu_t(F)$ imply, 
in particular, that for every $F\in C(G)$, and 
for {\it every} $ t > s > 0$, the following identities hold :  
$$\nu_t(F)=\int_0^t \frac{d}{dr}\nu_r (F)dr 
\text{             and            }  
\nu_t(F)-\nu_s(F)=
\int_s^t \frac{d}{dr} \nu_r (F)dr\,\,.$$
It follows that corresponding 
 equalities hold between the underlying
 (signed) measures on $G$, namely 
$$\nu_t= \int_0^t\frac{d}{dr}\nu_r dr = 
\int_0^t 
\frac{m_r(S_r)}{m_G(N_r)}
\left(\partial \nu_r-\nu_r\right)dr\,\,.$$

Thus the derivative  $\frac{d}{dr}\nu_r$ is a multiple 
of the difference between 
two probability measures on $G$ (for almost every $r$).  
Every bounded measure on $G$ naturally gives rise to   
a bounded operator on the representation space. 
We can thus conclude the following relations 
between the corresponding operators defined in any $G$-action.

\begin{corollary}\label{diff vectors}{\bf Differentiable vectors.}
\begin{enumerate}
\item For any strongly continuous 
 unitary representation, and any vector 
$u\in \cH$,  $r\mapsto \pi(\nu_r)u$ is almost surely  
differentiable in $r$ (strongly, namely in the norm topology), 
and the following holds: 
$\pi(\nu_t)u=\int_0^t \frac{d}{dr}\pi(\nu_r)u dr $ and 
$\pi(\nu_t)u-\pi(\nu_s)u=
\int_s^t \frac{d}{dr}\pi( \nu_r)u dr$. 
\item Consider a measurable $G$-action 
on a standard Borel probability space, 
and a function $u(x)\in L^p(X)$, $1 \le p < \infty$
 for which $g\mapsto
 u(g^{-1}x)$ is continuous in $g$ for almost every $x\in X$.  
Then the expression : $\pi(\nu_t)u(x)=\int_G u(g^{-1}x)d\nu_t(g)$ is 
differentiable in $t$ for almost every
$x\in X$ and almost every $r\in \RR_+$, 
and the following almost sure identities hold : 

$$\pi(\nu_t)u(x)=\int_0^t \frac{d}{dr}\pi(\nu_r)u(x) dr $$ 
and 
$$\pi(\nu_t)u(x)-\pi(\nu_s)u(x)=
\int_s^t\frac{d}{dr}\pi( \nu_r)u(x) dr\,\,.$$ 
\end{enumerate}
\end{corollary}

\begin{remark}
Note that in Corollary \ref{diff vectors}(2), the space 
of vectors $u$ satisfying the assumptions is norm dense 
in the corresponding Banach space. Indeed, the subspace contains  
$C_c(G)\ast L^\infty(X)$ which is clearly norm dense in $L^p(X)$. 
\end{remark} 




Our spectral approach uses direct integral decomposition 
for the representation of $G$ in $L^2(X)$, and we 
thus assume that $G$ is a group of type I. As is well-known, this  
assumption satisfied by all $S$-algebraic groups. We note further 
that typically, for an Iwasawa group $G=KP$, $K$ is large in $G$, 
namely in every {\it irreducible} representation $\pi$ of $G$, 
the space 
of $(K,\tau)$-isotypic vectors in $\cH_\pi$ is 
finite dimensional for every irreducible 
representation $\tau$ of $K$. Again, this property holds for 
every $S$-algebraic group.

\subsection{Ergodic theorems in the absence of a spectral gap, III}

We can now state the following convergence 
theorem for admissible families
of averages.

\begin{theorem}\label{acceptable}
{\bf Pointwise 
convergence on a dense subspace for admissible families.}
Let $G$ be an lcsc group of type I, $K$ a compact subgroup. 
Let $\nu_t$ be an admissible family of averages on $G$. 

Consider a $K$-finite vector $u$ in the $\tau$-isotypic component under
$K$, in an irreducible infinite-dimensional unitary
representation $\pi$ of $G$ (we do not assume $\tau$ is irreducible). 
Assume that for every such $\pi$  
there exists $\delta=\delta_\pi > 0$ and a positive constant 
 $C_\pi(\tau,\delta)$ (both independent of $u$), such that   

$$
 \int_{t_0}^\infty e^{\delta r} \left(\norm{\pi(\nu_r)u}+
\norm{\pi(\partial\nu_r)u}\right)^2 dr \le C_\pi(\tau,\delta)
\norm{u}^2\,\,.$$

Then in any measure-preserving weak-mixing action of 
$G$ on $(X,\mu)$ 
there exist closed subspaces 
$\cH_{\tau,\delta}\subset L^2_0(X)$ where 
$\pi(\nu_t)f(x)\to 0$  
almost surely for
$f\in \cH_{\tau,\delta}$. The convergence 
is of course also in the $L^2$-norm.  
Furthermore the union 
$$\cup_{\delta > 0,\tau\in
 \widehat{K}} \cH_{\tau,\delta}$$
 is dense in $L^2_0(X)$.  

\end{theorem}


\begin{proof} 
Our proof of Theorem \ref{acceptable} is divided into two parts, 
as follows.

1) First part of proof  : Direct integrals.

Any unitary representation $\pi$ is 
of the form $\pi=\int^\oplus_{z\in
\Sigma_\pi} \pi_z dE(z)$, where $\Sigma_\pi\subset \widehat{G}$ 
is the spectrum of the representation, and $E$ the corresponding
(projection valued) measure. Furthermore the Hilbert space of the 
representation admits a direct integral decomposition 
$\cH_\pi=\int^\oplus_{z\in
\Sigma_\pi} \cH_z dE(z)$. In particular, any vector $u\in
\cH$ can be identified with a measurable section of the family 
$\set{\cH_z\,;\, z\in \Sigma_\pi}$, namely 
$u=\int^\oplus_{z\in \Sigma_\pi} u_z dE(z)$, where $u_z\in \cH_z$ 
for $E$-almost all $z\in \Sigma_\pi$. Clearly, $u$ belongs to the
$\tau$-isotypic component of $\pi$ if and only if $\pi_z$ belongs 
to the $\tau$-isotypic of $\pi_z$ for $E$-almost all $z$. 
To see this note that $u$ is characterized by the equation
$u=\pi(\chi_{\tau})u$, where $\chi_{\tau}$ 
is the character of the
representation $\tau$ on $K$, 
namely as being in the range of a self-adjoint projection operator. The 
 projection operator commutes with all
the spectral projections, since the latter 
commute with all the unitary operator
$\pi(g)$, $g\in G$ and hence also with their linear combinations.
 This
implies $\pi(\chi_{\tau})u_z=u_z$ $E$-almost surely. Given another
vector $v$ in the $\tau$-isotypic component we conclude that the
following spectral representation is valid : 
$$\inn{\pi(g)u,v}=\int_{z\in \Sigma_\pi}
 \inn{\pi_z(g)u_z,v_z}dE_{u,v}(z)$$
where $E_{u,v}$ is the associated (scalar) spectral measure. 
Thus the $K$-finite vectors of $\pi$ (with variance
$\tau$ under $K$) are integrals
(w.r.t. the spectral measure), of $K$-finite vectors  
with the same variance, associated with irreducible unitary
representations $\pi_z$ of $G$. 

2) Second part of proof : Sobolev space argument. 

Our second step is a Sobolev space argument, 
following \cite[\S 7.1]{N1}. 
We assume the estimate stated
in Theorem \ref{acceptable} for every irreducible
non-trivial representation $\pi$. Given $\tau\in \widehat{K}$, for each $\pi$ we define
$\delta_\pi$ to be one half of the supremum of all  $\delta$
that satisfy the estimate stated in Theorem \ref{acceptable}, 
with some finite constant $C_\pi(\tau,\delta)$, 
(for all $u$ in the $\tau$-isotypic component).

We note that this function of $\pi$ is measurable
w.r.t. the spectral measure, and therefore for this choice of 
$\delta=\delta_\pi$, the estimator (for $u$ in the 
$\tau$-isotypic component)  
$$C_\pi(\tau,\delta_\pi)=2\sup_{\norm{u}=1} 
\int_{t_0}^\infty e^{r\delta_\pi } 
\left(\norm{\pi(\nu_r)u}+\norm{\pi(\partial \nu_r)u}\right)^{2}dr$$
is also measurable w.r.t. the spectral measure. 
Therefore we can consider the measurable sets 
$$ A(\tau,\delta,N)
=\set{z\in \Sigma_\pi\,;\,
\delta_{\pi_z} > \delta\,\,,\,\, C_{\pi_z}
(\tau,\delta_{\pi_z})\le N}$$
and the corresponding closed spectral subspaces 
$$\cH(\tau,{\delta},N)=\int^\oplus_{A(\tau,\delta,N)}
 \cH_z dE(z)\,\,.$$
 Thus in particular in these subspaces 
the decay of $\tau$-isotypic 
$K$-finite matrix coefficients is exponentially fast, 
with at least a fixed
positive rate, determined by $\delta$.

Now note that the subspace of differentiable vectors 
in $A(\tau,\delta,N)$ which are invariant under 
the projection $\pi(\chi_\tau)$ is norm 
dense in the $\tau$-isotypic subspace. Indeed, the subspace 
$\pi(\chi_{\tau}\ast C_c(G))( A(\tau,\delta,N))$ consists 
of differentiable $\tau$-isotypic vectors and is dense in the $\tau$-isotypic subspace. 
Furthermore, given a differentiable
 vector $u$ in the $\tau$-isotypic subspace, 
$u_z$ is also differentiable,
 for $E$-almost all $z\in \Sigma_\pi$ 
(w.r.t. the spectral measure) since for $f\in C_c(G)$ 
$$\pi(\chi_{\tau}\ast f\ast \chi_\tau)u
=\int_{z\in A(\tau,\delta,N)}
 \pi_z(\chi_{\tau}\ast f\ast \chi_\tau)u_z dE(z)\,.$$

We can now use the first part of  
Corollary \ref{diff vectors}, together with standard spectral 
theory,  and conclude that 
for a differentiable vector $u\in A(\tau,\delta,N)$,  
the following spectral identity holds, 
for every $t > s > 0$, and for every  $u,v\in L^2(X)$ :

$$\inn{\left(\pi(\nu_t)-\pi(\nu_s)\right)u,v}=
\int_{z\in \Sigma_\pi}
\inn{\left(\pi_z(\nu_t)
-\pi_z(\nu_s)\right)u_z,v_z}dE_{u,v}(z)=$$
$$=\int_{z\in \Sigma_\pi}
\int_s^{t}\inn{\frac{d}{dr}\pi_z(\nu_r) u_z,v_z} dr\, dE_{u,v}(z)\,\,.$$

%
%
%

Now using Corollary \ref{diff vectors}(2) and the fact that 
$v$ above is allowed to range over $L^2(X)$, for each $t $ and $s$
we have the following equality of functions in $L^2(X)$, namely for
almost all $x\in X$ :

$$\pi(\nu_t)u(x)-\pi(\nu_s)u(x)=
\int_s^t \frac{d}{dr}\pi( \nu_r)u(x) dr $$
so that for any $t> s \ge M $, for almost all $x\in X$ : 
$$\abs{\pi(\nu_t)u(x)-\pi(\nu_s)u(x)}\le 
\int_M^\infty \abs{\frac{d}{dr}\pi( \nu_r) u(x)} dr \,\,.$$
The averages $\nu_t$ form an continuous family in the $L^1(G)$-norm, 
consisting of absolutely continuous measures on $G$, 
and the function $t\mapsto \pi(\nu_t)u(x)$ therefore 
is a continuous function 
of $t$ for almost every $x\in X$. Restricting attention to these
points $x$, we conclude that for all $M > 0$ and almost every $x$ 
$$\limsup_{t, s\to \infty}\abs{\pi(\nu_t)u(x)-\pi(\nu_s)u(x)}\le  
 \int_M^\infty \abs{\frac{d}{dr}\pi( \nu_r) u(x)} dr $$
and thus the set 
$$\set{x \,;\,\limsup_{t, s\to
\infty}\abs{\pi(\nu_t)u(x)-\pi(\nu_s)u(x)} > \zeta }$$
 is contained in the set 
 $$\set{x\,;\, \int_M^\infty 
\abs{\frac{d}{dr}\pi( \nu_r) u(x)} dr > \zeta}\,. $$

We estimate the measure of the latter set by integrating over $X$, 
and using the Cauchy-Schwartz inequality.  
We obtain, for any $\zeta > 0$
and $M > 0$, the following 
estimate: 
 $$\mu \set{x\,;\, \int_M^\infty 
\abs{\frac{d}{dr}\pi( \nu_r) u(x)} dr > \zeta} $$
$$\le \frac{1}{\zeta}\int_X\left( \int_M^\infty 
\abs{\frac{d}{dr}\pi( \nu_r) u(x)} dr 
\right)d\mu(x)\le 
 \frac{1}{\zeta}\int_M^\infty
\norm{\frac{d}{dr}\pi(\nu_r)u}_{L^2(X)}dr $$
$$\le \frac{\exp(-M\delta/4)}{\zeta}
\int_M^\infty 
e^{-r\delta/4 }
e^{r\delta/2 }
\norm{\frac{d}{dr}\pi(\nu_r)u}_{L^2(X)}dr  
$$
$$
\le \frac{2\exp(-M\delta/2)}{\zeta \sqrt{M\delta/2}}
\left(\int_M^\infty e^{r\delta }
\norm{\frac{d}{dr}\pi(\nu_r)u}^2_{L^2(X)}dr\right)^{1/2}\,\,.
$$

Using Proposition \ref{derivative},  it suffices to show the 
finiteness of the following 
expression 
$$
\int_M^\infty e^{r\delta }
\norm{\frac{m_r(S_r)}{m_G(N_r)}
\left(\partial\nu_r -\nu_r\right)u}^2_{L^2(X)} dr\,\,. $$
Using 
our assumption that the $K$-finite vector 
$u$ has its spectral support in the set $A(\tau,\delta,N)$, 
we can write the last expression as 
$$
\int_M^\infty e^{\delta r} \int_{z\in
A(\tau,\delta,N)}
 \norm{\frac{m_r(S_r)}{m_G(N_r)}\pi_z(\nu_r-\partial
\nu_r)u_z}_{\cH_z}^2 dE_{u,u}(z)dr\,\,.   
$$
Using the uniform bound given in Corollary \ref{unif bound}, 
we can estimate by

$$
\le \int_M^\infty e^{\delta r} \int_{z\in
A(\tau,\delta,N)}
c^2\left(\norm{\pi_z(\nu_r)u_z}+
\norm{\pi_z(\partial\nu_r)u_z}\right)^2
dE_{u,u}(z)dr   
$$
and thus by definition of the space $A(\tau,\delta,N)$, and the 
fact that $u$ is spectrally supported in this subspace, the last 
expression is bounded by 
$c^2 N \norm{u}^2 < \infty$.

We have established that $\pi(\nu_t)u(x)$ 
converges almost surely (exponentially fast) to the ergodic mean, namely to zero. 
The fact that 
$$\cup_{\delta > 0,\tau\in \widehat{K}} \cH_{\tau,\delta}$$
 is dense in $L^2_0(X)$ is a standard fact in spectral theory. 

This concludes the proof of Theorem \ref{acceptable}.
\end{proof}

To complete the proof of Theorem
 \ref{th:semisimple-1}, we now need to verify that the 
assumptions of Theorem \ref{acceptable} are satisfied by 
$S$-algebraic groups. We begin with following 

\begin{theorem}\label{partial estimate} 
Let $G=G_1\cdots G_N$ be an $S$-algebraic group as in Definition 
\ref{alg-gps}. Let $G_r$ be any  
family of bounded Borel sets,
 and $\nu_r$ the Haar-uniform probability measures. 
Let $\pi=\pi_1\otimes\cdots \otimes \pi_N$, where  each $\pi_i$ is an irreducible 
unitary representation of $G_i$ without $G_i^+$-invariant unit vectors. Then 
there exist $\delta=\delta_\pi > 0$ and a constant $C_1$  
(depending only on $G$ and the family $G_r$) such that
 for every $\tau$-isotypic 
vector $u$ 
   
\begin{enumerate}
\item When $G_r$ is a coarsely admissible
 $1$-parameter family (or sequence), we have  
$$\norm{\pi(\nu_r)u}\le C_1(\dim \tau)
 e^{-\delta r}\norm{u}\,\,,$$ 
\item When $G_r$ is a left-$K$-radial admissible family, $K$ a good maximal compact subgroup, for almost every $r$ we have 
$$\norm{\pi(\partial\nu_r)u}\le C_1(\dim \tau) e^{-\delta r}
\norm{u}\,\,,$$ 
\item In particular, when $G_r$ is left-radial and 
admissible, there exists a constant 
$C_\pi(\tau,\delta) < \infty $ such that 
$$\int_{t_0}^\infty e^{\delta r}
\left(\norm{\pi(\nu_r)u}
+\norm{\pi(\partial \nu_r)u}\right)^2 dr 
\le C_\pi(\tau,\delta)\norm{u}^2\,.$$
\end{enumerate}

\end{theorem}
\begin{proof}

1) By Theorem 
\ref{Lp representations}, $\pi$ is strongly $L^p$ 
for some $p=p(\pi)< \infty$, and then if $n \ge p/2$ then 
$\pi^{\otimes n}\subset \infty \cdot \lambda_G$. Assume 
without loss of generality that $n$ is even, and then 
 (see \cite[Thm. 1.1]{N3}), since $\inn{\pi(g)u,u}$ 
is real-valued, using 
Jensen's inequality we obtain
$$\norm{\pi(\nu_r)u}^{2n}=
\left(\int_G \inn{\pi(g)u,u}d(\nu_r^\ast\ast\nu_r)\right)^n
$$
$$\le \int_G \left(\inn{\pi(g)u,u}\right)^n d(\nu_r^\ast\ast\nu_r)=
\int_G\inn{\pi^{\otimes n}(g)u^{\otimes n},u^{\otimes n}}
d(\nu_r^\ast\ast\nu_r)
$$
$$\le  \dim(\tau)^n \norm{u}^{2n} \int_G \Xi(g)
d(\nu_r^\ast\ast\nu_r) 
$$ 
where we have used the estimate given in 
Theorem \ref{Kfinite estimate}(1) for $K$-finite matrix  
coefficients in representations (weakly) 
contained in the regular representation. 
Now $\Xi$ is non-negative, and 
$G_r$ are assumed coarsely admissible and hence 
$(K,C)$-radial. Thus we can multiply 
the last estimate by $C^2$ and then replace $\nu_r$  
by their radializations $\tilde{\nu}_r$. 
 
Since we are considering 
$S$-algebraic groups, we can assume without 
loss of generality 
that $K$ is a good maximal compact subgroup so 
that $(G,K)$ is a 
Gelfand pair. Then
 $\Xi$ defines a
homomorphism of the commutative convolution 
algebra of bi-$K$-invariant 
functions $L^1(K\setminus G/K)$. We can therefore conclude that 
$$\norm{\pi(\nu_r)u}^{2n}
\le C^2 (\dim \tau)^n \norm{u}^{2n}
\left( \int_G\Xi(g)d\tilde{\nu}_r(g)
\right)^2\,\,.$$ 
Coarse admissibility implies the property of minimal 
growth for $G_r$ and their radializations 
$\tilde{G}_r$, namely 
$S^n\subset G_{an+b}$, for a compact generating set $S$. 
Thus 
the desired result follow from the 
standard estimates of the 
$\Xi$-function of an $S$-algebraic group, which shows that the integral 
of $\Xi_G$ on $G_r$ decay exponentially in $r$.

2) Now consider the case of $\partial \nu_r$, which 
is a singular measure on $G$, supported on the ``sphere'' 
$S_r$. Arguing as in (1) before 
we still have   
$$\norm{\pi(\partial\nu_r)u}^{2n}\le  
 \dim(\tau)^n \norm{u}^{2n}
 \int_G \Xi(g)
d(\partial\nu_r^\ast\ast\partial\nu_r)
$$ 
Now $\Xi$ is bi-$K$-invariant for a good maximal compact 
subgroup $K$, and $G_r$ (and thus $\nu_r$)
 are assumed to be 
left-$K$-invariant. In follows that $m_K\ast \partial \nu_r=\partial \nu_r$  for almost every $r$, and hence 
$$ \int_G \Xi(g)
d(\partial\nu_r^\ast\ast\partial\nu_r)=
 \int_G m_K\ast \Xi \ast m_K(g)d
(\partial\nu_r^\ast\ast m_K\ast \partial\nu_r)$$
$$=
\left(\int_G\Xi(g)d\partial\nu_r\right)^2
$$ 
But $\partial \nu_r=m_r/m_r(G_r)$ is a 
probability measure 
supported on $S_r$, and clearly 
the property of minimal growth 
for $G_r$ implies that $S_r$ is 
contained in the complement of 
$S^{a_1[r]+b_1}$ for some $a_1 > 0$. Therefore again the 
standard estimates of the $\Xi$-function yield the 
desired result.

3) The last part is an immediate consequence 
of the previous two.

\end{proof}


{\it Completion of the proof of Theorem \ref{th:semisimple-1}.}

The last step in the proof of Theorem \ref{th:semisimple-1} is 
to consider the various alternatives stated in its assumptions. 

If the action is irreducible and totally weak mixing, 
then any irreducible unitary representation $\pi_z$ of $G$ 
appearing in the direct integral decomposition of $\pi_X^0$ 
is indeed strongly 
$L^p$ for some finite $p$. This follows from Theorem 
\ref{Lp representations} since $\pi_z$  is then 
a tensor product of 
infinite dimensional irreducible 
representations of the simple constituent groups. In that case   
Theorem \ref{partial estimate}, parts (1) and (2)   
apply, and the proof of the mean and the pointwise 
ergodic theorem for left-radial admissible 
$1$-parameter families in irreducible actions is complete, 
taking into account also that 
the maximal inequality is covered in all cases by Theorem 
\ref{ample sets}. 
 
Note that, still in the irreducible case, we can apply 
part (1) 
of  Theorem \ref{partial estimate} to a coarsely 
admissible sequence, and this immediately 
yields the mean ergodic theorem and pointwise convergence 
almost surely on the dense subspace of vectors appearing there. 
Again using Theorem \ref{ample sets}, 
this completes the proof of the mean and 
pointwise ergodic theorem for coarsely admissible sequences in  
irreducible actions. 

Otherwise the action may be reducible, and we seek 
to prove the mean theorem when the left-radial  
averages are balanced and 
the pointwise theorem when they are standard radial and well-balanced. 
In the present
 case, each $\pi_z$ is a tensor product of infinite 
dimensional irreducible representations of some of the simple 
subgroups, and 
the trivial representations of the others. We can repeat the 
argument used in the first part of the proof of Theorem
 \ref{partial estimate}, and establish that 
 $\norm{\pi_z(\beta_t)u}\to 0$ 
using the assumption 
that $\beta_t$ and hence $\tilde{\beta}_t$ are balanced, and 
 $\norm{\pi_z(\beta_t)u}\to 0$ exponentially fast when 
$\beta_t$ are well-balanced.
 Indeed, instead of 
integrating against $\Xi_G(g)$ 
we will now be integrating against 
the $\Xi$-function lifted from some 
simple factor group, using the argument in 
 the second part of the proof of 
Theorem \ref{Exp decay}. By the balanced or 
well-balanced assumption, 
the standard estimates of the $\Xi$-function 
yield the desired norm decay conclusion.

 The last argument required to complete the proof of the pointwise theorem is the estimate 
of $\norm{\pi(\partial \beta_t)}$, when the averages are standard radial, well-balanced and boundary-regular.
 In this case each distance $\ell$ (or $d)$ on a factor group $L$ 
 obeys the estimate  provided by Theorem \ref{CAT}(ii), namely 
 $m_t(\partial G_t\cap L_{\alpha  t})\le Ce^{-\beta t} m_t(\partial G_t)$. 
 The total measure $m_t(S_t)$ 
on  $\partial G_t\subset G$ is obtained as an interated integral over the factor groups. Integrating against the $\Xi$-function lifted from a factor group, and using the decay of the $\Xi$-function, the required estimate follows. 

This concludes the proof of all parts of Theorem 
\ref{th:semisimple-1} (and of course also Theorem \ref{th: Mean Lie groups}).

%

\begin{remark}
\begin{enumerate}
\item 
In principle, our analysis  applies to a general almost surely 
differentiable family of averages 
$\nu_t$ (absolutely continuous w.r.t. Haar measure), and 
not only those arising from 
Haar uniform averages on admissible sets $G_t$ as 
in Theorem \ref{acceptable}. 
\item We need only assume that 
the irreducible representations of $G$ giving rise to the spectral decomposition of $L^2(X)$
satisfy the spectral estimates we have employed, and not necessarily all representations of $G$.
 This is useful  when considering a homogeneous space $X=G/\Gamma$, when $G$ is 
 an adele group, for example.  
\end{enumerate}
\end{remark}

\begin{remark}{\bf Singular averages.}
An important problem that arises naturally here is to extend the foregoing 
analysis to averages which are singular w.r.t. Haar measure.  
An obvious first step would be to establish a pointwise ergodic theorem for  the family of ``spherical averages'' 
supported on the boundaries $\partial G_t$ of the sets $G_t$. 
However, to prove such results it is necessary
 to establish estimates for the 
{\it derivatives} of the $\tau$-spherical functions. 
 While the 
matrix coefficients 
themselves obey uniform decay estimates which are independent of the 
representation (provided, say, that it is $L^p$, see
Theorems \ref{tensor} and \ref{Kfinite estimate}), this
is no longer the case for their derivatives. For example, consider the
principal series representations $\Ind_{MAN}^G 1\otimes i\eta$ induced
from a unitary character of $A$ and the trivial representation of 
$MN$. These representations 
have matrix coefficients whose derivatives
exhibit explicit dependence on the character $\eta$ parametrizing 
the representation. Consequently, sufficiently sharp 
{\it derivative} estimates for matrix coefficients  
are inextricably tied up with classification, or at least
parametrization, of the irreducible 
unitary representations of the
group (see \cite{CN} for more on this point). 

We have avoided appealing to classification theory and refrained 
from establishing 
such derivative estimates in the present paper. 
Instead we have utilized the fact 
that restricting to {\it Haar uniform averages}  
on admissible sets, the distribution $\frac{d}{dt}\nu_t$ is a signed  
measure, so that we need only use estimates of 
the spherical functions themselves in order to estimate it.

\end{remark}

\subsection{The invariance principle, and stability 
of admissible averages}

\subsubsection{The set of convergence}

It will be essential in our argument below to establish that 
for a family of admissible averages, the set where pointwise 
convergence of $\pi(\beta_t)f(x)$ holds 
contains a $G$-invariant set, for each 
fixed function $f$.

Let $G$ be a locally compact second countable group with left Haar measure $m_G$.
Consider a measure-preserving action of $G$ on a standard Borel 
space $(X,\mathcal{B},\nu)$.
For Borel subsets $G_t\subset G$ and $g\in G$, consider probability measures
$$
\beta_t^g=\frac{1}{m_G(G_t)}\int_{gG_t} \delta_h\, dm_G(h)\quad\hbox{and}\quad \beta_t=\beta_t^e.
$$

Let us formulate the following 
invariance principle which applies to 
all quasi-uniform  families. 
This result 
generalizes  \cite{BR}, where the case 
of ball averages on $SO(n,1)$ was considered.

\begin{theorem}\label{th:g_inv}
Let $G$ be an lcsc group, and 
suppose that $\{G_t\}_{t>0}$ is a quasi-uniform family, with 
$\beta_t$  satisfying 
the pointwise ergodic theorem in $L^p(\nu)$.
Then for every $f\in L^p(\nu)$, there exists a
 $G$-invariant measurable set $\Omega(f)$ of full measure 
such that for every $x\in \Omega(f)$,
$$
\lim_{t\to\infty} \pi(\beta_t)f(x)=\int_X f\, d\nu. 
$$
In particular, this holds for admissible 
$1$-parameter families and 
admissible sequences on $S$-algebraic groups.  
\end{theorem}

\begin{proof}
Writing $f=f^+-f^-$ for $f^+,f^-\in L^p(\nu)$, $f^+,f^-\ge 0$, and assuming that the theorem holds for $f^+$ and $f^-$,
we can take
$$
\Omega(f)=\Omega(f^+)\cap\Omega(f^-).
$$
Hence, without loss of generality, we may assume that $f\ge 0$.

Consider then the conull measurable set of convergence: 
$$
C=\left\{x\in  X:\, \lim_{t\to\infty} \pi(\beta_t)f(x)=\int_X f\, d\nu \right\}.
$$
Take a countable dense set 
$\{g_i\}_{i\ge 1}\subset G$ and let 
$$
\Omega=\bigcap_{i\ge 1}g_i C.
$$
Then $\Omega$ is a measurable  set of full measure, and for every 
$x\in \Omega$ and every $g_i$, we have $g_i^{-1}x\in C$. 
Let $\delta>0$ and take $\vre>0$ and $\mathcal{O}$
 as in (\ref{eq:01})
and (\ref{eq:02}). We may also assume that 
$\mathcal{O}$ is symmetric.
Then for any $g\in G$ there 
exists $g_i$ such that $g_i\in g\mathcal{O}$.
Hence, for sufficiently large $t$,
$$
g_i G_{t-\vre}\subset g G_t\subset g_i G_{t+\vre}.
$$
Therefore, for every $x\in X$,
\begin{align*}
&\pi(\beta_t)f(g^{-1}x)
=\frac{1}{m_G(G_t)}\int_{G_t}f(h^{-1}g^{-1}x)dm_G(h)\\
&=\frac{1}{m_G(G_t)}\int_{gG_t}f(u^{-1}x)dm_G(u)\le \frac{1}{m_G(G_t)}\int_{g_iG_{t+\vre}}f(u^{-1}x)dm_G(u)\\
&\le  \frac{1+\delta}{m_G(g_iG_{t+\vre})}
\int_{g_iG_{t+\vre}}f(u^{-1}x)dm_G(u)
=(1+\delta)\pi(\beta_{t+\vre})f(g_i^{-1}x).
\end{align*}
This implies that for every $g\in G$ and $x\in\Omega$, since 
$g_i^{-1}x\in C$
$$
\limsup_{t\to\infty}\pi(\beta_t)f(g^{-1}x)\le (1+\delta)\int_X f\, d\mu
$$
for every $\delta>0$. Similarly, we show that
$$
\liminf_{t\to\infty}\pi(\beta_t)f(g^{-1}x)\ge (1+\delta)^{-1}\int_X f\, d\mu.
$$
Therefore, let us take $\Omega(f)=G\cdot\Omega$. Then 
$\Omega(f)$ is strictly invariant under $G$, namely 
$g\Omega(f)=\Omega(f)$ for every $g\in G$, and the complement of 
$\Omega(f)$ is a null set. Thus $\Omega(f)$ is a strictly invariant 
measurable set in the Lebesgue $\sigma$-algebra, namely in the completion 
of the standard Borel structure on $X$ with respect to the measure $\mu$.

\end{proof}

An immediate corollary of the foregoing considerations
 is the following 
\begin{corollary}\label{gGth}
Let $G$ be an lcsc group, and 
suppose that $\{G_t\}_{t>0}$ is a quasi-uniform
family, with  
$\beta_t$  satisfying 
the pointwise ergodic theorem in $L^p(\nu)$. Then the Haar-uniform 
averages on $gG_th$ also satisfy it, for any fixed $g,h\in G$.
In particular, this holds for admissible 
$1$-parameter families and 
admissible sequences on $S$-algebraic groups.   
\end{corollary}

\subsubsection{Stability of admissible averages 
under translations}

When is the family $gG_t h$  
itself already admissible if $G_t$ is ? This was asserted 
in Definition \ref{admissible-Lie} for connected Lie groups. 
In this subsection we note that in the general case, 
the property of admissibility is stable under 
two-sided translations. Indeed,  the sets $\mathcal{O}_\vre$ 
we used to defined admissibility on $S$-algebraic groups satisfy the following. 
For every $g\in G$, there exists a positive constant $c(g)> 0$ 
such that $g\mathcal{O}_\vre g^{-1}$ contains $\mathcal{O}_{c(g)\vre}$ 
for all $0 < \vre < \vre_0$.  We can therefore easily conclude :

\begin{lemma}\label{stability}
{\bf Stability under translations.} Let $G_t$ be a
 $1$-parameter 
family of coarsely admissible averages on an $S$-algebraic group as in Definition \ref{alg-gps}. 
Then for any 
$g,h\in G$, the family 
$gG_t h$ is also coarsely admissible. If $G_t$ is admissible, then so is $gG_t h$. 
\end{lemma}
\begin{proof}

To see that for coarsely admissible averages $G_t$, 
the averages $gG_t h$ are also coarsely admissible
note that 
for any bounded set $B$, 
$$BgG_t hB \subset B^\prime G_t B^\prime  \subset G_{t+c^\prime}$$
and in addition 
$g^{-1}G_{t+c^\prime}h^{-1}\subset G_{t+c''}$, so that 
$$BgG_t hB \subset G_{t+c^\prime}\subset gG_{t+c''}h\,.$$
As to the second condition of coarse admissibility, 
by unimodularity, 
$m_G(gG_{t+c} h)\le dm_G(gG_t h)$ and so 
 $gG_t h$ is coarsely admissible. 


Now let $G_t$ be an admissible 
$1$-parameter family, and $h,g$ be fixed. 
For every open set $\mathcal{O}_\vre$ in the basis, the set 
$g\mathcal{O}_\vre g^{-1}\cap h^{-1}\mathcal{O}_\vre h$ 
is open and contains 
$\mathcal{O}_{\eta(\vre)}$. 
By definition of an appropriate basis, we can choose 
$\eta(\vre)\ge c_0\vre$, for some fixed 
positive $c_0=c_0(g,h)< 1 $, 
uniformly for all $0< \vre <  \vre(g,h)$. 
 Then, checking the 
conditions in the definition of admissibility :
$$\mathcal{O}_{\eta(\vre)} gG_th \mathcal{O}_{\eta(\vre)}\subset 
g\mathcal{O}_\vre G_t \mathcal{O}_\vre h\subset gG_{t+c\vre} h$$
so that for all $0< \vre < \vre^\prime(g,h)$  
$$\mathcal{O}_\vre gG_t h \mathcal{O}_\vre \subset gG_{t+\vre c/c_0}h\,.$$

When $G$ is totally disconnected, and $G_t$ satisfies 
$K G_t K=G_t$,  clearly $gG_th$ 
is also invariant under translation by the 
compact open subgroup $K^\prime=gKg^{-1}\cap hKh^{-1}$.

As to the Lipschitz continuity of 
the measure of the family, 
we have of course, since $G$ is unimodular
$$m_G(gG_{t+\vre} h)=m_G(G_{t+\vre})\le (1+c\vre) m_G(G_t)=$$
$$=
(1+c\vre)m_G(gG_t h)\,.$$ 
\end{proof}

\section{Proof of ergodic theorems for lattice actions}

\subsection{Induced action}\label{sec:induced}

We now turn to consider an lcsc group $G$ 
 and a discrete lattice $\Gamma$ in $G$.
The existence  of a lattice implies that $G$ is unimodular, and 
we denote Haar measure by $m_G$, as before. 
Denote by $m_{G/\Gamma}$ the corresponding measure on $G/\Gamma$.
We normalize $m_G$ so that $m_{G/\Gamma}(G/\Gamma)=1$.

For a family of Borel subsets $\{G_t\}_{t>0}$, 
we consider the averages 
$\lambda_t$ uniformly distributed on $G_t\cap \Gamma$. We will  
use the mean, maximal and pointwise ergodic theorems  
established for 
the averages $\beta_t$ acting in a $G$-action, in order 
to establish similar ergodic theorems for the averages 
$\lambda_t$ acting in a $\Gamma$-action. The fundamental link 
used to implement this reduction is of course 
the well-known construction of the
 {\it induced $G$-action} defined for 
a measure-preserving action of $\Gamma$, to which we now turn.

Thus let $\Gamma$ act on a standard Borel space $(X,\mathcal{B},\mu)$, 
preserving the probability measure $\mu$. Let
$$
\tilde Y\stackrel{def}{=}G\times X.
$$
Define the right action of $\Gamma$ on $\tilde Y$:
\begin{equation}\label{eq:G_act}
(g,x)\cdot \gamma=(g\gamma,\gamma^{-1} x)\quad (g,x)\in \tilde Y,\; \gamma\in \Gamma,
\end{equation}
and the left action of $G$:
\begin{equation}\label{eq:Gamma_act}
g_1\cdot (g,x)=(g_1g,x), \quad (g,x)\in \tilde Y,\; g_1\in G.
\end{equation}
The space $\tilde Y$ is equipped with the 
product measure $m_G\otimes \mu$, which is preserved by these actions.
Since the actions (\ref{eq:G_act}) and (\ref{eq:Gamma_act}) commutes,
there is a well-defined action of $G$ on the factor-space
$$
Y\stackrel{def}{=}\tilde Y/\Gamma.
$$
We denote by $\pi$ the projection map $\pi : \tilde Y\to Y$. Note 
that $Y$ admits a natural map $j: (g,x)\Gamma \mapsto g\Gamma$ 
onto $G/\Gamma$. This map is measurable and $G$-equivariant, and 
thus $Y$ is a bundle over the homogeneous space $G/\Gamma$, with 
the fiber over each point $g\Gamma$ identified with $X$. 

For a bounded measurable function 
$\chi: G\to \mathbb{R}$ with compact support and a measurable
function $\phi:X\to \mathbb{R}$, we define $\tilde{F}: 
\tilde{Y}\to \RR$
by $\tilde{F}(g,x)=\chi(g)\phi(x)$. We then define 
$F: Y\to \mathbb{R}$ by summing over $\Gamma$-orbits  
\begin{equation}\label{eq:F}
F(y)=F((g,x)\Gamma)=\sum_{\gamma\in\Gamma} \chi(g\gamma)\phi(\gamma^{-1} x)=
\sum_{\gamma\in \Gamma}\tilde{F}((g,x)\gamma)\,\,.
\end{equation}
There is a unique $G$-invariant Borel measure $\nu$ on $Y$ such that
\begin{equation}\label{eq:F_int}
\int_Y F\,d\nu=\left(\int_G\chi\, dm_G\right) \left(\int_X \phi\, d\mu\right).
\end{equation}

For $F$ defined above, we have the following expression for the averaging 
operators we will consider below. Let $(h,x)\Gamma=y\in Y$, and then  
$$\pi_Y(\beta_t)F(y)=\frac{1}
{m_G(G_t)}\int_{G} F(g^{-1}y)d\beta_t(g)$$
$$=\sum_{\gamma\in \Gamma}  \frac{1}{m_G(G_t)} \left
(\int_{G_t}\chi(g^{-1}h\gamma)dm_G(g)\right)\phi(\gamma^{-1}x)
$$
The latter expression will serve as the basic link between the averaging 
operators $\beta_t$ on $G$ acting on $L^p(Y)$, and the averaging operators 
$\lambda_t$ acting on $L^p(X)$.

We now recall the following fact regarding 
induced actions, which will  
play an important role below. Namely, it will allow us 
to deduce results 
about the pointwise behaviour of the averages $\lambda_t$ on the 
$\Gamma$-orbits 
in $X$ from the pointwise behaviour of the averages $\beta_t$ 
on $G$-orbits in $Y$. 

Consider the factor map $j: (Y,\nu)\to (G/\Gamma,m_{G/\Gamma})$, which is 
a Borel measurable, everywhere defined, $G$-equivariant and measure-preserving. 
For a Lebesgue measurable set $B\subset Y$, the set $B_{y\Gamma}=
j^{-1}(y\Gamma)\cap B$ is a Lebesgue measurable subset of $X$ for every 
$y\Gamma\in G/\Gamma$. (Recall that the Lebesgue $\sigma$-algebra is the completion of the Borel $\sigma$-algebra w.r.t. the measure at hand, namely $\nu$ on $Y$ or $\mu$ on $X$).  

Any set $B$ can be written as the disjoint 
union $B=\coprod_{y\Gamma\in G/\Gamma}B_{y\Gamma}$. 
Furthermore, the $G$-action is given by 
$$gB=\coprod_{y\Gamma\in G/\Gamma} \alpha(g,y\Gamma)B_{y\Gamma}$$
where $\alpha : G\times G/\Gamma\to \Gamma$ is a 
Borel cocycle associated with 
a Borel section of the canonical projection $G\to G/\Gamma$.  

We can now state the 
following well-known fact, 
 whose proof is included for completeness.  
\begin{lemma}\label{strict}
If $B\subset Y$ is a  Lebesgue measurable set with $\nu(B)=1$, which is strictly 
$G$-invariant ($gB=B$ for all $g\in G$)  
then $\mu(B_{y\Gamma})=1$, for {\bf every} $y\Gamma\in G/\Gamma$ 
(and not only for almost every $y\Gamma$). 
\end{lemma}
\begin{proof} The map $b : G/\Gamma\to \RR_+$ given by 
$y\Gamma\mapsto \mu(B_{y\Gamma})$ is everywhere defined, 
Lebesgue measurable, 
and strictly $G$-invariant, namely $b(gy\Gamma)=b(y\Gamma)$ for all 
$g\in G$ and $y\Gamma\in G/\Gamma$. 
Since $G/\Gamma$ is a transitive $G$-space,  
$b(y\Gamma)$ is strictly a constant, and 
this constant is of course $1$.  
\end{proof} 

We conclude the introduction to induced actions with following 
simple fact. 
\begin{lemma}\label{l:lp}
Let $1\le p\le \infty$ and $Q$ be a compact subset of $G$. 
\begin{enumerate}
\item[(a)] There exists $a_{p,Q}>0$ such that for every $\phi\in L^p(\mu)$ and
a bounded $\chi: G\to \mathbb{R}$ such 
that $\supp(\chi)\subset Q$, with $F$ defined 
as in (\ref{eq:F})
$$
\|F\|_{L^p(\nu)}\le a_{p,Q}\cdot \|\chi\|_{L^p(m_G)}\cdot\|\phi\|_{L^p(\mu)}.
$$
Moreover, if $Q$ is contained in a sufficiently 
small neighborhood of $e$, then
$$
\|F\|_{L^p(\nu)}=\|\chi\|_{L^p(m_G)}\cdot \|\phi\|_{L^p(\mu)}.
$$
\item[(b)] There exists $b_{p,Q}>0$ such that for any 
measurable $F: Y\to \mathbb{R}$,
$$
\|F\circ\pi\|_{L^p(m_G\otimes\mu|_{Q\times X})}\le b_{p,Q}\cdot \|F\|_{L^p(\nu)}.
$$
When $Q=\mathcal{O}_\vre$ we denote $b_{p,Q}=b_{p,\vre}$. 
\end{enumerate}
\end{lemma}

\subsection{Reduction theorems}
We now  turn to formulate the fundamental result reducing the 
ergodic theory of 
the lattice subgroup $\Gamma$ to that of the enveloping group $G$. 

Such a result necessarily involves  an approximation argument 
based on smoothing, and thus the metric properties of a shrinking family of 
neighbourhoods in $G$ come into play. The crucial property 
is finiteness of the upper local dimension of $G$ (see 
Definition \ref{ULD}), namely 
$$\varrho_0\stackrel{def}{=}
\limsup_{\vre\to 0^+}
 \frac{\log m_G(\mathcal{O}_\vre)}{\log\vre}<\infty$$.  

We will assume this condition when considering admissible sets, throughout our discussion below.  
Note that for $S$-algebraic groups as in Definition \ref{alg-gps}, and for the sets $\mathcal{O}_\vre$ we chose in that case, $\rho$ is simply the real dimension of the Archimedian factor, and thus vanishes 
for totally disconnected groups. 
  
Let us note that the induced representation  of $G$ on $L^p(Y)$, $1 \le p \le \infty$, contains the representation of $G$ on $L^p(G/\Gamma)$ as a subrepresentation. Thus, whenever a maximal 
inequality, exponential maximal inequality, norm decay estimate, spectral gap condition, mean or pointwise ergodic theorem hold for $\pi_Y(\beta_t)$ acting on $L^p(Y)$, they also hold for $\pi_{G/\Gamma}(\beta_t)$ acting on $L^p(G/\Gamma)$.

We now formulate the following reduction theorem, and 
emphasize that it is valid {\it for 
every lattice subgroup of every 
lcsc group}.

\begin{theorem}\label{reduction}{\bf Reduction Theorem.} 
Let $G$ be an lcsc group, $\mathcal{O}_\vre$ of finite upper local dimension, 
$G_t$ an increasing family of 
bounded Borel sets, 
and $\Gamma$ a lattice subgroup. 
Let $p\ge r\ge 1$, and consider the averages 
$\beta_t$ on $G_t$ and $\lambda_t$ on $\Gamma\cap G_t$ as above. 
Then 
\begin{enumerate}
\item If the family $\{G_t\}_{t>0}$ is coarsely 
admissible, then the strong maximal inequality for $\beta_t$
in $(L^p(\nu),L^r(\nu))$ implies the strong maximal 
inequality for $\lambda_t$ in $(L^p(\mu),L^r(\mu))$.
\item If the family $\{G_t\}_{t>0}$ is admissible, 
then the mean ergodic theorem for $\beta_t$
in $L^p(\nu)$ implies
 the mean ergodic theorem for $\lambda_t$ in $L^p(\mu)$.
\item If the family $\{G_t\}_{t>0}$ is quasi-uniform, 
and the pointwise ergodic theorem holds for $\beta_t$
in $L^p(\nu)$, 
then the pointwise ergodic 
theorem holds for $\lambda_t$ in $L^p(\mu)$.
\item If the family $\{G_t\}_{t>0}$ is admissible 
and $r>\varrho_0$, 
then the exponential mean ergodic theorem for $\beta_t$
in $(L^p(\nu),L^r(\nu))$ implies the 
exponential mean 
ergodic theorem for $\lambda_t$ in $(L^p(\mu),L^r(\mu))$ (but the rate may change).
\item Let the family $\{G_t\}_{t>0}$ be admissible, 
$p\ge r>\varrho_0$, and assume $\beta_t$ satisfies 
the exponential mean ergodic theorem in 
$(L^p(\nu),L^r(\nu))$, as well as 
the strong maximal 
inequality in $L^q(\nu)$, for 
$q > 1$.  
Then $\lambda_t$ satisfies the exponential 
strong maximal inequality in $(L^{p'},L^{r'})$
with $p'$, $r'$ such that $1/{p'}=(1-u)/q$ 
and $1/{r'}=(1-u)/q+u/r$ for some $u\in (0,1)$.
\end{enumerate}
\end{theorem}

The proof of Theorem \ref{reduction} will  
occupy the rest of \S 6, and will be divided to 
a sequence of separate statements.

One basic ingredient in the proof of 
Theorem \ref{reduction} is as follows.

\begin{theorem}\label{reduct-lattice-points}
Let $G$, $G_t$, $\beta_t$ and $\lambda_t$ be as in Theorem 
\ref{reduction}. Then 
\begin{enumerate}
\item Suppose that the family $\{G_t\}_{t>0}$ is coarsely 
admissible and $\beta_t$ satisfies the strong maximal
inequality in $(L^p(m_{G/\Gamma}),L^r(m_{G/\Gamma}))$ 
for some $p\ge r\ge 1$. Then for some $C>0$ and all sufficiently large $t$,
$$
C^{-1}\cdot m_G(G_t)\le |\Gamma\cap G_t|\le C\cdot m_G(G_t).
$$
\item Suppose that the family $\{G_t\}_{t>0}$ is admissible and 
$\beta_t$ satisfies 
the mean ergodic theorem in $L^p(m_{G/\Gamma})$ for some $p\ge 1$. Then
$$
\lim_{t\to \infty}\frac{|\Gamma\cap G_t|}{ m_G(G_t)}=1. 
$$
\item Suppose that the family $\{G_t\}_{t>0}$ is quasi-uniform and 
$\beta_t$ satisfies   
the pointwise ergodic theorem in $L^\infty(m_{G/\Gamma})$. Then
$$
\lim_{t\to\infty}\frac{|\Gamma\cap G_t|}{m_G(G_t)}=1. 
$$
\item Suppose that the family $\{G_t\}_{t>0}$ is admissible and $\beta_t$ satisfies 
the exponential mean ergodic theorem 
in $(L^p(m_{G/\Gamma}),L^r(m_{G/\Gamma}))$ 
for some $p\ge r\ge 1$. Then for some $\alpha>0$ (made explicit below),
$$
\frac{|\Gamma\cap G_t|}{m_G(G_t)}=1+O(e^{-\alpha t}).
$$
\end{enumerate}
\end{theorem}

\subsection{Strong maximal inequality}

We now prove some results necessary for the proof of Theorem \ref{reduct-lattice-points}. In this subsection we assume 
that the family $\{G_t\}_{t>0}$ is coarsely admissible, 
and as usual set $\Gamma_t=G_t\cap \Gamma$.

\begin{lemma}\label{l:gamma_g}
\begin{enumerate}
\item $|\Gamma_t|\le C m_G(G_t)$.
\item Assuming the strong maximal 
inequality for $\beta_t$ in 
$(L^p(m_{G/\Gamma}),L^r(m_{G/\Gamma}))$ for some $p\ge r\ge 1$,
we have $|\Gamma_t|\ge C^\prime m_G(G_t)$ for sufficiently large $t$.
\end{enumerate}
\end{lemma} 

\begin{proof}
Let $B\subset G$ be a bounded measurable subset of positive measure, 
and we assume that $B$ is small enough so 
that all of its right translates by elements of $\Gamma$ are 
pairwise disjoint. 
Then by (\ref{eq:B_max}) and (\ref{eq:B_max2}),
\begin{align*}
|\Gamma_t|&=\frac{1}{m_G(B)}
\sum_{\gamma\in\Gamma_t} 
m_G(B\gamma)=\frac{1}{m_G(B)} m_G\left(\bigcup_{\gamma\in\Gamma_t}
 B\gamma\right)\\
&\le \frac{1}{m_G(B)} m_G(G_{t+c})\le C m_G(G_t).
\end{align*}
This proves the first part of the lemma.

To prove the second part, we first show

\begin{claim}
There exists a compact set $Q\subset G/\Gamma$ 
and $x_0\in G/\Gamma$ such that
$$
\liminf_{t\to \infty} \pi_{G/\Gamma}(\beta_t)\chi_Q(x_0)=
\liminf_{t\to\infty}\frac{m_G(\{g\in G_t: gx_0\in Q\})}{m_G(G_t)}>0.
$$
\end{claim}
\begin{proof}
Suppose that the claim is false. For a compact set $Q\subset G/\Gamma$, 
denote by $\psi$
the characteristic function of the set $(G/\Gamma)\setminus Q$, the complement 
of $Q$. Then for every $x\in G/\Gamma$,
$$
\sup_{t\ge t_0}\pi_{G/\Gamma}(\beta_t)\psi(x)\ge 
\limsup_{t\to \infty} \pi_{G/\Gamma}(\beta_t)\psi(x)=1.
$$
On the other hand, 
$$
\|\psi\|_{L^p(G/\Gamma)}=m_{G/\Gamma}((G/\Gamma)\setminus Q)^{1/p},
$$
and it can be made arbitrary small by increasing $Q$. 
This contradicts the strong maximal inequality and proves the claim.
\end{proof}
\medskip

Continuing with the proof of Lemma \ref{l:gamma_g}, 
denote by $\chi_Q$ the characteristic function of the set $Q$. 
Then for some $x_0\in G/\Gamma$ 
$$
\liminf_{t\to\infty} \pi_{G/\Gamma}(\beta_t)\chi_Q(x_0)=C_0 > 0.
$$

There exists 
a non-negative measurable function $\tilde{\chi}:G\to\mathbb{R}$ with compact 
support such that  
$$
\chi_Q(g\Gamma)=\sum_{\gamma\in\Gamma} \tilde{\chi}(g\gamma)
$$
since the projection $C_c^+(G)\to C_c^+(G/\Gamma)$ by summing over 
$\Gamma$-orbits is onto. 

Letting $x_0=g_0\Gamma$, we conclude that 

$$
\int_{G_t} \sum_{\gamma\in\Gamma} \tilde{\chi}(g^{-1}g_0\gamma) dm_G(g)
\ge \frac12 C_0 m_G(G_t)
$$
for sufficiently large $t$. Now, if 
$\tilde{\chi}(g^{-1}g_0\gamma)\ne 0$ for some $g\in G_t$, then
$$
\gamma\in g_0^{-1}\cdot G_t\cdot (\supp\, \tilde{\chi})\subset G_{t+c}
$$
by (\ref{eq:B_max}). Hence,
\begin{align*}
\int_{G_t} \sum_{\gamma\in\Gamma} \tilde{\chi}(g^{-1}g_0\gamma)\, 
dm_G(g)&\le \sum_{\gamma\in \Gamma_{t+c}} \int_{G_t^{-1}g_0\gamma} \tilde{\chi}\, dm_G\\
&\le |\Gamma_{t+c}|\cdot\int_G\tilde{\chi}\,dm_G.
\end{align*}
Now Lemma \ref{l:gamma_g} follows from (\ref{eq:B_max2}).
\end{proof}
We now prove the following result, reducing the 
maximal inequality for $\lambda_t$ to the 
maximal inequality for $\beta_t$, under the assumption of coarse admissibility. 

\begin{theorem}\label{red. max}
Suppose that $\beta_t$ satisfies the strong maximal inequality in $(L^p(\nu),L^r(\nu))$,
then $\lambda_t$ satisfies the strong maximal inequality $(L^p(\mu),L^r(\mu))$.
\end{theorem}

\begin{proof}
Take $\phi\in L^p(\mu)$.

First, we observe that it suffices to prove the theorem for $\phi\ge 0$.
Write
$$
\phi=\phi^+-\phi^-
$$
where $\phi^+,\phi^-:X\to\mathbb{R}_+ $ are Borel functions such that 
$$
\max\{\phi^+,\phi^-\}\le |\phi|.
$$ 
Assuming that the strong maximal inequality holds for $\pi_X(\lambda_t)\phi^+$ and $\pi_X(\lambda_t)\phi^-$,
we get
\begin{align*}
\left\|\sup_{t\ge t_0} |\pi_X(\lambda_t)\phi|\right\|_{L^r(\mu)}&\le \left\|\sup_{t\ge t_0} |\pi_X(\lambda_t)\phi^+|\right\|_{L^r(\mu)}+\left\|\sup_{t\ge t_0} |\pi_X(\lambda_t)\phi^-|\right\|_{L^r(\mu)}\\
&\le C \|\phi^+\|_{L^p(\mu)}+C\|\phi^-\|_{L^p(\mu)}\le 2C\|\phi\|_{L^p(\mu)}.
\end{align*}
Hence, we can assume that $\phi\ge 0$.

Let $B$ be a positive-measure compact subset of $G$, 
small enough so that all of its 
right translates under $\Gamma$ are disjoint, and let  
$$
\chi=\frac{\chi_{B}}{m_G(B)},
$$
and $F:Y\to \mathbb{R}$ be defined as in (\ref{eq:F}). 

\begin{claim}
 There exists $c,d>0$ such that for all 
sufficiently large $t$, and every $h\in B$ and $x\in X$,
$$
\pi_X(\lambda_t)\phi(x)\le d\cdot\pi_Y(\beta_{t+c})F(\pi(h,x)).
$$
\end{claim}

\begin{proof}
For $(h,x)\in G\times X$, we have 
\begin{align*}
\pi_Y(\beta_t)F(\pi(h,x))&=\frac{1}{m_G(G_t)} \int_{G_t}\left(\sum_{\gamma\in\Gamma} \chi(g^{-1}h\gamma)\phi(\gamma^{-1}\cdot x)\right)\,dm_G(g)\\
&=\frac{1}{m_G(G_t)}\sum_{\gamma\in\Gamma} \left(\int_{G_t} \chi(g^{-1}h\gamma)\,dm_G(g)\right)\phi(\gamma^{-1}\cdot x).
\end{align*}
By (\ref{eq:B_max}), for $\gamma\in\Gamma_{t}$ and $h\in B$,
$$
\supp(g\mapsto \chi(g^{-1}h\gamma))=h\gamma\,\supp(\chi)^{-1}\subset G_{t+c}.
$$
Hence,
$$
\int_{G_{t+c}} \chi(g^{-1}h\gamma)\,dm_G(g)=1.
$$
Also, by Lemma \ref{l:gamma_g} and (\ref{eq:B_max2}),
$$
|\Gamma_t|\ge C^\prime m_G(G_{t+c}).
$$

Applying the previous arguments to $\pi(\beta_{t+c})$, and summing 
only on $\gamma\in \Gamma_t$,  
we conclude that for $(h,x)\in B\times X$, 
\begin{align*}
\pi_Y(\beta_{t+c})F(\pi(h,x))&\ge\frac{1}{m_G(G_{t+c})}\sum_{\gamma\in\Gamma_{t}} \left(\int_{G_{t+c}} \chi(g^{-1}h\gamma)\,dm_G(g)\right)\phi(\gamma^{-1}\cdot x)\\
&=\frac{1}{m_G(G_{t+c})}\sum_{\gamma\in\Gamma_{t}} \phi(\gamma^{-1}\cdot x)\ge
C'' \pi_X(\lambda_t)\phi(x).
\end{align*}
This proves the claim.
\end{proof}
\medskip

Continuing with the proof of Theorem \ref {red. max}, 
we now take the supremum over $t$ on both sides. 
Let us lift $\pi_X(\lambda_t)\phi$ to be defined on $B\times X$ (depending only 
the second coordinate). By the claim, 
 for sufficiently large $t_0'>0$, integrating over
$h\in B$ we obtain 
\begin{align*}
\left\|\sup_{t\ge t_0'} |\pi_X(\lambda_t)\phi)|\right\|_{L^r(\mu)}&=m_G(B)^{-1/r}\left\|\sup_{t\ge t_0'} |\pi_X(\lambda_t)\phi|\right\|_{L^r(m_G\otimes\mu|_{B\times X})}\\
&\le C^\prime
 \left\|\sup_{t\ge t_0'} |\pi_Y(\beta_{t+c})(F\circ\pi)|\right\|_{L^r(m_G\otimes\mu|_{B\times X})}.
\end{align*}
Now $\pi_Y(\beta_t)(F\circ \pi)=\left(\pi_Y(\beta_t)F\right)\circ \pi$, 
since the left $G$-action on $G\times X$ 
commutes with the right $\Gamma$-action. Hence, by Lemma 
\ref{l:lp}(b) and the strong maximal 
inequality for $\beta_t$ in $(L^p(\nu), L^r(\nu))$,
\begin{align*}
\left\|\sup_{t\ge t_0'} |\pi_X(\lambda_t)\phi)|\right\|_{L^r(\mu)}
&\le C^\prime b_{r,B}
\left\|\sup_{t\ge t_0'} |\pi_Y(\beta_{t+c})F|\right\|_{L^r(\nu)}\le 
C'' \|F\|_{L^p(\nu)}\\
&= C'' \|\chi\|_{L^p(m_G)}\cdot \|\phi\|_{L^p(\mu)}\le C
 \|\phi\|_{L^p(\mu)}.
\end{align*}
where the equality uses
 the fact that $B$ has disjoint right translates under $\Gamma$ 
and Lemma \ref{l:lp}(a). 

This concludes the proof of Theorem \ref{red. max}.
\end{proof}

\subsection{Mean ergodic theorem}

We now turn from maximal inequalities to establishing convergence 
 results for averages on $\Gamma$, using smoothing to    
approximate discrete averages by absolutely 
continuous ones, and thus utilizing the finiteness of 
the upper local dimension of $G$. 
Generalizing the definition of upper local dimension 
 somewhat, consider a base of neighborhoods 
$\{\mathcal{O}_\vre\}_{0<\vre<1}$ of $e$ in $G$ such that 
$\mathcal{O}_\vre$'s are symmetric, 
bounded, and increasing with $\vre$.
We assume that the family $\{G_t\}_{t>0}$ satisfy the following conditions:
\begin{itemize}
\item There exists $c>0$ such that for every small $\vre>0$ and $t\ge t(\vre)$,
\begin{align}\label{eq:point}
\mathcal{O}_\vre\cdot G_t\cdot \mathcal{O}_\vre &\subset G_{t+c\vre}.
\end{align}
\item For 
$$
\delta_\vre=\limsup_{t\to\infty}\frac{m_G(G_{t+\vre}-G_t)}{m_G(G_t)}\,,
$$
and for some $p\ge 1$, we have : 
\begin{equation}\label{eq:delta_e}
\delta_\vre^p \cdot m_G(\mathcal{O}_\vre)^{-1}\to 0\quad\hbox{as}\quad \vre\to 0^+.
\end{equation}
\end{itemize}
Note that if the family $\{G_t\}_{t>0}$ is admissible and
$$
\varrho_0=\limsup_{\vre\to 0^+} \frac{\log m_G(\mathcal{O}_\vre)}{\log\vre}<\infty,
$$
then (\ref{eq:delta_e}) holds for $p>\varrho_0$.

Note that (\ref{eq:delta_e}) implies that $\delta_\vre\to 0$ as $\vre\to 0^+$. 
For every $\delta>\delta_\vre$ and for sufficiently large $t$,
\begin{equation}\label{eq:delta_e2}
m_G(G_{t+\vre})\le (1+\delta)m_G(G_t).
\end{equation}

\begin{lemma}\label{l:gG}
Under condition (\ref{eq:delta_e}), 
if the mean ergodic 
theorem holds for $\beta_t$ in 
$L^q(m_{G/\Gamma})$ for some $q\ge 1$, then
$$
|\Gamma_t|\sim m_G(G_t)\quad\hbox{as}\quad t\to\infty.
$$
\end{lemma}

\begin{proof}
Let
$$
\chi_\vre=\frac{\chi_{\mathcal{O}_\vre}}{m_G(\mathcal{O}_\vre)}
$$
and
$$
\phi_\vre(g\Gamma)=\sum_{\gamma\in \Gamma} \chi_{\vre}(g\gamma).
$$
Note that $\phi$ is a measurable bounded function on $G/\Gamma$ with compact support,
$$
\int_G\chi_\vre\,dm_G=1,\quad\hbox{and}\quad\int_{G/\Gamma}\phi_\vre\,dm_{G/\Gamma}=1.
$$
It follows from the mean ergodic theorem that for every $\delta>0$,
$$
m_{G/\Gamma}(\{g\Gamma\in G/\Gamma:\, |\pi_{G/\Gamma}(\beta_t)\phi_\vre(g\Gamma) -1|>\delta\})\to 0\quad\hbox{as}\quad t\to\infty.
$$
In particular, for sufficiently large $t$, there exists $g_t\in\mathcal{O}_\vre$ such that 
$|\pi_{G/\Gamma}(\beta_t)\phi_\vre(g_t\Gamma) -1|\le\delta$, or equivalently 
\begin{equation}\label{eq:G_Gamma}
1-\delta \le \frac{1}{m_G(G_t)}
\int_{G_t}\phi_\vre(g^{-1}g_t\Gamma)dm_G(g)\le 1+\delta
\end{equation}
Thus let us now prove the following 

\begin{claim}
Given $0 < \vre \le \vre_0$, for every $t\ge t_0+\vre_0$ 
and for every $h\in \mathcal{O}_\vre$,
$$
\int_{G_{t-c\vre}}\phi_\vre(g^{-1}h\Gamma)\,dm_G(g)\le |\Gamma_t|\le \int_{G_{t+c\vre}}\phi_\vre(g^{-1}h\Gamma)\,dm_G(g).
$$
\end{claim}

Indeed, if $\chi_\vre(g^{-1}h\gamma)\ne 0$ for some $g\in G_{t-c\vre}$ and $h\in\mathcal{O}_\vre$, then
$$
\gamma\in h^{-1}\cdot G_{t-c\vre}\cdot (\supp\, \chi_\vre)\subset G_t.
$$
Hence, 
$$
\int_{G_{t-c\vre}}\phi_\vre(g^{-1}h\Gamma)\,dm_G(g)\le \sum_{\gamma\in \Gamma_t} \int_{G_t} \chi_\vre(g^{-1}h\gamma)\,dm_G(g)\le |\Gamma_t|.
$$
In the other direction, for 
$\gamma\in\Gamma_t$ and $h\in\mathcal{O}_\vre$,
$$
\supp(g\mapsto \chi_\vre(g^{-1}h\gamma))=h\gamma(\supp\,\chi_\vre)^{-1}\subset G_{t+c\vre}.
$$
Since $\chi_\vre\ge 0$,
$$
\int_{G_{t+c\vre}}\phi_\vre(g^{-1}h\Gamma)\,dm_G(g)\ge \sum_{\gamma\in\Gamma_t} \int_{G_{t+c\vre}}\chi_\vre(g^{-1}h\gamma)\,dm_G(g)\ge |\Gamma_t|\,.
$$
and this establishes the claim.\qed

\medskip
Continuing with the proof of Lemma \ref{l:gG}, 
let us take $h=g_t$ defined above. By the claim and
 (\ref{eq:G_Gamma}),
$$
|\Gamma_t|\le (1+\delta)m_G(G_{t+\vre}),
$$
and the upper estimate on $|\Gamma_t|$ follows from (\ref{eq:delta_e2}). The lower estimate is proved similarly.
\end{proof}

We now generalize Lemma \ref{l:gG}, 
and prove the following result,
 reducing the mean ergodic theorem for $\lambda_t$ to the mean 
ergodic theorem for $\beta_t$.

\begin{theorem}\label{th:mean}
Under condition (\ref{eq:delta_e}), 
if the mean ergodic theorem holds for 
$\beta_t$ in $L^p(\nu)$,
then the mean ergodic 
theorem holds for $\lambda_t$ in $L^p(\mu)$.

\end{theorem}

\begin{proof}
Take small $\vre>0$ and $\delta\in (\delta_\vre,1)$, where $\delta_\vre$ 
(as well as $p$) are defined by (\ref{eq:delta_e}).

We need to show that for every $\phi\in L^p(\mu)$
$$
\left\|\pi_X(\lambda_t)\phi-\int_X\phi\, d\mu\right\|_{L^p(\mu)}\to 0\quad\hbox{as}\quad t\to\infty,
$$
and without loss of generality, we may assume that $\phi \geq 0$.
Let
$$
\chi_\vre=\frac{\chi_{\mathcal{O}_\vre}}{m_G(\mathcal{O}_\vre)}
$$
and $F_\vre:Y\to \mathbb{R}$ be defined as in (\ref{eq:F}). Then $F_\vre\in L^p(\nu)$, and 
$$
\int_Y F_\vre\, d\nu=\int_X\phi\,d\mu.
$$

\begin{step1}
For every $(g,x)\in \mathcal{O}_\vre\times X$ and sufficiently large $t$, 
$$
(1+\delta)^{-1}\pi_Y(\beta_{t-c\vre})F_\vre(\pi(g,x))\le \pi_X(\lambda_t)\phi(x)\le (1+\delta)\pi_Y(\beta_{t+c\vre})F_\vre(\pi(g,x))
$$
\end{step1}

To prove the first inequality, 
note that by Lemma \ref{l:gG} and (\ref{eq:delta_e2}),
$$
(1+\delta)^{-1}m_G(G_{t+c\vre})<|\Gamma_t|<(1+\delta)m_G(G_{t-c\vre})
$$
for suffciently large $t$.

For $(h,x)\in G\times X$,
\begin{align*}
\pi_Y(\beta_t)F_\vre(\pi(h,x))&=\frac{1}{m_G(G_t)} \int_{G_t}\left(\sum_{\gamma\in\Gamma} \chi_\vre(g^{-1}h\gamma)\phi(\gamma^{-1}\cdot x)\right)\,dm_G(g)\\
&=\frac{1}{m_G(G_t)}\sum_{\gamma\in\Gamma} \left(\int_{G_t} \chi_\vre(g^{-1}h\gamma)\,dm_G(g)\right)\phi(\gamma^{-1}\cdot x).
\end{align*}
If $\chi_\vre(g^{-1}h\gamma)\ne 0$ for some $g\in G_t$ and $h\in \mathcal{O}_\vre$, then by (\ref{eq:point}),
$$
\gamma\in h^{-1}g\,\supp(\chi_\vre)\subset G_{t+c\vre}.
$$
Using that
\begin{equation}\label{eq_chi_e0}
\int_G \chi_\vre\, dm_G=1\quad\hbox{and}\quad \chi_\vre\ge 0,
\end{equation}
we deduce that for $(h,x)\in \mathcal{O}_\vre\times X$,
\begin{align*}
\pi_Y(\beta_t)F_\vre(\pi(h,x))&= 
\frac{1}{m_G(G_t)}\sum_{\gamma\in\Gamma_{t+c\vre}} \left(\int_{G_t} \chi_\vre(g^{-1}h\gamma)\,dm_G(g)\right)\phi(\gamma^{-1}\cdot x)\\
&\le\frac{1}{m_G(G_t)}\sum_{\gamma\in\Gamma_{t+c\vre}} \phi(\gamma^{-1}\cdot x)\le (1+\delta)\pi_X(\lambda_{t+c\vre})\phi(x).
\end{align*}

To prove the second inequality, note that by
 (\ref{eq:1}), for $\gamma\in\Gamma_{t-c\vre}$ 
and $h\in \mathcal{O}_\vre$,
$$
\supp(g\mapsto \chi_\vre(g^{-1}h\gamma))=h\gamma\,\supp(\chi_\vre)^{-1}\subset G_t.
$$
By (\ref{eq_chi_e0}), this implies that for $(h,x)\in \mathcal{O}_\vre\times X$,
\begin{align*}
\pi_Y(\beta_{t})F_\vre(\pi(h,x))&\ge\frac{1}{m_G(G_{t})}\sum_{\gamma\in\Gamma_{t-c\vre}} \left(\int_{G_t} \chi_\vre(g^{-1}h\gamma)\,dm_G(g)\right)\phi(\gamma^{-1}\cdot x)\\
&=\frac{1}{m_G(G_t)}\sum_{\gamma\in\Gamma_{t-c\vre}} \phi(\gamma^{-1}\cdot x)\ge (1+\delta)^{-1}\pi_X(\lambda_{t-c\vre})\phi(x).
\end{align*}
Using Lemma \ref{l:gG} and (\ref{eq:delta_e2}) again, and 
then shifting indices completes the proof of Step 1.

\medskip

We now continue with the proof of Theorem \ref{th:mean}. 
To simplify notations, 
we write for a measurable function $\Psi:Y\to\mathbb{R}$ 
$$
\|\Psi\|_{p,\vre}\stackrel{def}{=}
\|\Psi\circ\pi\|_{L^p(m_G\otimes\mu|_{\mathcal{O}_\vre\times X})}
$$

Now by Lemma \ref{l:lp}(b),
for each fixed $\vre > 0$
\begin{equation}\label{eq:psi_vre}
\|\Psi\|_{p,\vre}\le  b_{p,\vre}   
 \|\Psi\|_{L^p(\nu)}
\end{equation}
and clearly, if $\vre^\prime < \vre$ we may take 
$b_{p,\vre^\prime}\le b_{p,\vre}$.

\begin{step1}
For every sufficiently small fixed $\vre>0$,
$$
\|\pi_Y(\beta_{t+c\vre})F_\vre-\pi_Y(\beta_{t})F_\vre\|_{p,\vre}
\to 0\quad\hbox{as}\quad t\to\infty
$$
and
$$
\limsup_{t\to\infty}\|\pi_Y(\beta_{t})F_\vre\|_\vre\le b_{p,\vre}
 \|\phi\|_{L^1(\mu)}.
$$
\end{step1}

For the proof, let us 
note that 
by the triangle inequality and (\ref{eq:psi_vre}),
\begin{align*}
\|\pi_Y(\beta_{t+c\vre})F_\vre-\pi_Y(\beta_{t})F_\vre\|_{p,\vre}&\le b_{p,\vre} 
\left\|\pi_Y(\beta_{t+c\vre})F_\vre-\int_Y F_\vre\, d\nu\right\|_{L^p(\nu)}\\
&+b_{p,\vre} \left\|\pi_Y(\beta_{t})F_\vre-\int_Y F_\vre\, d\nu\right\|_{L^p(\nu)}.
\end{align*}
Similarly,
\begin{align*}
\|\pi_Y(\beta_{t})F_\vre\|_{p,\vre} &\le b_{p,\vre}
\left\|\pi_Y(\beta_{t})F_\vre-\int_Y F_\vre\, d\nu\right\|_{L^p(\nu)}+
b_{p,\vre} \int_Y F_\vre\, d\nu\\
&= b_{p,\vre}\left\|\pi_Y(\beta_{t})F_\vre-\int_Y F_\vre\, 
d\nu\right\|_{L^p(\nu)}+ b_{p,\vre}\int_X \phi\, d\mu.
\end{align*}
Hence, Step 2 follows from the mean ergodic theorem for $\beta_t$ in $L^p(\nu)$.

\medskip

To complete the proof of Theorem \ref{th:mean}, we need to estimate
\begin{align*}
\left\|\pi_X(\lambda_t)\phi -\int_X \phi\, d\mu\right\|_{L^p(\mu)}&=
m_G(\mathcal{O}_\vre)^{-1/p}\left\|\pi_X(\lambda_t)\phi -\int_X \phi\, d\mu\right\|_{L^p(m_G\otimes\mu|_{\mathcal{O}_\vre\times X})}
\end{align*}
where we have extended $\pi_X(\lambda_t)\phi$ to a function on $\mathcal{O}_\vre\times X$ in the obvious manner. By the triangle inequality,
\begin{align*}
&\quad\left\|\pi_X(\lambda_t)\phi -\int_X \phi\, d\mu\right\|_{L^p(m_G\otimes\mu|_{\mathcal{O}_\vre\times X})}\\
&\le \left\|\pi_X(\lambda_t)\phi -(1+\delta)^{-1}\pi_Y(\beta_{t-c\vre})(F_\vre\circ\pi)\right\|_{L^p(m_G\otimes\mu|_{\mathcal{O}_\vre\times X})}\\
&+\left\|(1+\delta)^{-1}\pi_Y(\beta_{t-c\vre})(F_\vre\circ\pi)-\int_X \phi\, d\mu\right\|_{L^p(m_G\otimes\mu|_{\mathcal{O}_\vre\times X})}.
\end{align*}
We estimate the last two summands as follows. 
First, using Step 1, we estimate the first summand by 
\begin{align*}
&\quad\left\|\pi_X(\lambda_t)\phi -(1+\delta)^{-1}\pi_Y(\beta_{t-c\vre})
(F_\vre\circ\pi)\right\|_{L^p(m_G\otimes\mu|_{\mathcal{O}_\vre\times X})}\\
&\le\left\|(1+\delta)\pi_Y(\beta_{t+c\vre})F_\vre 
-(1+\delta)^{-1}\pi_Y(\beta_{t-c\vre})F_\vre\right\|_{p,\vre}\\
&\le \left\|\pi_Y(\beta_{t+c\vre})F_\vre -\pi_Y(\beta_{t-c\vre})
F_\vre\right\|_{p,\vre}+\delta\left(\left\|\pi_Y(\beta_{t+c\vre})F_\vre
\right\|_{p,\vre}+\left\|\pi_Y(\beta_{t-c\vre})F_\vre\right\|_{p,\vre}\right).
\end{align*}
Hence, it follows from Step 2 that
\begin{align*}
\limsup_{t\to\infty}\left\|\pi_X(\lambda_t)\phi -(1+\delta)^{-1}\pi_Y(\beta_{t-c\vre})(F_\vre\circ\pi)\right\|_{L^p(m_G\otimes\mu|_{\mathcal{O}_\vre\times X})}
\le 2 b_{p,\vre} \delta\|\phi\|_{L^1(\mu)}.
\end{align*}
Second, observing that for  $\delta < 1$ 
\begin{align*}
\left\|(1+\delta)^{-1}\pi_Y(\beta_{t-c\vre})F_\vre-\int_X \phi\, 
d\mu\right\|_{p,\vre}
&\le \left\|\pi_Y(\beta_{t-c\vre})F_\vre-\int_X \phi\, d\mu\right\|_{p,\vre}\\
&+2\delta\left\|\pi_Y(\beta_{t-c\vre})F_\vre\right\|_{p,\vre},
\end{align*}
we deduce from step 2 that the second summand is estimated by 
$$
\limsup_{t\to\infty}\left\|(1+\delta)^{-1}\pi_Y(\beta_{t-c\vre})
F_\vre-\int_X \phi\, d\mu\right\|_{p,\vre}\le 4b_{p,\vre} \delta\|\phi\|_{L^1(\mu)}\,,
$$
where $b_{p,\vre}$ are uniformly bounded. 

We have thus shown that for every $\delta\in (\delta_\vre, 1)$,
and a constant $B$ independet of $\delta$ and $\vre$, 
$$
\limsup_{t\to\infty}\left\|\pi_X(\lambda_t)
\phi -\int_X \phi\, d\mu\right\|_{L^p(\mu)}\le 
B \delta m_G(\mathcal{O}_\vre)^{-1/p}\|\phi\|_{L^1(\mu)}.
$$ 
By our choice of $\delta_\vre$ in 
(\ref{eq:delta_e}), this concludes the proof  of Theorem 
\ref{th:mean}.
\end{proof}

\subsection{Pointwise ergodic theorem}
In this section, we assume only 
that the family $\{G_t\}_{t>0}$ is 
quasi-uniform. Recall that we showed in 
Corollary \ref{gGth} that then $gG_th$ satisfy the 
pointwise ergodic theorem if $G_t$ does.  
\begin{lemma}
Suppose that the pointwise ergodic 
theorem holds in $L^\infty(G/\Gamma)$ 
for the quasi-uniform family $\{gG_t\}$, 
for every  $g\in G$.
Then
$$
|\Gamma_t|\sim m_G(G_t)\quad\hbox{as}\quad t\to\infty.
$$
\end{lemma}

\begin{proof}
Let $f$ be any measurable bounded function on $G/\Gamma$.  
$G/\Gamma$ being a homogeneous $G$-space, it follows 
from Theorem \ref{th:g_inv} that the 
pointwise ergodic theorem holds for {\it every point}
 in $G/\Gamma$. In particular, this holds for the point $e\Gamma$, and so 
$$
\frac{1}{m_G(G_t)}\int_{G_t} f(g^{-1}
\Gamma)dm_G(g)\to\int_{G/\Gamma} f\, dm_{G/\Gamma}
$$
for every measurable bounded $f$. The lemma is then 
proved as Proposition 6.1 in [GW].
\end{proof}

\begin{theorem}\label{th:pointwise}
If the pointwise ergodic theorem holds for the 
quasi-uniform family 
$\beta_t$ in $L^p(\nu)$, then the 
pointwise ergodic theorem holds for $\lambda_t$ in $L^p(\mu)$.
\end{theorem}

\begin{proof}

We need to show that for every $\phi\in L^p(\mu)$,
$$
\pi_X(\lambda_t)\phi(x)\to \int_X\phi\, d\mu\quad\hbox{as}\quad t\to\infty
$$
for $\mu$-a.e. $x\in X$  and without loss of generality, we may assume that $\phi \geq 0$.

Take $\delta>0$ and let $\vre>0$ and $\mathcal{O}$ be as in (\ref{eq:01}) and (\ref{eq:02}).
Let
$$
\chi=\frac{\chi_{\mathcal{O}}}{m_G(\mathcal{O})}
$$
and $F:Y\to \mathbb{R}$ be defined as in (\ref{eq:F}). Then $F\in L^p(\nu)$ and
$$
\int_Y F\, d\nu=\int_X\phi\,d\mu.
$$
Recalll that it follows from Corollary \ref{gGth} that 
the pointwise ergdoic theorem holds for the family $gG_t$ in 
$L^\infty\subset L^p$, for every $g\in G$.  
Using also the assumption of Theorem \ref{th:pointwise}, 
and Theorem \ref{th:g_inv},
\begin{equation}\label{eq:beta_F}
\pi_Y(\beta_t)F(y)\to \int_Y\phi\, d\nu\quad\hbox{as}\quad t\to\infty
\end{equation}
for $y$ in a $G$-invariant subset of $Y$ of full measure. 
Then it follows from Lemma \ref{strict} that (\ref{eq:beta_F})
holds for $y=\pi(e,x)$ for $x$ in a $\Gamma$-invariant subset of $X$ of full measure. Arguing as in the proof of Proposition 2.1
in [GW], one shows that for every such $x$,
\begin{align*}
\limsup_{t\to\infty} \pi_X(\lambda_t)\phi(x)&\le (1+\delta)\int_X\phi\, d\mu,\\
\liminf_{t\to\infty} \pi_X(\lambda_t)\phi(x)&\ge (1+\delta)^{-1}\int_X\phi\, d\mu,
\end{align*}
for every $\delta>0$. This completes the proof of Theorem 
\ref{th:pointwise}.
\end{proof}

\subsection{Exponential mean ergodic theorem}

In this section we assume that  the family $\{G_t\}_{t>0}$ is admissible, and as usual 
w.r.t. a family $\mathcal{O}_\vre$ of finite upper local dimension $\varrho_0$.
By definition for $\varrho>\varrho_0$ 
and small $\vre>0$,
\begin{equation}\label{eq:epsilon_rho}
m_G(\mathcal{O}_\vre)\ge C_\rho \vre^{\varrho}.
\end{equation}

\begin{theorem}\label{l:gG2}
If the exponential mean 
ergodic theorem holds for the admissible family 
$\beta_t$ in $(L^p(m_{G/\Gamma}),L^r(m_{G/\Gamma}))$ for some $p\ge r\ge 1$, then
$$
\frac{|\Gamma_t|}{m_G(G_t)}=1+O(e^{-\alpha t})\text{    where    }
\alpha=\frac{\theta_{p,r}}{\varrho(1+r-\frac{r}{p})+r}\,\,.$$ 
When the estimate $\norm{\pi_{G/\Gamma}^0(\beta_t)}_{L^2_0\to L^2_0}\le Ce^{-\theta t}$ holds, 
we can take $\theta_{p,r}=r\theta$. 
\end{theorem}

\begin{proof}
As usual let 
$$
\chi_\vre=\frac{\chi_{\mathcal{O}_\vre}}{m_G(\mathcal{O}_\vre)}
$$
and
$$
\phi_\vre(g\Gamma)=\sum_{\gamma\in \Gamma} \chi_{\vre}(g\gamma), 
$$
so that that $\phi_\vre$ is a bounded function on $G/\Gamma$ with compact support,
$$
\int_G\chi_\vre\,dm_G=1,
\quad\hbox{and}\quad\int_{G/\Gamma}\phi_\vre\,dm_{G/\Gamma}=1.
$$
It follows from the norm estimate given by the 
$(L^p,L^r)$-exponential mean ergodic 
theorem for $\beta_t$ acting on $G/\Gamma$ 
that for some fixed $\theta_{p,r}>0$ and $C > 0$, and for 
every $\delta>0$, $t > 0$ and 
$\vre> 0$
$$
m_{G/\Gamma}(\{x\in G/\Gamma:\, |\pi_{G/\Gamma}(\beta_t)
\phi_\vre(x) -1|>\delta\})\le C \delta^{-r}\|\phi_\vre\|_{L^p(G/\Gamma)}^r 
e^{-\theta_{p,r}t}\,.
$$
 By Lemma \ref{l:lp}(a), for sufficiently small $\vre > 0$ 
$$
\|\phi_\vre\|_{L^p(G/\Gamma)}= m_G(\mathcal{O}_\vre)^{1/p-1}.
$$
We will choose both of 
 the parameters $\vre$ and $\delta$ as a function of $t$, and 
begin by requiring that the following condition holds:
\begin{equation}\label{eq:cond_exp}
C \delta^{-r}m_G(\mathcal{O}_\vre)^{r/p-r}e^{-\theta_{p,r} t}
 = \frac12 
m_G(\mathcal{O}_\vre).
\end{equation}
Then for sufficiently large $t$,
\begin{align*}
m_{G/\Gamma}(\{x\in G/\Gamma:\, |\pi_{G/\Gamma}(\beta_t)\phi_\vre(x) -1|>\delta\})
& \le \frac12 m_G(\mathcal{O}_\vre).
\end{align*}
Now as soon as $\mathcal{O}_\vre$ maps injectively into $G/\Gamma$, 
we have 
$$
m_{G/\Gamma}(\mathcal{O}_\vre\Gamma)=m_G(\mathcal{O}_\vre)
$$
and so we deduce that for every sufficiently large $t$,
there exists $g_t\in\mathcal{O}_\vre$ such that
$$
|\pi_{G/\Gamma}(\beta_t)\phi_\vre(g_t\Gamma) -1|\le\delta\,.
$$
Then using the claim from Lemma \ref{l:gG} and (\ref{eq:2}), 
we have 
for sufficiently large $t$ 
$$
|\Gamma_t|\le (1+\delta)m_G(G_{t+c\vre})\le (1+\delta)(1+c\vre)m_G(G_t)\,\,,
$$
provided only that $\vre$, $\delta$ and $t$ satisfy condition 
(\ref{eq:cond_exp}). 

In order to balance the two significant 
parts of the error estimate 
let us take $c\vre=\delta$. Then (\ref{eq:cond_exp}) together 
with (\ref{eq:epsilon_rho}) yield 
$$C^\prime e^{-\theta_{p,r} t}=\vre^r 
m_G(\mathcal{O}_\vre)^{1+r-\frac{r}{p}}\ge C'' 
\vre^{\varrho(1+r-\frac{r}{p})+r}\,\,.$$
Thus we take both $\vre$ and $\delta$ to be constant multiples of 
$$\exp\left(-\frac{\theta_{p,r} t}{\varrho(1+r-\frac{r}{p})+r}
\right)$$ and conclude that 
$$
\abs{\frac{|\Gamma_t|}{m_G(G_t)}-1}\le B
\exp\left(-\frac{\theta_{p,r}t}{\varrho(1+r-\frac{r}{p})+r}\right)$$ 
where $B$ is independent of $t$. 

The lower estimate is proved similarly.

The last statement of the theorem follows immediately 
from Riesz-Thorin interpolation. 
\end{proof}


We now turn to the mean ergodic theorem with exponential 
rate of convergence.

\begin{theorem}\label{th:mean_exp}
Let $\beta_t$ be an admissible family, $\mathcal{O}_\vre$ of finite upper local dimension, and let $p\ge r>\varrho> \varrho_0$.
If the exponential mean ergodic theorem holds 
for $\beta_t$ in $(L^p(\nu),L^r(\nu))$, then
the exponential mean ergodic theorem holds for $\lambda_t$ in 
$(L^p(\mu),L^r(\mu))$.
\end{theorem}

\begin{proof}
We need to show that for some $\zeta_p>0$, $C_p > 0$ and every 
$\phi\in L^p(\mu)$
$$
\left\|\pi_X(\lambda_t)\phi-\int_X\phi\, d\mu\right\|_{L^r(\mu)}
\le C e^{-\zeta t}\|\phi\|_{L^p(\mu)}.
$$
Without loss of generality, we may assume that $\phi \geq 0$.
Let
$$
\chi_\vre=\frac{\chi_{\mathcal{O}_\vre}}{m_G(\mathcal{O}_\vre)}
$$
and $F_\vre:Y\to \mathbb{R}$ be defined as in (\ref{eq:F}). 
Then $F_\vre\in L^p(\nu)$, and
$$
\int_Y F_\vre\, d\nu=\int_X\phi\,d\mu.
$$
In particular,
\begin{equation}\label{eq:mean_exp}
\left\|\pi_Y(\beta_t)F_\vre-\int_Y F_\vre\, 
d\nu\right\|_{L^r(\nu)}\le C e^{-\theta t}\|F_\vre\|_{L^p(\nu)}
\end{equation}
for some $\theta=\theta_{p,r}>0$.


Repeating the arguments in Steps 1-2 of the proof of 
Theorem \ref{th:mean}, we derive that
\begin{enumerate} 
\item[(a)] For sufficiently large $t$ and every 
$(g,x)\in \mathcal{O}_\vre\times X$, 
$$
\eta^{-1}\pi_Y(\beta_{t-c\vre})F_\vre(\pi(g,x))\le \pi_X(\lambda_t)\phi(x)\le \eta\pi_Y(\beta_{t+c\vre})F_\vre(\pi(g,x))
$$
where
$$
\eta=(1+c\vre)(1+O(e^{-\alpha t}))\,,
$$
$\alpha$ the error estimate in the lattice point count from Theorem 
\ref{l:gG2}. 

\item[(b)] For sufficiently large $t$ and small $\vre>0$,
$$
\|\pi_Y(\beta_{t+c\vre})F_\vre-\pi_Y(\beta_{t})F_\vre\|_{L^r(m_G\otimes\mu|_{\mathcal{O}_\vre\times X})}\ll e^{-\theta t}\|F_\vre\|_{L^p(\nu)}
$$
and
$$
\|\pi_Y(\beta_{t})F_\vre\|_{L^r(m_G\otimes\mu|_{\mathcal{O}_\vre\times X})} \ll e^{-\theta t}\|F_\vre\|_{L^p(\nu)}+\|\phi\|_{L^1(\mu)}.
$$
\end{enumerate}

Using (a),(b) and (\ref{eq:mean_exp}),  we deduce as in the proof of Theorem \ref{th:mean} that
\begin{align*}
&\quad\left\|\pi_X(\lambda_t)\phi -\int_X \phi\, d\mu\right\|_{L^r(\mu)}\\
&\le C^\prime m_G(\mathcal{O}_\vre)^{-1/r}\left( e^{-\theta t}\|F_\vre\|_{L^p(\nu)}+(\vre+e^{-\alpha t})\|\phi\|_{L^1(\mu)}\right)
\end{align*}
for every small $\vre>0$.

Using Lemma \ref{l:lp}(a) and (b), we can estimate 
$$\norm{F_\vre}_{L^p(\nu)}\le c_{p,\vre}
m_G(\mathcal{O}_\vre)^{\frac1p -1} \norm{\phi}_{L^p(\mu)}\,.
$$

Let now $\varrho$ be such that $\varrho_0<\varrho<r.$
 By H\"older's inequality, $\norm{\phi}_{L^1(\mu)}\le 
\norm{\phi}_{L^p(\mu)}$, and collecting terms, we obtain 
\begin{align*}
&\quad\left\|\pi_X(\lambda_t)\phi -\int_X
 \phi\, d\mu\right\|_{L^r(\mu)}\\
&\le C'' m_G(\mathcal{O}_\vre)^{-1/r}\left( e^{-\theta t}m_G(\mathcal{O}_\vre)^{1/p-1}+\vre+e^{-\alpha t}\right)\|\phi\|_{L^p(\mu)}\\
& \le (e^{-\theta t}\vre^{-\varrho(1+1/r-1/p)} +\vre^{1-\varrho/r}+e^{-\alpha t})\|\phi\|_{L^p(\mu)}.
\end{align*}
Setting
$$
\varepsilon=\exp\left(-\frac{\theta_{p,r} t}{1+\varrho-\varrho/p}\right),
$$
we have
$$
\left\|\pi_X(\lambda_t)\phi -\int_X \phi\, d\mu\right\|_{L^p(\mu)}\ll e^{-\zeta t}\|\phi\|_{L^p(\mu)}
$$
with
$$
\zeta=\min\left\{\alpha,\frac{(1-\varrho/r)
\theta_{p,r}}{1+\varrho-\varrho/p}\right\}>0.
$$
This concludes the proof of Theorems \ref{th:mean_exp}. 
\end{proof}

\subsection{Exponential strong  maximal inequality}

We now prove that the 
exponential decay of the norms $\pi_Y^0(\beta_t)$ in 
$L_0^2(Y)$, together with the ordinary strong maximal 
inequality for $\pi_Y(\beta_t)$ 
in some $L^q(Y)$  imply  
the exponential strong maximal inequality
for $\lambda_t$ following the method developed in \cite{MNS}\cite{N3}.

\begin{theorem}\label{th:exp.max}
Let $\beta_t$ be admissible averages w.r.t. $\mathcal{O}_\vre$ of finite upper local dimension $\varrho_0$, satisfying 
\begin{itemize}
\item  the exponential mean ergodic theorem in 
$(L^p(\nu),L^r(\nu))$ with $\varrho_0<r\le p$,
\item the strong maximal inequality in $L^q(\nu)$ 
for some $q\ge 1$.
\end{itemize}
Then $\lambda_t$ satisfies the exponential strong maximal inequality in $(L^v(\mu), L^w(\mu))$
for $v$, $w$ such that $1/v=(1-u)/q$ and $1/w=(1-u)/q+u/r$ for some $u\in (0,1)$.
\end{theorem}

\begin{proof}
By Theorem \ref{th:mean_exp}, $\pi_X(\lambda_t)$ satisfies the exponential mean ergodic theorem:
for every $f\in L^p_0(\mu)$, and $\theta=\theta_{p,r}> 0$
\begin{equation}\label{eq:exp_proof1}
\|\pi_X(\lambda_t)f\|_{L^r(\mu)}\le C e^{-\theta t}\|f\|_{L^p(\mu)}.
\end{equation}
Consider an increasing sequence $\{t_i\}$ that 
contains all positive integers and  divides each interval of the form
$[n,n+1]$, $n\in \mathbb{N}$, into $\lfloor e^{r\theta n/4}\rfloor$ subintervals of equal length. Then
\begin{equation}\label{eq:t_i}
t_{i+1}-t_i\le e^{-r\theta \lfloor t_i\rfloor/4},
\end{equation}
and
$$
\sum_{i\ge 0} e^{-r\theta t_i/2}\le \sum_{n\ge 0} \lfloor e^{r\theta n/4}\rfloor e^{-r\theta n/2}<\infty.
$$
Hence, it follows from (\ref{eq:exp_proof1}) that
$$
\int_X\left(\sum_{i\ge 0} e^{r\theta t_i/2} |\pi_X(\lambda_{t_i})f(x)|^r\right)\, d\mu(x)\le  C^\prime\|f\|_{L^p(\mu)}^r.
$$
Setting
$$
B_r(x,f)\stackrel{def}{=}
 \left(\sum_{i\ge 0} e^{r\theta t_i/2} |\pi_X(\lambda_{t_i})f(x)|^r\right)^{1/r},
$$
we have, for some $C''$ independent of $f$ :
\begin{align}
|\pi_X(\lambda_{t_i})f(x)|&\le B_r(x,f)
e^{-\theta t_i/2},\label{eq:exp_proof2}\\
\|B_r(\cdot,f)\|_{L^r(\mu)}&\le C''\|f\|_{L^p(\nu)}.\label{eq:exp_proof3}
\end{align}

\begin{claim}
For any sufficiently large $t$, there exists $t_i$ such that $|t-t_i|\ll e^{-r\theta \lfloor t\rfloor/4}$, and 
then for every $f\in L^\infty(\mu)$,
$$
|\pi_X(\lambda_t)f-\pi_X(\lambda_{t_i})f|\le e^{-\eta t}\|f\|_{L^\infty(\mu)}
$$
with some $\eta>0$, independent of $t$ and given explicitly below.
\end{claim}

To prove the claim, note that 
it follows from (\ref{eq:t_i}) 
that $t_i$ satisfying the first property exists.
Without loss of generality we suppose that $t>t_i$. Then, 
\begin{align*}
&\quad\pi_X(\lambda_t)f(x)-\pi_X(\lambda_{t_i})f(x)\\
&=\frac{1}{|\Gamma_t|\cdot|\Gamma_{t_i}|}
\left(|\Gamma_{t_i}|\sum_{\gamma\in\Gamma_t} f(\gamma^{-1}\cdot x)- |\Gamma_{t}|\sum_{\gamma\in\Gamma_{t_i}} f(\gamma^{-1}\cdot x)\right)\\
&= -\frac{|\Gamma_t-\Gamma_{t_i}|}{|\Gamma_t|\cdot|\Gamma_{t_i}|}\sum_{\gamma\in\Gamma_{t_i}} f(\gamma^{-1}\cdot x)
+\frac{1}{|\Gamma_t|}\sum_{\gamma\in\Gamma_t-\Gamma_{t_i}} 
f(\gamma^{-1}\cdot x)\\
&\le \frac{2|\Gamma_t-\Gamma_{t_i}|}{|\Gamma_t|}\|f\|_{L^\infty(\mu)}.
\end{align*}
Applying the estimate provided by Theorem \ref{l:gG2} and (\ref{eq:2}), we get
\begin{align*}
\frac{|\Gamma_t-\Gamma_{t_i}|}{|\Gamma_t|}&=1-\frac{|\Gamma_{t_i}|}{|\Gamma_t|}
=1-\frac{1+O(e^{-\alpha t_i})}{1+O(e^{-\alpha t})}\cdot\frac{m_G(G_{t_i})}{m_G(G_t)}\\
&\le 1- (1+O(e^{-\alpha t}))\cdot \frac{1}{1+c(t-t_i)}.
\end{align*}
This implies the claim. \qed

\medskip

Continuing with the proof of Theorem \ref{th:exp.max}, we use 
 (\ref{eq:exp_proof2}) and the Claim, and 
deduce that for $\delta=\min\{\theta/2,\eta\}$ 
and every $f\in L_0^\infty(\mu)$,
\begin{align*}
|\pi_X(\lambda_t)f(x)|&\le |\pi_X(\lambda_{t_i})f(x)|
+|\pi_X(\lambda_t)f(x)-\pi_X(\lambda_{t_i})f(x)|\\
&= \left(B_r(x,f)+\|f\|_{L^\infty(\mu)}\right)e^{-\delta t}.
\end{align*}
Hence, by (\ref{eq:exp_proof2}), (\ref{eq:exp_proof3})
 for some $t_0, C>0$, and every $f\in L_0^\infty(\mu)$ 
\begin{equation}\label{eq:exp_proof4}
\left\| \sup_{t\ge t_0} e^{\delta t} |\pi_X(\lambda_t)f|\right\|_{L^r(\mu)}\le C \|f\|_{L^\infty(\mu)}\,\,
\end{equation}
where we have used the fact that since $\mu$ 
is a probability measure, for $f\in L^\infty(\mu)$, 
$\norm{f}_{L^r(\mu)}\le \norm{f}_{L^\infty(\mu)}$.  

Now for a measurable function
 $\tau: X\to [t_0,\infty)$ and 
$z\in \mathbb{C}$, we consider the linear operator
$$
U_z^\tau f(x)= e^{z\delta\tau(x)}\left(\pi_X(\lambda_{\tau(x)})f(x)-\int_X f\,d\mu\right).
$$
By the strong maximal inequality for $\lambda_t$ in $L^q(\mu)$ (which holds 
using Theorem  \ref{reduction}(1) and our second assumption), 
when $\hbox{Re}\, z=0$, the operator  
$$U_z^\tau: L^q(\mu)\to L^q(\mu)$$
 is bounded. By (\ref{eq:exp_proof4}), when $\hbox{Re}\, z=1$, the operator  
$$U_z^\tau: L^\infty(\mu)\to L^r(\mu)$$
 is also bounded, 
 with bounds independent of the function $\tau$. 
Hence, by the complex interpolation theorem (see e.g. \cite{MNS} for a 
fuller discussion) 
for every $u\in (0,1)$ and $v$, $w$ such that 
$$
1/v=(1-u)/q\quad\hbox{and}\quad 1/w=(1-u)/q+u/r,
$$
we have 
$$
\left\| \sup_{t\ge t_0} 
e^{u\delta t} \left|\pi_X(\lambda_t)f-
\int_X f\,d\mu\right|\right\|_{L^w(\mu)}\le C \|f\|_{L^v(\mu)}.
$$
This completes the proof of Theorem \ref{th:exp.max}. 
\end{proof}


\subsection{Completion of proofs of ergodic theorems for 
lattices}\text{                             }

{\it 1) Completion of the proof 
of Theorem \ref{th:semisimple_lattice-1}}.

Clearly, parts (1), (2) and (3)
 of Theorem \ref{reduction} together imply 
Theorem \ref{th:semisimple_lattice-1}, 
provided only that the admissible 
averages $\beta_t$ do indeed satisfy the 
mean, maximal and pointwise 
ergodic theorems in $L^p(\nu)$ (and as a result also in  
$L^p(m_{G/\Gamma})$). This follows 
immediately from Theorem \ref{th:semisimple-1}, taking also 
into account the fact that since we are considering the action 
induced to $G^+$, the action 
is necessarily totally weak-mixing, 
since $G^+$ has no non-trivial finite-dimension 
representations.\qed

 {\it 2) Completion of the proof of Theorem 
\ref{th:semisimple_lattice-T}}.

The formulation of Theorem \ref{reduction}(4) 
 incorporates the assumption that 
$p \ge r > \varrho_0$, where 
$\varrho_0$ is the upper local dimension.
Thus in order to complete the proof of Theorem \ref{th:semisimple_lattice-T} 
we must remove this restriction.

Let us consider the exponentially fast mean ergodic theorem 
in $(L^p, L^r)$ 
first. By Theorem \ref{th:semisimple-T} $\beta_t$ 
on $G$ satisfies this theorem for all  $p=r> 1$, 
and hence by Theorem \ref{th:mean_exp}, we obtain that 
$\lambda_t$ satisfies it if $p > \varrho $. But clearly 
$\lambda_t-\int_X d\mu$ has norm bounded by $2$ in every 
$L^p$, $1 \le p \le \infty$. By Riesz-Thorin interpolation, 
it follows that that 
$\lambda_t$ satisfies the exponentially fast 
mean ergodic theorem in every $L^p$, $1 < p < \infty$, and 
hence in $(L^p,L^r)$, $p\ge r\ge 1$, $(p,r)\neq (1,1)$.

As to the exponential maximal inequality, 
note first that by Theorem  \ref{th:semisimple-T} 
$\beta_t$ satisfies the 
$(L^\infty,L^2)$ exponential-maximal inequality in every action 
of $G$ where $\norm{\pi(\beta_t)}_{L^2}\le C\exp(-\theta t)$. 
It follows that it satisfies the 
exponential-maximal inequality 
in $(L^\infty, L^r)$, for a finite $r> \varrho_0$.
By Theorem \ref{ample sets}, $\beta_t$ also satisfy the standard 
strong maximal inequality in every $L^q$, $q > 1$. 
Thus by Theorem \ref{th:exp.max}
the exponential maximal inequality in $(L^v,L^w)$
 holds for the averages $\lambda_t$ 
(provided the norm exponential decay condition holds 
in the induced action, which is the case under our assumptions). 
By their explicit formula it is clear that we can choose 
$v$ to be as close as we like to $1$, thus determining a 
 consequent $w < v$ and some positive rate of exponential  
decay.

This completes the proof of Theorem 
\ref{th:semisimple_lattice-T}. \qed

Finally, we remark that the ergodic theorems stated for 
connected semisimple Lie groups and their lattices in 
\S \S 1.2, 1.3, and 1.4 all follow from 
Theorem \ref{th:semisimple-1}, Theorem \ref{th:semisimple-T}
Theorem \ref{th:semisimple_lattice-1} 
and 
Theorem \ref{th:semisimple_lattice-T}, 
together with Theorem \ref{l:gG2} (or more precisely 
Corollary \ref{Lie-error} below). 
This is a straightforward 
verification, bearing in mind that every unitary representation 
of a connected semisimple Lie group with finite center is
 totally weak-mixing.
 
 The only comment necessary  is regarding Theorem 
\ref{killing}, which asserts that $\lambda_t$ defined using the Riemannnian averages associated with the Killing form satisfy 
the pointwise ergodic action in any ergodic action of any 
lattice $\Gamma$, even if the induced action is reducible. Theorem 
\ref{th:semisimple_lattice-1} establishes this result based  
 on Theorem \ref{th:semisimple-1}, since the Riemannian averages are indeed admissible and well-balanced. This fact follows from the discussion in \cite{MNS} (see also the Appendix below).   
It follows that an exponential  decay estimate holds for the Riemannian sphere averages $\norm{\pi(\partial\nu_t)}$, even 
  in the case of a reducible action, and thus pointwise convergence on a dense subspace holds. 

An alternative, more direct argument for the latter conclusion is the fact that in  
\cite[Thm. 1]{MNS} the pointwise ergodic theorem for $\beta_t$ is proved in 
full generality, even for reducible actions, and thus 
the result for $\lambda_t$ follows from Theorem 
\ref{reduction}(3).

In might be worth commenting that  for the {\it sequences} of averages $\beta_n$ and $\lambda_n$   
the entire Sobolev space argument is superfluous, of course. Since the maximal inequality holds, as well as pointwise convergence on dense subspace, it follows that the pointwise 
ergodic theorem holds for both sequences, even in the reducible case.

This concludes the proofs of the ergodic theorems for lattice actions. \qed
 
 We now turn to discuss equidistribution.

\subsection{Equidistribution in isometric actions}

Let us prove the following 
generalization of Theorem \ref{th:equidistribution}.
\begin{theorem}\label{equi-char-0}

Let $G$ be an $S$-algebraic group as in 
Definition \ref{alg-gps}, over fields of characteristic zero.
Let $\Gamma\subset G^+$ be a lattice and $G_t\subset G^+$ an 
admissible $1$-parameter family or sequence.  
Let $(S,m)$ be an isometric action of $\Gamma$ on 
a compact metric space $S$, 
preserving an ergodic probability measure $m$ of full support. 
If $\Gamma$ is an irreducible lattice, then 
$$\lim_{n\to \infty}\max_{s\in S} \abs{\pi_S(\lambda_t)
 f(s)-\int_S f dm}=0$$
 and in particular 
$\pi_S(\beta_t)f(s)\to \int_S fdm $ for {\it every} 
$s\in S$.  

For any $S$-algebraic  $G$ and lattice $\Gamma$, the same results holds provided 
that $G_t$ are left-radial and balanced. 

\end{theorem}
\begin{proof} 
When $G$ is defined over fields of characteristic zero, 
an action of $G$ induced by an isometric ergodic action of an 
irreducible lattice is an irreducible action of $G$, as shown 
in \cite{St}.
 Then, according to Theorem \ref{th:semisimple_lattice-1}, 
$\lambda_t$ satisfies the mean ergodic theorem in $L^2(S,m)$. 
The mean ergodic theorem holds also 
for general $G$ and $\Gamma$, provided 
only that $\lambda_t$ is balanced and left-radial, as asserted in Theorem \ref{th:semisimple_lattice-1}. 

The proof is therefore complete using the following proposition.
\end{proof}

\begin{proposition}\label{mean-equi}(see \cite{G})
Let the discrete group $\Gamma$ act isometrically on a compact metric space 
$(S,m)$, preserving a probability measure of full support. If a 
$1$-parameter family (or 
sequence) of averages $\lambda_t$ on 
$\Gamma$ satisfy the mean ergodic theorem 
in $L^2(S,m)$, 
then $\pi_S(\lambda_t)f$ converges uniformly to the constant 
$\int_S fdm$ for every continuous $f\in C(S)$. 
\end{proposition}
\begin{proof}
Our proof is a straightforward 
 generalization of \cite{G} (where the case 
of free groups is considered), and is brought 
here for completeness.
 
Given a function $f\in C(S)$, consider the set $ C(f)$ 
consisting of $f$ together with 
$\pi_S(\lambda_t)f$, $t \in 
\RR_+$. $C(f)$ constitutes 
an equicontinuous family of functions, since $\Gamma$ acts 
isometrically on $S$.
 Thus $C(f)$ has compact 
closure in $C(S)$, w.r.t. the uniform 
norm. Let $f_0\in C(S)$ be the 
uniform limit of $\pi_S(\lambda_{t_i})f$ for some 
subsequence $t_i\to \infty$. Then $f_0$ is of course also the limit of 
$\pi_S(\lambda_{t_i})f$ in the $L^2(S,m)$-norm. Given that  $\lambda_t$ 
satisfy the mean ergodic theorem, it follows that 
$f_0(s)=\int_S fdm $ for $m$-almost all  
$s\in S$. Since $m$ has full support, the last equality holds on a dense 
subset of $S$, and since $f_0$ is continuous, it must hold everywhere.
Thus $\pi_S(\lambda_{t_i})f$ converges uniformly to the constant $\int_S fdm$, 
and this holds for every subsequence $t_i\to \infty$. It follows that 
the latter constant in the unique limit point in the uniform closure 
of the family $\pi_S(\lambda_t)f$. Hence 
$$\lim_{t\to\infty}\norm{\pi_S(\lambda_t)f-\int_S fdm }_{C(S)}= 0\,\,.$$

\end{proof}

\section{Comments and complements}

\subsection{Explicit error term} 

Let us start by stating the following error estimate.

\begin{corollary}\label{Lie-error}
Let $G$ be an $S$-algebraic group as in Definition 
\ref{alg-gps}. Let 
$\Gamma$ be a lattice subgroup 
and $G_t$ an admissible 
family, both contained in $G^+$. If $G_t$ are well-balanced, 
or the lattice is irreducible, then 
 the number of lattice points in $G_t\cap \Gamma$ 
is estimated by (for all $\vre > 0$)
$$
\abs{\frac{|\Gamma_t|}{m_G(G_t)}-1}\le B_\vre 
\exp\left(-\frac{(\theta-\vre) t}{\varrho_0+1}\right)\,\,,
$$ 
where $\theta$ need only satisfy (for all $\vre > 0$)
$$\norm{\pi_{G/\Gamma}(\beta_t)}_{L^2_0(G/\Gamma)}\le 
C_\vre e^{-(\theta-\vre) t}\,\,. $$
If $\pi=\pi^0_{G/\Gamma}$ has a strong spectral gap, then  
$\pi^{\otimes n}\subset \infty \cdot\lambda_G$ for some $n$, and 
the spectral parameter $\theta$ is given explicitly in term of the rate of volume growth of $G_t$ by   
$$\theta=\frac{1}{2n}\limsup_{t\to \infty}\frac1t \log m_G(B_t)\,.$$
\end{corollary}
We note that Corollary \ref{Lie-error} is an 
immediate corollary of Theorem \ref{l:gG2}, 
taking $r=p=2$ and $\varrho_0$ to be the upper local dimension.
In addition one uses Remark \ref{KS}, and the fact that 
 $S$-arithmetic groups as in Definition \ref{alg-gps} have the 
Kunze-Stein property. This is well known in the real case 
\cite{Co1}, and was 
proved by A. Veca \cite[Thm. 1]{V} in the totally disconnected, simply connected case. 

 We remark that a 
 somewhat weaker rate can be established using (in effect)  just the radial 
Kunze-Stein phenomenon (which is much easier to establish). This proceeds  by using the fact that admissible 
averages are $(K,C)$-radial, and estimating the norm 
of the corresponding radialized averages directly using 
Proposition  \ref{radial estimate}(2) and the standard estimate of the $\Xi$-function.

Of course, $\varrho_0=\dim_\RR G$ when 
$G$ is a connected Lie group and 
$\mathcal{O}_\vre$ are Riemannian balls. Furthermore, 
$\varrho_0=0$ when $G$ is a totally disconnected $S$-algebraic group.

\begin{remark}{\bf Lattice point counting problem.}
 Corollary \ref{Lie-error} constitutes a quantitative solution to the lattice point counting problem in admissible domains. For a systematic discussion of quantitative counting results for more general domains, and more general groups, together with many applications,   we refer to \cite{GN}. 
\end{remark}

%

%
%
%
%
%



\subsection{Exponentially fast convergence versus 
equidistribution}

In this section we give an example of a connected semisimple Lie group $H$
without compact factors acting by translations on a homogeneous space $G/\Gamma$
of finite volume and Haar uniform admissible averages $\beta_t$ 
on $H$ such that
\begin{itemize}
\item equidistribution of $H$-orbits fails (i.e., there exist dense
  orbits for which the averages do not converge to the Haar
  measure),
\item An exponentially fast pointwise ergodic theorem holds
(i.e., for almost all starting points the averages converge to the Haar
  measure exponentially fast).
\end{itemize}

This example was originally constructed in [GW], Section 12.3.

Let
$$
H=H_1\times H_2=\hbox{SL}(2,\mathbb{R})\times\hbox{SL}(2,\mathbb{R})
$$
and
$$
r:H\to\hbox{SL}(2l,\mathbb{R})
$$
be the representation of $H$ which is a tensor product of irreducible
representations of $H_1$ and $H_2$ of dimensions $2$ and $l> 2$
respectively. We fix a norm $\|\cdot\|$ on $\hbox{M}(2l,\mathbb{R})$
and define
$$
H_t=\{h\in H:\, \|r(h)\|<e^t\}.
$$
Note that the sets $H_t$ are not balanced, as shown in \cite{GW}.

Let $G=\hbox{SL}(2l,\mathbb{R})$ and
$\Gamma=\hbox{SL}(2l,\mathbb{Z})$.
For $x\in G/\Gamma$ and $t>0$, consider the Radon probability measure
$$
\mu_{x,t}(f)=\frac{1}{m_H(H_t)} \int_{H_t} f(h^{-1}x)dm_H(h),\quad f\in
C_c(G/\Gamma).
$$

\begin{proposition}
\begin{enumerate}
\item There exists $x\in G/\Gamma$ such that $\overline{Hx}=G/\Gamma$, but
 Haar measure $m_{G/\Gamma}$ is not an accumulation point of the sequence
 $\mu_{x,t}$, $t\to \infty$, in the weak$^*$ topology.

\item For a.e. $x\in G/\Gamma$, $\mu_{x,t}\to m_{G/\Gamma}$ as
  $t\to\infty$ in the weak$^*$ topology.
Moreover, for every $p>r\ge 1$, there exists $\theta>0$ such that
for $f\in L^p(m_{G/\Gamma})$ and a.e. $x\in G/\Gamma$,
$$
\left|\mu_{x,t}(f)-\int_{G/\Gamma} f\,d\mu_{G/\Gamma}\right|\le
C(x,f)e^{-\theta_{p,r} t}
$$
with
$$
\|C(\cdot, f)\|_{L^r(\mu_{G/\Gamma})}\le C\|f\|_{L^p(\mu_{G/\Gamma})}.
$$
\end{enumerate}
\end{proposition}

\begin{proof}
Part (1) was proved in [GW], Section 12.3.

To prove part (2), it suffices to observe that the representation of $G$ on
$L^2_0(G/\Gamma)$ has spectral gap. Being simple, the spectral gap is strong, 
and so some 
tensor power of the 
representation embeds in $\infty\cdot \lambda_G$, as 
 follows from the spectral
transfer principle (see Theorem \ref{tensor}, or \cite{N3}). 
Thus the same tensor power of 
  the representation of $H$ on
$L^2_0(G/\Gamma)$ (restricted to $H$) 
embeds in $\infty \cdot \lambda_H$ and so $H$ 
has a strong spectral gap as well. Therefore the desire result follows 
from Theorem \ref{th:semisimple-T}. 
\end{proof}

\subsubsection{Remark about balanced sets}

The last example owes its existence to the fact that the averages 
considered are not balanced. Let us therefore give an easy 
geometric criterion for a family of sets defined by
a matrix norm on a product of simple groups to be balanced. 

Let $G=G_1\cdots G_s$ be a connected semisimple Lie group where $G_i$'s
are the simple factors and
\begin{equation}\label{eq:a_decomp}
\mathfrak{a}=\mathfrak{a}_1\oplus\cdots \oplus \mathfrak{a}_s
\end{equation}
a Cartan subalgebra of $G$ where $\mathfrak{a}_i$'s are Cartan
subalgebras of $G_i$'s. We fix a system of simple roots
$\Phi=\Phi_1\cup\cdots\cup\Phi_s$ for $\mathfrak{a}$ where $\Phi_i$ is a system of simple
roots for $\mathfrak{a}_i$ and denote by
$$
\mathfrak{a}^+=\mathfrak{a}_1^+\oplus\cdots \oplus \mathfrak{a}_s^+
$$
the corresponding positive Weyl chamber.

Let $r:G\to \hbox{GL}(d,\mathbb{R})$ be a representaion of $G$.
For a norm $\|\cdot\|$ on $M_d(\mathbb{R})$, let 
$$
G_t=\{g\in G:\, \|r(g)\|<t\}.
$$
Let $\Psi_r$ be the set of weights of $\mathfrak{a}$, and 
$$
\mathfrak{p}_r=\{H\in\mathfrak{a}^+:\, \lambda(H)\le
1\quad\hbox{for all $\lambda\in\Psi_r$}\},
$$
Finally, let 
$$
\delta=\max\{\rho(H):\, H\in \mathfrak{p}_r\}
$$
where $\rho$ denotes the half sum of the positive roots of
$\mathfrak{a}$.

We can now formulate the following 
\begin{proposition}
The sets $G_t$ are balanced iff the set $\{\rho=\delta\}\cap
  \mathfrak{p}_\rho$ is not contained in any proper subsum of the
  direct sum (\ref{eq:a_decomp}).
\end{proposition}
The Proposition follows from [GW], Section 7.

\section{Appendix : volume estimates and volume regularity}

The appendix is devoted to establishing  
admissibility or H\"older-admissibility of the standard radial averages and more general ones, as well as to establishing conditions sufficient for the averages to be  balanced or 
 well balanced. We will also discuss boundary-regularity and differentiability properties of volume functions for some metrics, particularly $CAT(0)$-metrics. 
 
 \subsection{Admissibility of standard radial averages}
 We begin with a proof of  Theorem \ref{Appendix:vol}, whose statement we recall. 
\vskip0.1cm

\textbf{Theorem 3.14.}  {\it For  an $S$-algebraic group $G=G(1)\cdots G(N)$ as in Definition \ref{alg-gps}, the following families of sets $G_t\subset G$ are admissible, where $a_i$ are any positive constants.
\begin{enumerate}
\item Let $S$ consist of infinite places, and let $G(i)$ be a closed
subgroup of the isometry group of a symmetric space $X_i$ of nonpositive curvature
equipped with the Cartan--Killing metric. For $u_i,v_i\in  X_i$, define
$$
G_t=\{(g_1,\ldots,g_N): \sum_i   a_i d_i(u_i, g_i\cdot v_i)<t\}.
$$
\item Let $S$ consist of infinite places, and let $\rho_i:G(i)\to \hbox{\rm GL}(V_i)$
be proper rational representations. For norms $\|\cdot \|_i$ on $\hbox{\rm End}(V_i)$, define
$$
G_t=\{(g_1,\ldots,g_N): \sum_i a_i \log \|\rho_i(g_i)\|_i<t\}.
$$
\item For infinite places, let $X_i$ be the symmetric space of $G(i)$ equipped with the Cartan-Killing
 distance $d_i$, and  for finite places, let $X_i$ be 
the Bruhat-Tits building of $G(i)$ equipped with the path metric $d_i$ on its $1$-skeleton. 
For $u_i\in  X_i$, define
$$
G_t=\{(g_1,\ldots,g_N): \sum_i a_i d_i(u_i, g_i\cdot u_i)<t\}.
$$
\item Let $\rho_i:G(i)\to \hbox{\rm GL}(V_i)$
be proper representations, rational over the fields of definition $F_i$. For infinite places, let
$\|\cdot \|_i$ be a Euclidean norm on $End(V_i)$, 
 and assume 
that $\rho_i(G(i))$ is self-adjoint : $\rho_i(G(i))^t=\rho_i(G(i))$.  
For finite places, let $\|\cdot \|_i$ be the $\max$-norm on $\hbox{\rm End}(V_i)$.
Define
$$
G_t=\{(g_1,\ldots,g_N): \sum_i  a_i \log \|\rho_i(g_i)\|_i<t\}.
$$
\end{enumerate}
}


The proof is divided into several propositions.
To handle Archimedian groups we will employ some arguments originating in \cite{DRS} and \cite{EMS}, 
and in the general case of $S$-algebraic groups we will also employ 
convolution arguments which will  be developed below. We note that the latter arguments will utilize knowledge of the behavior of $\vol (B_t)$ for {\it all } $t > 0$, and we will thus consider below the behavior for $t$ large and for $t$ small,  separately.

\begin{proposition}\label{p:metric}
Let $G$ be a connected semisimple group with finite center
and $X$ the corresponding symmetric space equipped with the
Cartan--Killing metric $d$. For $u\in X$, set
$$
G_t=\{g\in G;\, d(u,g\cdot u)<t\}.
$$
Then there exists $c>0$ such that for all $t\ge 0$ and $\epsilon\in (0,1)$,
$$
\vol(G_{t+\epsilon})-\vol(G_t) \le c\epsilon\max\{1, \vol(G_t)\}.
$$
\end{proposition}

\begin{proof}
Note that the stabilizer of $u$ is a maximal compact subgroup $K$ of $G$,
and for a Cartan subgroup $A$ of $G$, the map $a\mapsto a\cdot u$, $a\in A$,
is an isometry. We introduce polar coordinates $(r,\omega)\in\mathbb{R}^+\times S^+$ on the Lie
algebra of $A$. With respect to the Cartan decomposition $G=KA^+K$, the Haar
measure on $G$ is given by $\xi(r,\omega)\,drd\omega dk$ with a nonnegative smooth density function $\xi$.
Then
$$
\vol(G_{t+\epsilon})-\vol(G_t)=\int_{S^+}\int_t^{t+\epsilon} \xi(r,\omega)\,drd\omega=
\epsilon\int_{S^+} \xi(\sigma(\omega),\omega)\,d\omega
$$
for some $\sigma(\omega)\in [t,t+\epsilon]$. This implies the claim for $t\in [0,1]$.
To establish the claim for $t>1$, we use the following property of the function $\xi$ (see
\cite{EMS}, Lemma A.3): there exists $c>0$ such
that for every $r>1$ and $\omega\in S^+$,
\begin{equation}\label{eq:ems}
\xi(r,\omega)\le c\int_0^r\xi(s,\omega)ds.
\end{equation}
\end{proof}

\begin{proposition}\label{p:metric2}
Let $d$ be the Cartan--Killing metric on a symmetric space $X$ of
nonpositive curvature, $u,v\in X$, and $G$ a closed connected
semisimple subgroup of the isometry group of $X$.
Define the sets
$$
G_t=\{g\in G;\, d(u,g\cdot v)<t\}.
$$
Then there exist $c,t_0>0$ such that for all $t>t_0$ and
$\epsilon\in (0,1)$,
$$
\vol(G_{t+\epsilon})- \vol(G_t) \le c\epsilon\, \vol(G_t).
$$
\end{proposition}

\begin{proof}
It follows from Mostow's theorem that there exist a maximal compact subgroup $K$ of
$G$ and an associated Cartan subgroup $A$ such that
$K\subset\hbox{Stab}(v)$ and the map $a\mapsto a\cdot v$, $a\in A$, is
an isometry. Consider polar coordinates
$(r,\omega)\in\mathbb{R}^+\times S^+$ on the Lie
algebra of $A$, and set
$$
S_t(k)=\{ (r,\omega);\, d(k^{-1}u,\exp(r\omega)v)<t\}.
$$
Then
$$
\vol(G_t)=\int_K\int_{S_t(k)} \xi(r,\omega)\,dr d\omega dk.
$$
For $p\in Ku$, we consider the function 
$$
f_p(x)=\frac{1}{2}d(p,x)^2,\quad x\in X.
$$
We have
$$
\hbox{grad}\, f_p=-\exp^{-1}_x(p),
$$
and for the unit-speed geodesic ray $\gamma_\omega(r)=\exp(r\omega)v$,
\begin{align*}
\frac{d}{dr}f_p(\gamma_\omega(r)) &=\langle
(\hbox{grad}\,f_p)_{\gamma_\omega(r)},\gamma'_\omega(r)\rangle_{\gamma_\omega(r)},\\
\frac{d^2}{dr^2}f_p(\gamma_\omega(r))&=\langle
\grad_{\gamma'_\omega(r)}
(\hbox{grad}\,f_p)_{\gamma_\omega(r)},\gamma'_\omega(r)\rangle_{\gamma_\omega(r)}.
\end{align*}
Since the space $X$ has nonpositive sectional curvature, we deduce
(see \cite{j}, Theorem 4.6.1) that
$$
\langle \grad_w (\hbox{grad}\,f_p)_x,w\rangle_x\ge
\|w\|_x^2\quad\hbox{for every $w\in T_x X$}.
$$
Hence,
\begin{equation}\label{eq:deriv}
\frac{d^2}{dr^2}f_p(\gamma_\omega(r))\ge 1.
\end{equation}
Therefore, there exists $r_0>0$ and $\alpha>0$ such that for every $r>r_0$, $p\in
Ku$, $\omega\in S^d$, we have
$$
\frac{d}{dr}f_p(\gamma_\omega(r))\ge \alpha r.
$$
For $\epsilon>0$ and $r>r_0$,
\begin{align*}
d(p,\exp((r+\epsilon)\omega)v)&=\sqrt{2f_p(\gamma_\omega(r+\epsilon))}\ge
\sqrt{2f_p(\gamma_\omega(r))+2\alpha r\epsilon}\\
&=d(p,\exp(r\omega)v)\sqrt{1+2\alpha r\epsilon/d(p,\gamma_\omega(r))^2}.
\end{align*}
Since it follows from the triangle inequality that for some $c>0$,
$$
r-c\le d(p,\gamma_\omega(r))\le r+c\quad\hbox{for all $p\in Ku$ and
  $\omega\in S^+$},
$$
we conclude that for some $\beta>0$,
\begin{equation}\label{eq:beta}
d(p,\exp((r+\epsilon)\omega)v)\ge d(p,\exp(r\omega)v)+\beta\epsilon.
\end{equation}
for sufficiently large $r$ and sufficiently small $\epsilon>0$.
This implies that for sufficiently large $t$, the sets
$S_t(k)$ are star-shaped. 
Let $r_t(k,\omega)$ be the unique solution of the equation.
$$
d(k^{-1}u,\exp(r\omega)v)=t.
$$
Note that, by (\ref{eq:beta}), for sufficiently large $t$ and $\epsilon \in (0,1)$,
$$
r_{t+\epsilon}(k,\omega)\le r_t(k,\omega)+\beta^{-1}\epsilon.
$$
Setting
$$
m_t(k,\omega)=\int_0^{r_t(k,\omega)} \xi(r,\omega) dr,
$$
we have
\begin{align*}
m_{t+\epsilon}(k,\omega)-m_t(k,\omega)
&\le
\int_0^{r_t(k,\omega)+\beta^{-1}\epsilon} \xi(r,\omega)r
dr
-\int_0^{r_t(k,\omega)} \xi(r,\omega) dr\\
&=\beta^{-1}\epsilon \cdot \xi(\sigma,\omega).
\end{align*}
for some $\sigma\in [r_t(k,\omega),r_t(k,\omega)+\beta^{-1}\epsilon]$.
Now it follows from (\ref{eq:ems}) that
$$
m_{t+\epsilon}(k,\omega)-m_t(k,\omega)\le (\beta^{-1}c)\epsilon \cdot m_{t+\epsilon}(k,\omega).
$$
This shows that
$$
\vol(G_{t+\epsilon})
=\int_K\int_{S^d}m_{t+\epsilon}(k,\omega)d\omega dk
\le (\beta^{-1}c)\epsilon\cdot \vol(G_{t+\epsilon})+ \vol(G_t),
$$
which implies the proposition.
\end{proof}

\begin{proposition}(cf \cite[Appendix]{EMS})\label{p:e}
Let $\rho:G\to \hbox{\rm GL}(V)$ is a proper representation of 
a connected semisimple Lie group $G$,
$\|\cdot\|$ a norm on $\hbox{\rm End}(V)$, and
$$
G_t=\{g\in G;\, \log\|\rho(g)\|< t\}.
$$
Then there exist $c,t_0>0$ such that for all $t>t_0$ and
all $\epsilon\in (0,1)$,
$$
\vol(G_{t+\epsilon})- \vol(G_t) \le c\epsilon\, \vol(G_t).
$$
\end{proposition}

\begin{proof}
We employ the argument from Appendix of \cite{EMS}, but since this argument is not quite complete (see (\ref{eq:eskin2}) below), we provide a sketch which indicates that in our setting it is
indeed applicable and provides the Lipschitz estimate.

We fix a Cartan decomposition $G=KA^+K$ and use polar coordinates
$(r,\omega)$ on the Lie algebra of $A$. Let
$$
S_t(k_1,k_2)=\{(r,\omega);\,\|\rho(k_1 \exp(r\omega) k_2)\|< e^t\}.
$$
Then
$$
\vol(G_t)=\int_{K\times K} \int_{S_t(k_1,k_2)} \xi(r,\omega) drd\omega dk_1dk_2.
$$
Since $\rho(A)$ is (simultaneously) diagonalizable over $\mathbb{R}$, there exist
$v_i\in\hbox{End}(\mathbb{R}^n)$
such that
$$
\rho(\exp(r\omega))=\sum_i e^{r\lambda_i(\omega)}v_i.
$$
Since all norms are equivalent and $\rho(K)$ is compact, there exist $c_1,c_2>0$ such that for
every $k_1,k_2\in K$, $r>0$, and $\omega\in S^+$, 
\begin{equation}\label{eq:eskin2}
c_1 \exp(r\max_i\lambda_i(\omega))\le \|\rho(k_1 \exp(r\omega)
k_2)\|\le c_2 \exp(r\max_i\lambda_i(\omega)).
\end{equation}
Note that the lower estimate needs further argument in the  generality of \cite{EMS},  since there the $v_i$'s depend on $\omega$. But for constant $v_i$, 
using these estimates, the argument from \cite{EMS}, Lemma A.4, shows
that 
\begin{enumerate}
\item There exists $t_0>0$ such that for $t>t_0$ the sets
$S_t(k_1,k_2)$ are star-shaped.
\item There exists $r_0>0$ such that for
every $r>r_0$, $\omega\in S^d$, and $\epsilon\in [0,1)$, we have
$$
\|\rho(k_1 \exp((r+\epsilon)\omega) k_2)\|\ge g(\epsilon)\cdot \|\rho(k_1 \exp(r\omega)
k_2)\|
$$
where $g:[0,1)\to [1,\infty)$ is explicit smooth function such that
$g(0)=1$ and $g'>0$. In particular, there exists $\beta>0$ such that
$g(\epsilon)\ge e^{\beta\epsilon}$.
\end{enumerate}

Let $r_t(k_1,k_2,\omega)$ denote the unique solution of the equation
$$
\|\rho(k_1 \exp(r\omega) k_2)\|= e^t.
$$ 
Then it follows that  
$$
r_{t+\epsilon}(k_1,k_2,\omega)\le r_t(k_1,k_2,\omega)+\beta^{-1}\epsilon.
$$
Finally, the Lipschitz property of the sets $G_t$ can be proved as in Proposition \ref{p:metric} above.
\end{proof}

Let us note the following consequence of the foregoing arguments. Let $H$ be a symmetric subgroup of a connected semisimple Lie group with finite center, embedded a Zariski closed $G$-orbit in a linear space $V$.  
\begin{corollary}
Proposition \ref{p:e}  applies to subsets of affine symmetric varieties $G/H$ defined by an arbitrary norm on the ambient vector space.
\end{corollary}
\begin{proof}
 Indeed, for symmetric varieties one has a decomposition of the form  $G=KAH$ where $A$ is simultaneously diagonalizable, and the arguments utilized in the proof of Proposition \ref{p:e} apply without any material changes.  
\end{proof}

\begin{proposition}\label{p:e2}
Let $\rho:G\to \hbox{\rm GL}(V)$ is a proper representation of 
a connected semisimple Lie group $G$ such that ${}^t\rho(G)=\rho(G)$,
$\|\cdot\|$ the Euclidean norm on $\hbox{\rm End}(V)$, and
$$
G_t=\{g\in G;\, \log\|\rho(g)\|< t\}.
$$
Then there exist $c>0$ such that for all $t>0$ and
all $\epsilon\in (0,1)$,
$$
\vol(G_{t+\epsilon})- \vol(G_t) \le c\epsilon\, \max\{1,\vol(G_t)\}.
$$
\end{proposition}

\begin{proof}
Let $t_0$ be as in Proposition \ref{p:e}. It remains to prove the claim for $t\le t_0$.

Since $\rho(G)$ is self-adjoint, there exist a maximal compact subgroup $K$ such that
$\rho(K)\subset \hbox{SO}(V)$ and a Cartan subgroup such that $\rho(A)$ is diagonal.
For $k_1,k_2\in K$ and $a\in \hbox{Lie}(A)$,
$$
\|\rho(k_1\exp(a)k_2)\|^2=\sum_i e^{2\lambda_i(a)}.
$$
where $\lambda_i$'s are characters of $\hbox{Lie}(A)$
such that $\sum_i\lambda_i=0$. We use polar coordinates $(r,\omega)$ on $A$
and set
$$
f_\omega(r)=\sum_i e^{2r\lambda_i(\omega)}.
$$
Then
$$
\vol(G_t)=\int_{(r,\omega): \log f_\omega(r)<2t} \xi(r,\omega)\, drd\omega.
$$
Since $f_\omega''>0$ and $f_\omega'(0)=0$, the function $\log f_\omega$ is increasing.
Let $r_\omega$ be the inverse function of $\log f_\omega$ and $r_0=\max\{r_\omega(2t_0)\}$.
By the mean value theorem, there exists $\alpha>0$ such that
$$
f_\omega'(r)\ge \alpha r\quad\hbox{for $r\in [0,r_0]$.}
$$
Then for some $\beta>0$,
\begin{equation}\label{eq:der}
r_\omega'(t)\le \beta r_\omega(t)^{-1}\quad\hbox{for $r\in [0,2t_0]$.}
\end{equation}
Since for some $c>0$,
$$
\xi(r,\omega)\le c\, r\quad\hbox{for $r\in [0,r_0]$ and $\omega\in S^d$,}
$$
we have
\begin{align*}
\vol(G_{t+\epsilon})-\vol(G_t)&=\int_{S^d} \int_{r_\omega(2t)}^{r_\omega(2t+2\epsilon)}\xi(r,\omega)\,drd\omega\\
&\le c\int_{S^d}(r_\omega(2t+2\epsilon)^2-r_\omega(2t)^2)\, d\omega.
\end{align*}
Now the proposition follows from (\ref{eq:der}).
\end{proof}

\subsection{Convolution arguments}
We now turn to discuss convolution arguments, which together with the foregoing result will complete 
the proof of Theorem \ref{Appendix:vol}. 
Let $G_i$, $i=1,2$, be locally compact noncompact groups,  and let
$d_i:G_i\to [t_0,\infty)$, $i=1,2$, be proper continuous functions.
We set
\begin{align*}
v_i(t)&=\vol(\{g_i\in G_i;\, d_i(g_i)<t\}),\\
v(t)&=\vol(\{(g_1,g_2)\in G_1\times G_2;\, d_1(g_1)+d_2(g_2)<t\}).
\end{align*}

\begin{proposition}\label{p:prod1}
Suppose that there exist $c>0$ and $s_0>t_0$ such that  
for all sufficiently small $\epsilon>0$ and for all $t>s_0$,
$$
v_i(t+\epsilon)\le (1+c\epsilon)\, v_i(t),\quad i=1,2.
$$
Then for all sufficiently small $\epsilon>0$ and for all $t>2s_0+2$,
$$
v(t+\epsilon)\le (1+3c\epsilon)\, v(t).
$$
\end{proposition}

\begin{proof}
Let
\begin{align*}
w_1(t)&:=\vol(\{(g_1,g_2);\, d_1(g_1)+d_2(g_2)<t, d_1(g_1)\ge s_0+1, d_2(g_2)<s_0+1\}),\\
w_2(t)&:=\vol(\{(g_1,g_2);\, d_1(g_1)+d_2(g_2)<t, d_1(g_1)< s_0+1, d_2(g_2)\ge s_0+1\}),\\
w_3(t)&:=\vol(\{(g_1,g_2);\, d_1(g_1)+d_2(g_2)<t, d_1(g_1)\ge s_0+1, d_2(g_2)\ge s_0+1\}).
\end{align*}
For sufficiently small $\epsilon$ and for all $t$, we have
\begin{align*}
&w_1(t+\epsilon)-w_1(t)\\
\le& \int_{g_2:\, d_2(g_2)<s_0+1} \vol(\{g_1: \max\{s_0+1,t-d_2(g_2)\}\le d_1(g_1)<t+\epsilon-d_2(g_2)\})\,dg_2\\
\le & \int_{g_2:\, t-d(g_2)>s_0} (v_1(t-d_2(g_2)+\epsilon)-v_1(t-d_2(g_2))\,dg_2\\
\le & c\epsilon\, \int_{g_2:\, t-d(g_2)>s_0} v_1(t-d_2(g_2))\,dg_2\le c\epsilon\,v(t).
\end{align*}
Using similar argument, one shows that
$$
w_i(t+\epsilon)-w_i(t)\le c\epsilon\,v(t),\quad i=1,2,3,
$$
for sufficiently small $\epsilon$ and for all $t$.
Since for $t>2s_0+2$,
$$
v(t+\epsilon)-v(t)=(w_1(t+\epsilon)-w_1(t))+(w_2(t+\epsilon)-w_2(t))+(w_3(t+\epsilon)-w_3(t)),
$$
this implies the claim.
\end{proof}
A very similar argument establishes the following :

\begin{proposition}\label{p:prod2}
Suppose that there exist $c>0$ and $s_0>t_0$ such that  
for all $t>s_0$,
$$
v_i(t+1)\le c\, v_i(t),\quad i=1,2.
$$
Then for all all $t>2s_0+2$,
$$
v(t+1)\le (1+3c)\, v(t).
$$
\end{proposition}
%

%
\begin{proposition}\label{p:prod3}
Suppose that there exist $c>0$ such that  
for all sufficiently small $\epsilon>0$ and for all $t\ge t_0$,
$$
v_1(t+\epsilon)-v_1(t)\le c\epsilon \max\{v_1(t),1\}.
$$
Then there exists $s_0\ge 0$ such that 
for all sufficiently small $\epsilon>0$ and for all $t\ge 2t_0$,
$$
v(t+\epsilon)-v(t)\le c\epsilon\, v(t+s_0).
$$
\end{proposition}

\begin{proof}
Since
$$
v(t)=\int_{G_2} v_1(t-d_2(g_2))\, dg_2=\int_{g_2:\,t-d_2(g_2)\ge t_0} v_1(t-d_2(g_2))\, dg_2,
$$
it follows that for all sufficiently small $\epsilon>0$ and $t\ge 2t_0$,
\begin{align*}
v(t+\epsilon)-v(t) &\le c\epsilon\, \int_{g_2:\,t-d_2(g_2)\ge t_0}\max\{v_1(t-d_2(g_2)),1\}\, dg_2\\
&\le c\epsilon\, \int_{g_2:\,t-d_2(g_2)\ge t_0}\max\{v_1(t-d_2(g_2)),1\}\, dg_2.
\end{align*}
Since $G_1$ is noncompact, $v(t)\to\infty$ as $t\to\infty$ and 
there exists $s_0>0$ such that for all $t>s_0+t_0$, we have $v_1(t)>1$.
Then
\begin{align*}
v(t+\epsilon)-v(t) &\le c\epsilon\, \int_{g_2:\,t-d_2(g_2)\ge t_0}v_1(s_0+t-d_2(g_2))\, dg_2\\
&\le c\epsilon\, v(t+s_0).
\end{align*}
\end{proof}

\begin{proof}[Proof of Theorem \ref{Appendix:vol}]
(1) follows from Propositions \ref{p:metric} and \ref{p:prod1}.
(2) follows from Propositions \ref{p:e} and \ref{p:prod1}.
(3) follows from Propositions \ref{p:metric2}, \ref{p:prod2} and \ref{p:prod3}.
(4) follows from Propositions \ref{p:e2}, \ref{p:prod2} and \ref{p:prod3}.
\end{proof}

As we saw, the properties of balancedness and well-balancedness play an important role 
in the proofs of  the ergodic theorems.  To complete our discussion of the averages discussed in
in Theorem \ref{Appendix:vol}, let us note the following.   

First, the following criterion is sufficient to establish  that  the averages $\left(\sum_i a_i d_i^p\right)^p$ are in fact well-balanced, provided $1 < p < \infty$.

\begin{proposition}\label{p:balan}
Suppose that for some $s_0>t_0$, $a_i,b_i>0$, $u_i\ge 0$ and $w_i>0$,
$$
a_i t^{u_i} e^{w_it}\le v_i(t)\le  b_i t^{u_i} e^{w_it}\quad \hbox{for $t>s_0$ and $i=1,2$.}
$$
Then the sets 
$$
G_t=\{(g_1,g_2):\, d_1(g_1)^{p}+d_2(g_2)^{p}< t^p\}
$$
are well-balanced.
\end{proposition}

\begin{proof}
For $s\in [0,1]$,
$$
\vol(G_t)\ge v_1((1-s^p)^{1/p}t)v_2(st)\gg (1-s^p)^{u_1/p}s^{u_2}t^{u_1+u_2}\exp (\kappa(s)t) 
$$
where $\kappa(s)=w_1 (1-s^p)^{1/p}+w_2s$. It follows from convexity of $\kappa$ that for some
$s_0\in (0,1)$, $\kappa(s_0)> w_1, w_2$. This implies the claim.
\end{proof}

Typically, the averages defined by the distances $\sum_i a_i d_i$  are not balanced. However, we have

\begin{proposition}\label{p:balan2}
Under the assumptions of Proposition \ref{p:balan},  there exist $\alpha_1,\alpha_2>0$
such that the sets 
$$
G_t=\{(g_1,g_2):\, \alpha_1 d_1(g_1)+\alpha_2 d_2(g_2)< t\}
$$
are balanced.
\end{proposition}

\begin{proof}
Choosing $\alpha_1,\alpha_2>0$ suitably and rescaling the distance functions
we may assume that $\alpha_1=\alpha_2=1$ and $w_1=w_2$. We have
$$
\vol(G_t)\ge v_1(t/2)v_2(t/2)\gg t^{u_1+u_2}\exp (w_1t). 
$$
This implies the claim unless $u_1=u_2=0$. In this case,
$$
\vol(G_t)\gg \int_{t-d_2(g_2)>s_0} e^{w_1(t-d_2(g_2))}\,dg_2.
$$
Since $v_2(t)\gg e^{w_2 t}$, we have $\int_{G_2} e^{-w_2d_2(g)} dg=\infty$.
This implies the proposition.
\end{proof}

\subsection{Admissible, well-balanced, boundary-regular families} 
The present subsetion is devoted to the proof of Theorem \ref{CAT}, whose formulation we recall.    
\vskip0.1cm
{\bf Theorem 3.17.} {\it Let $G=G_1\cdots G_s$ be an $S$-algebraic group and
$\ell_i$ denote the standard $CAT(0)$-metric on either
the symmetric space $X_i$ or the Bruhat--Tits building $X_i$
associated to $G_i$. For $p>1$ and $u_i\in X_i$, define
$$
G_t=\{(g_1,\ldots, g_s):\, \sum_i \ell_i(u_i, g_iu_i)^p<t^p\}.
$$
Let $m$ be a Haar measure $G$.

 \begin{enumerate}
\item[(i)] 
There exist $\alpha,\beta>0$ such that for every nontrivial projection $\pi:G\to L$,
$$
m(G_t\cap \pi^{-1}(L_{\alpha t}))\ll e^{-\beta t}\cdot m_t(G_t)\,,
$$
namely the averages are well-balanced. 
\item[(ii)] If $G$ has at least one Archimedian factor, then the family $G_t$ is admissible, and 
writing $m=\int_0^\infty m_t\, dt$ where $m_t$
is a measure supported on $\partial G_t$, the following estimate holds : 

There exist $\alpha,\beta>0$ such that for every nontrivial projection $\pi:G\to L$,
$$
m_t(\partial G_t\cap \pi^{-1}(L_{\alpha t}))\ll e^{-\beta t}\cdot m_t(\partial G_t)\,,
$$
namely the averages are boundary-regular.

\end{enumerate}
}

In the proof, we use the following lemma:

\begin{lemma}\label{l:v}
With notation as in Theorem \ref{CAT},
there exists $\eta >0$ such that for every $\epsilon>0$ and $t\gg 0$,
$$
e^{(\eta-\epsilon)t}\ll_\epsilon m(G_t)\ll_\epsilon e^{(\eta+\epsilon) t}.
$$
\end{lemma}

\begin{proof}
Note that the stabilizer of $u_i$ is a maximal compact subgroup $K_i$ of $G_i$.
For archimedian factors, one can choose
a Cartan subgroup $A_i$ of $G_i$ equipped with a scalar product such that
the map $a\mapsto a\cdot u_i$, $a\in A_i$, is an isometry. 
Then Cartan decomposition $G_i=K_iA^+_iK_i$ holds, and
a Haar measure on $G_i$ is given by $dk_1d\nu_i(a)dk_2$ where
\begin{equation}\label{eq:arch}
d\nu_i(a)
=\prod_{\alpha\in \Sigma_i^+} \sinh (\alpha(a))^{n_{i,\alpha}}da,
\end{equation}
$\Sigma_i^+$ denotes the set of positive roots, and $n_{i,\alpha}$ is the dimension of the root space.
For non-Archimedean factors, there exists a lattice $A_i$ in the centralizer of a maximal split torus,
equipped with a scalar product such that the map $a\mapsto a\cdot u_i$, $a\in A_i$, is an isometry 
and $G_i=K_iA^+_iK_i$. Let
$$
d\nu_i=\sum_{a\in A_i^+} \hbox{vol}(K_iaK_i)\delta_a.
$$
Note that
\begin{equation}\label{eq:rho}
q_i^{2\rho_i(a)}\ll\hbox{vol}(K_iaK_i)\ll q_i^{2\rho_i(a)}
\end{equation}
where $q_i$ is the order of the residue field and $2\rho_i$ is the sum of positive roots.
Consider a measure $\nu=\otimes_i \nu_i$ on $A^+=\prod_i A^+_i$,
and set $d(a)=(\sum_i \|a_i\|_i^p)^{1/p}$. 
Let $\tilde A^+=\prod_i (A_i^+\otimes \mathbb{R})$,
$d\tilde\nu_i(a)=q_i^{2\rho_i(a)}da$ for non-Archimedean factors,
$\tilde\nu_i=\nu_i$ for Archimedean factors, and
$\tilde\nu=\otimes_i \tilde\nu_i$.
It follows from (\ref{eq:rho}) that there exists $c>0$ such that for $t\gg 0$,
$$
\tilde \nu(\{a\in \tilde A: d(a)<t-c\})\ll\nu(\{a\in A: d(a)<t\})\ll \tilde \nu(\{a\in \tilde A: d(a)<t+c\}).
$$
Hence, the proof is reduced to estimation of the integral given by 
$I(t):=\int_{a\in \tilde{A}^+:\,d(a)<t} d\tilde \nu(a)$.
Let 
$$
2\rho=\sum_i (\log q_i)\, 2\rho_i\quad\hbox{and}\quad
\eta=\sup\{2\rho(a):\, a\in \tilde A^+, d(a)< t\}
$$
(here we set $q_i=e$ for Archimedean factors). We have
\begin{equation}\label{up}
I(t)\ll t^{(\dim \tilde A^+-1)} e^{\eta t}.
\end{equation}
For every $\epsilon>0$, there exists a closed cone $C$ 
contained in the interior of $\tilde A^+$ such that for $a\in C$, $2\rho(a)>(\eta-\epsilon)d(a)$.
Then
$$
I(t)\gg \int_{a\in C:\,d(a)<t} e^{2\rho(a)}\, da\gg e^{(\eta-\epsilon)t}.
$$
This implies the claim.
\end{proof}

\begin{proof}[Proof of Theorem \ref{CAT}]
To prove part (i), 
let $G=LL'$ be a nontrivial decomposition of $G$. 
Using Lemma \ref{l:v}, we deduce that for $s,\epsilon\in (0,1)$ and $t\gg 0$,
$$
m(G_t)\ge m_L(L_{(1-s^p)^{1/p}t})m_{L'}(L'_{st})\gg_\epsilon \exp (\kappa(s,\epsilon)t) 
$$
where $\kappa(s,\epsilon)=(\eta-\epsilon) (1-s^p)^{1/p}+(\eta^\prime-\epsilon)s$.
It follows from convexity of $\kappa(\cdot,\epsilon)$ that for some
$s_0,\epsilon_0\in (0,1)$, we have $\kappa(s_0,\epsilon_0)> \eta, \eta^\prime$.
Hence, there exists $\beta>0$ such that for every $t\gg 0$,
\begin{equation}\label{e1}
m_L(L_{t})\ll e^{-\beta t}\cdot m(G_t). 
\end{equation}
For  $\alpha>0$, and $t\gg 0$, we have
$$
m(G_t\cap L_{\alpha t}L'))\le m_L(L_{\alpha t})m_{L'}(L'_t)\ll m_L(L_{\alpha t})e^{-\beta t} \cdot m(G_t).
$$
This implies that the averages in question are well-balanced.

As to part (ii), namely admissibility, the property that 
$$
\mathcal{O}_\epsilon G_t\mathcal{O}_\epsilon\subset G_{t+c\epsilon}\quad
\hbox{for some $c>0$ and every $\epsilon, t>0$}
$$
follows from the triangle inequalities for $d_i$'s and the $L^p$-norm. Now we show that
\begin{equation}\label{eq:adm}
m(G_{t+\epsilon})\le (1+c\epsilon)m(G_t)\quad
\hbox{for some $c>0$ and every $t\gg 0$, $\epsilon\in (0,1)$}.
\end{equation}
Let $d(g)=(\sum_i \ell_i(u_i,gu_i)^p)^{1/p}$.
Write $G=MN$ where $M$ is an archimedian factor and $N$ is its complement.
Setting
$$
v(t)=m_M(M_t)\quad\hbox{and}\quad w(t,n)=v((t^p-d(n)^p)^{1/p}),
$$
we have
$$
m(G_t)=\int_{d(n)<t} w(t,n)\, dn.
$$
We claim that the function $v$ is differentiable and
\begin{equation}\label{eq:vv}
v'(t)\ll\max\{t,v(t)\}\quad\hbox{for all $t>0$.}
\end{equation}
To prove this, we consider the Cartan decomposition $M=KA^+K$ (as in the proof of Lemma \ref{l:v})
and introduce polar coordinates $(r,\omega)\in\mathbb{R}^+\times S^+$ on the Lie
algebra of $A$. The Haar measure is given by $\xi(r,\omega)\,dk_1drd\omega dk_2$
with explicit density $\xi$ (see (\ref{eq:arch})). We get
\begin{equation}\label{comp}
m_M(M_{t+\epsilon})-m_M(M_t)=\int_{S^+}\int_t^{t+\epsilon} \xi(r,\omega)\,drd\omega=
\epsilon\int_{S^+} \xi(\sigma(\omega),\omega)\,d\omega
\end{equation}
for some $\sigma(\omega)\in [t,t+\epsilon]$. Since $\xi(r,\omega)\ll r$
for $r\in [0,1]$, this implies (\ref{eq:vv}) for $t\in [0,1]$.
To establish (\ref{eq:vv}) for $t>1$, we use the property of the function $\xi$ stated in (8.1). 
It follows from (\ref{eq:vv}) that
$$
w'(t,n)\ll\max\{1,w(t,n)\}.
$$
uniformly on $t>0$ and $n\in N$ satisfying $d(n)<t$, and we deduce the estimate
$$
m(G_{t+\epsilon})-m(G_{t})
\ll \epsilon\int_{d(n)<t} \max\{1,w(t+\epsilon,n)\}\, dn
$$
for every $t>0$ and $\epsilon\in (0,1)$. There exists $t_0>1$ such that
$w(t+t_0,n)>1$
for every $t>0$ and $n\in N$ such that $d(n)<t$. Then
$$
m(G_{t+\epsilon})-m(G_{t})
\ll \epsilon\int_{d(n)<t} w(t+t_0,n)\, dn\le \epsilon\,m(G_{t+t_0}).
$$
Now (\ref{eq:adm}) follows from Lemma \ref{l:v}. This proves that the sets $G_t$ are admissible.
We have decomposition of the Haar measure on $M$: $m_M=\int_0^\infty m_{M,t}\, dt$,
where $m_{M,t}$ is measure supported on $\partial M_t$.
Note that
$$
m_{M,t}(\partial M_t)=v'(t).
$$
The Haar measure on $G$ has a decomposition: $m=\int_0^\infty m_{t}\, dt$,
where
$$
dm_{t}(m,n)=t^{p-1}(t^p-d(n)^p)^{1/p-1}dm_{M,(t^p-d_N(n)^p)^{1/p}}(m)dn.
$$
Hence, we have
$$
m_{t}(\partial G_t)=\int_{d(n)<t} w'(t, n)dn.
$$
As in (\ref{up}),
$$
m_M(M_t)=\int_{S^+}\int_0^t \xi(s,\omega)dsd\omega\ll t^{\dim A-1} e^{\eta t}
$$
where $\eta=\max\{2\rho(\omega):\omega\in S^+\}$ and $2\rho$ is the sum of positive roots of $A$.
Choosing a small neighborhood $U$ of $\omega_0$ satisfying $2\rho(\omega_0)=\eta$, we deduce 
that using (\ref{comp}) that for every $\delta>0$ and $t\gg 0$,
$$
m_M(M_{t+\epsilon})-m_M(M_t)\ge 
\epsilon\int_{U} \xi(\sigma(\omega),\omega)\,d\omega
\gg \epsilon\cdot e^{(\eta-\delta)t},
$$
This implies that for every $\delta>0$ and $t\gg 0$,
$v'(t)\gg v((1-\delta) t).$
Hence, for $\alpha\in (0,1)$ and $t\gg 0$,
\begin{align}\label{eq:part}
m_{t}(\partial G_t)&\ge\int_{d(n)<\alpha t} w'(t, n)\,dn\gg \int_{d(n)<\alpha t} w(t, n)\,dn \ge m(G_{\alpha t}).
\end{align}
For $\alpha\in (0,1)$ and $t\gg 0$,
\begin{align*}
m_{t}(\partial G_t\cap M_{\alpha t}N)&=
\int_{(1-\alpha^p)^{-1/p}t<d(n)<t} w'(t, n)\,dn\\
&\ll \int_{(1-\alpha^p)^{-1/p}t<d(n)<t} \max\{1,w(t, n)\}\,dn\\
&\le m_N(N\cap G_t)+m(G_t\cap M_{\alpha t}N).
\end{align*}
Also, for every nontrivial simple factor $\pi:N\to L$,
\begin{align*}
m_{t}(\partial G_t\cap M \pi^{-1}(L_{\alpha t}))&=
\int_{n\in \pi^{-1}(L_{\alpha t}),d(n)<t} w'(t, n)\,dn\\
&\ll \int_{n\in \pi^{-1}(L_{\alpha t}),d(n)<t} \max\{1,w(t, n)\}\,dn\\
&\le m_N(\pi^{-1}(L_{\alpha t})\cap G_t)+m(G_t\cap M\pi^{-1}(L_{\alpha t})).
\end{align*}
Now boundary-regularity follows from (\ref{e1}), (i), (\ref{eq:part}).
%
\end{proof}

\subsection{Admissible sets on principal homogeneous spaces}

We now consider sets defined by a norm on principal homogeneous spaces,  which appear in the 
discussion of integral equivalence of forms in two or more variables in \S 2.3. 

\begin{proposition}\label{p:admis}
Let $G$ be a connected semisimple Lie group with finite center,
$\rho: G\to \hbox{\rm GL}(V)$ an irreducible representation, and $v_0\in V$ with compact stabilizer.
We fix a norm on $V$ and set
$$
G_t=\{g\in G: \log \|\rho(g)v_0\|< t\}
$$
Let $\pi$ denote the projection on the highest weight space.
If $0\notin \pi(\rho(K)v_0)$ for a maximal compact subgroup $K$ of $G$, then the sets $G_t$ are admissible.
\end{proposition}

In particular, the proposition applies to the following example:
the group $G=\hbox{SL}_k(\mathbb{R})$ acting on the space $W_{n,k}$ of homogeneous polynomials
of degree $n$ and $f_0\in W_{n,k}$ is such that $f(x)\ne 0$ for all $x\in\mathbb{R}^d\backslash \{0\}$.

In the proof we use the following lemma.

\begin{lemma}[cf. \cite{EMS}, Lemma A.4]\label{l:emss}
Let $\lambda_i(\omega)\in \mathbb{R}$ and $v_i(k)\in V$ depend on parameters $\omega, k$
and for every $s\ge 0$,
$$
\exp(s\max_i\lambda_i(\omega))\ll \left\|\sum_i e^{s\lambda_i(\omega)} v_i(k)\right\|\ll \exp(s\max_i\lambda_i(\omega))
$$
uniformly on $\omega$, $k$. Then there exists $T_0>0$ such that for every $T>T_0$, the set
$$
\left\{s\ge 0: \left\|\sum_i e^{s\lambda_i(\omega)} v_i(k)\right\|<T\right\}
$$
is an interval $[0, r(T,\omega, k)]$ and
\begin{equation}\label{eq:lip}
r((1+\epsilon)T,\omega, k)-r(T,\omega, k)\ll \epsilon
\end{equation}
uniformly on $\omega$, $k$, $\epsilon\in (0,1)$, and $T>T_0$.
\end{lemma}

We note that the statement of Lemma A.4 in \cite{EMS} is somewhat weaker than the statement above,
but the proof there implies the lemma in this generality.

\begin{proof}[Proof of Proposition \ref{p:admis}] 
We fix a Cartan decomposition $G=KA^+K$.
For a weight $\lambda$ of $\hbox{Lie}(A)$, we denote by $\pi_\lambda$ the projection on the weight space of $\lambda$.
Then for $g\in k_1\exp(a)k_1\in G$,
$$
gv_0=\sum_\lambda e^{\lambda(a)} k_1 \pi_\lambda(k_2v_0) 
$$
This implies that 
$$
\max_{\lambda,\; k_2} e^{\lambda(a)}\|\pi_\lambda(k_2v_0)\|\ll
\|gv_0\|\ll \max_{\lambda,\; k_2} e^{\lambda(a)}\|\pi_\lambda(k_2v_0)\|,
$$
and it follows from the assumption $0\notin \pi(Kv_0)$ that
\begin{equation}\label{eq:estt}
e^{\lambda_{max}(a)}\ll\|gv_0\|\ll e^{\lambda_{max}(a)}.
\end{equation}

It is straightforward to check that $\mathcal{O}_\epsilon G_t\subset G_{t+c\epsilon}$
for some $c>0$.  Now we show that $G_t\mathcal{O}_\epsilon\subset G_{t+c\epsilon}$ as well.
For $g=k_1\exp(a)k_2$ and $h$ in $G$, we have
\begin{align*}
\|ghv_0-gv_0\|&\ll \sum_\lambda e^{\lambda(a)}\|\pi_\lambda(k_2(hv_0-v_0))\|\ll e^{\lambda_{max}(a)}\|hv_0-v_0\|\\
&\ll \|gv_0\| d(h, e).
\end{align*}
This estimate implies the claim.

It remains to show that for sufficiently small $\epsilon>0$ and for sufficiently large $t$,
\begin{equation}\label{eq:ad}
m_G(G_{t+\epsilon})-m_G(G_t)\ll \epsilon.
\end{equation}
Since (\ref{eq:estt}) holds, we may apply Lemma \ref{l:emss} with
vectors $v_i's$ given by $k_1\pi_\lambda(k_2 v_0)$ with $k_1,k_2\in K$.
This implies the Lipschitz estimate (\ref{eq:lip}), and 
(\ref{eq:ad}) is deduced as in \cite[Proposition A.5]{EMS}.
\end{proof}

\begin{proposition}\label{p:admis2}
Under the assumptions of Proposition \ref{p:admis},
there exist $c>0$, $a>0$, and $b=1,\ldots, \hbox{\rm rank}(G)$ such that
$$
\hbox{\rm vol}(G_t)\sim ct^{b-1}e^{at}\quad\hbox{as $t\to\infty$.}
$$
\end{proposition}

\begin{proof}
Fix a Cartan decomposition $G=KA^+K$. 
Then the Haar measure on $G$ is given by $\xi(a)dk_1dadk_2$.
For $k_1,k_2\in K$, set
$$
A_t(k_1,k_2)=\{a\in A:\, \log \|k_1ak_2\|< t\}.
$$ 
We have
$$
\hbox{\rm vol}(G_t)=\int_{K\times K}\int_{A_t(k_1,k_2)} \xi(a)da dk_1dk_2.
$$
By \cite[\S 7]{GW},
$$
\int_{A_t(k_1,k_2)} \xi(a)da\sim c(k_1,k_2) t^{b(k_1,k_2)} e^{a(k_1, k_2)t}\quad\hbox{as $t\to\infty$.}
$$
Also, it follows from (\ref{eq:estt}) that
$$
t^{b-1} e^{at}\ll \int_{A_t(k_1,k_2)} \xi(a)da\ll t^{b-1} e^{at}
$$
for sufficiently large $t$ and $k_1,k_2\in K$.
Hence, the parameters $a(k_1,k_2)$ and $b(k_1,k_2)$ are constant,
and the claim follows from the dominated convergence theorem.
\end{proof}

\subsection{Tauberian arguments and H\"older continuity}

Finally, we will now establish the H\"older-admissibility property of averages defined 
by a regular proper function on an algebraic varieties with a regular volume form. 
Our approach uses the following Tauberian theorem
which is proved using the argument of \cite{CT}, Theorem A.1.

\begin{proposition}\label{p:tschinkel}
 Let $v: [0,\infty)\to [0,\infty)$ and $f(s)=\int_0^\infty x^{-s} v(x)\,dx$.
\begin{enumerate}
 \item Let $v(t)$ be increasing for sufficiently large $t$.
Assume that the integral $f(s)$ converges for $\Re(s)\gg 0$,
admits meromorphic continuation to $\Re(s)>a-\delta_0$,
and in this domain it has unique pole $s=a$ of multiplicity $b$
and satisfies
$$
\left|f(s)\frac{(s-a)^b}{s^b}\right|=O((1+\Im(s))^\kappa)
$$
for some $\kappa>0$.
Then
$$
v(t)=t^{a-1}P(\log t)+O(t^{a-1-\delta})\quad \hbox{as $t\to\infty$}
$$
for some nonzero polynomial $P$ and $\delta>0$.

\item 
Let $v(t)$ be increasing for sufficiently small $t$.
Assume that
the integral $f(s)$ converges for $\Re(s)\ll 0$,
admits meromorphic continuation to $\Re(s)<a+\delta_0$,
and in this domain it has unique pole $s=a$ of multiplicity $b$
and satisfies
$$
\left|f(s)\frac{(s-a)^b}{s^b}\right|=O((1+\Im(s))^\kappa)
$$
for some $\kappa>0$.
Then
$$
v(t)=t^{a-1}P(\log t)+O(t^{a-1+\delta})\quad \hbox{as $t\to 0^+$}
$$
for some nonzero polynomial $P$ and $\delta>0$.
\end{enumerate}
\end{proposition}

\begin{proof}
The first statement is essentially proved in \cite{CT} (in the context of Dirichlet series), and
the second statement is proved similarly. We give a sketch of the proof for the second statement.

For negative $a'<a$, define
$$
w_k(t)=\frac{(-1)^{k+1} k!}{2\pi i} \int_{a'+i\mathbb{R}} f(s) t^s\, \frac{ds}{s^{k+1}}.
$$
This integral is absolutely convergent for $k>\kappa$.
Applying Cauchy formula for the region $a'<\Re(s)<a+\delta_0/2$, $\Im(s)\le S$
with $S\to\infty$ we deduce that for some nonzero polynomial $P_k$,
\begin{equation}\label{eq:v_k}
w_k(t)=t^aP_k(\log t)+O(t^{a+\delta_0/2})\quad \hbox{as $t\to 0^+$.}
\end{equation}
It follows from the formula
$$
\int_{a'+i\mathbb{R}} \lambda^s\, \frac{ds}{s^{k+1}}=\left\{
\begin{tabular}{cl}
$-\frac{2\pi i}{k!}(\log\lambda)^k$ & $0<\lambda\le 1$,\\
0 & $\lambda<1$
\end{tabular}
\right.
$$
that
$$
w_k(t)=(-1)^k\int_t^\infty (\log (t/x))^k\, v(x)dx=\int_t^\infty (\log (x/t))^k\, v(x)dx.
$$
Now we derive asymptotic expansion for $w_{k-1}$ assuming that (\ref{eq:v_k}) holds.
By the intermediate value theorem, for every $t>0$ and $\eta\in (0,1)$,
\begin{align*}
w_{k-1}(t)
&\le \frac{\int_t^\infty \left((\log(x/t(1-\eta)))^k -(\log(x/t))^k\right)\,v(x)dx}{-k\log(1-\eta)}\\
&\le \frac{w_k(t(1-\eta))-w_k(t)}{-k\log(1-\eta)},
\end{align*}
and
\begin{align*}
w_{k-1}(t)
&\ge \frac{\int_{t(1+\eta)}^\infty \left((\log(x/t))^k -(\log(x/t(1+\eta))^k\right)\,v(x)dx}{k\log(1+\eta)}\\
&\quad +\frac{1}{\log(1+\eta)}\int_t^{t(1+\eta)} (\log(x/t))^k\, v(x)dx\\
&\ge \frac{w_k(t)-w_k(t(1+\eta))}{k\log(1+\eta)}.
\end{align*}
Taking $\eta=t^\epsilon$ with small $\epsilon>0$ and using (\ref{eq:v_k}), we deduce that
from the above estimates that
$$
w_{k-1}(t)=t^aP_{k-1}(\log t)+O\left((\log t)^{\deg P_k} t^{a+\epsilon}+ t^{a+\delta_0/2-\epsilon}\right)\quad \hbox{as $t\to 0^+$}
$$
for some nozero polynomial $P_{k-1}$ (see \cite{CT}, proof of Theorem A.1, for a detailed computation).
This implies that (\ref{eq:v_k}) holds for all $k\ge 0$.
To complete the proof, we observe that for small $t>0$ and $\eta\in (0,1)$,
$$
\frac{1}{t\eta}(w_0(t(1-\eta))-w_0(t))\le v(t)\le \frac{1}{t\eta}(w_0(t)-w_0(t(1+\eta))).
$$
Setting $\eta=t^\epsilon$ with small $\epsilon>0$, we deduce the required
asymptotic expansion for $v(t)$ from the asymptotic expansion for $w_0(t)$.
\end{proof}

We will now employ Proposition \ref{p:tschinkel}
 and prove the H\"older continuity of the volume function in the context 
of algebraic functions on algebraic varieties. 

\begin{theorem}\label{alg-holder}
Let $X$ be a real algebraic variety equipped with a regular volume form $\omega$
and $\Psi: X\to\mathbb{R}$ a nonconstant regular proper function.
Then the function $g(t)=\int_{\Psi(x)<t}d\omega$ is uniformly H\"older on finite intervals.
\end{theorem}

\begin{proof}
At regular values of $\Psi$,
$$
g'(t)=\int_{\Psi^{-1}(t)}\frac{d\nu_t}{\|\grad\Psi\|}
$$
where $\nu_t$ is the induced measure on the fiber $\Psi^{-1}(t)$. 
We claim that for $t$ in a neighborhood of an isolated critical value $t_0$,
\begin{equation}\label{eq:cl0}
g'(t)\ll |t-t_0|^{-r}
\end{equation}
for some $r>0$. Let $Z$ be the set of critical points of $\Psi$ in $\Psi^{-1}(t_0)$.
By Lojasiewicz's inequality (see, for example, \cite{BM}, Theorem 6.4),
$$
\|(\grad\Psi)_x\|\gg d(x, Z)^r
$$
for some positive $r$. For $x\in \Psi^{-1}(t)$ and $z\in \Psi^{-1}(t_0)$,
$$
d(x,z)\gg |\Psi(x)-\Psi(z)|=|t-t_0|.
$$
This implies (\ref{eq:cl0}), and by the intermediate value theorem,
for $t_0<t_1<t_2$ in a neighborhood of $t_0$,
\begin{equation}\label{eq:est1}
|g(t_2)-g(t_1)|\ll |t_1-t_0|^{-r}\cdot |t_2-t_1|.
\end{equation}

Next, we show that $g$ is H\"older at critical values of $\Psi$.
For $c>0$, we consider the function $v$ defined by
$v(t)=\int_{c\le \Psi <c+t} d\omega$
for $t<1$ and $v(t)=0$ for $t\ge 1$,
and its transform
$$
f(s)=\int_0^\infty t^{-s} v(t)dt=\int_{c\le \Psi(x)<c+1}\frac{1-(\Psi(x)-c)^{1-s}}{1-s}\, d\omega(x)
$$
which is absolutely convergent for $\Re (s)<1$. Applying Hironaka resolution of singularities to the function $F(g)=\Psi(x)-c$, we deduce
that there exists an atlas of maps $\phi_i:(-1,1)^d\to U_i$, $U_i$ is open in $X$,
such that $\phi_i$'s are diffeomorphisms on sets of full measure, and 
$$
F(\phi_i(x))=x^{\alpha_i}F_i(x)\quad\hbox{and}\quad d\omega(\phi_i(x))=x^{\beta_i}\rho_i(x)dx
$$
where $x^{\alpha_i}$ and $x^{\beta_i}$
denote monomials, and $F_i$ and $\rho_i$ are
positive smooth functions. Let $\{\eta_i\}$ be a partition of unity
subordinate to the cover $\{U_i\}$ such that
$$
\hbox{$\sum_i \eta_i=1$ on $\{c\le\Psi(x)\le c+1/2\}$ and
$\hbox{supp}(\eta_i)\subset \{\Psi(x)< c+1\}$}.
$$
We have
\begin{align}\label{eq:Long}
&\int_{c\le \Psi(x)<c+1}(\Psi(x)-c)^{1-s}\, dm_G(g)\\
=&\sum_i \int_{x\in (-1,1)^d:\,x^{\alpha_i}\ge 0} x^{\beta_i+(1-s)\alpha_i}F_i(x)^{1-s}\rho_i(x)\eta_i(\phi_i(x))\, dx+\xi(s)\nonumber
\end{align}
where $\xi(s)$ is an integral over the region $\Psi(x)>c+1/2$, hence, holomorphic,
and the other integrals can be meromorphic continued
integrating by parts:
\begin{align*}
&\int_{(0,1)^d} x^{\beta_i+(1-s)\alpha_i}F_i(x)^{1-s}\rho_i(x)\eta_i(\phi_i(x))\, dx\\
=&\left(\prod_j (1+\beta_i+ (1-s)\alpha_{i,j})\right)^{-1} \int_{(0,1)^d} x^{1+\beta_i+(1-s)\alpha_i}(F_i(x)^{1-s}\rho_i(x)\eta_i(\phi_i(x)))'\, dx.
\end{align*}
Therefore, (\ref{eq:Long}) implies that the conditions of Proposition \ref{p:tschinkel}(2)
are satisfied, and hence, 
$$
v(t)=t^{a-1}P(\log t)+O(t^{a-1+\delta})\quad \hbox{as $t\to 0^+$}
$$
for some $a\ge 1$. If $a=1$, then $\vol(\{x:\Psi(x)=c\})>0$, but
the set $\{x:\Psi(x)=c\}$ is a proper algebraic subvariety of $X$.
Hence, $a>1$ and we deduce the H\"older estimate
\begin{equation}\label{eq:H}
\vol(\{x: c\le \Psi(x)<c+t\})\ll t^{\alpha}.
\end{equation}
with $\alpha<a-1$.

Since $g$ is a polynomial function it has only finitely many critical values.
The function $g$ is $C^1$ on the set of regular values.
Hence, it remains to show that $g$ is H\"older in a neighborhood
of a critical values. For instance, consider the case 
when $t_1$ and $t_2$ are in a neighborhood of a critical value $t_0$
and $t_0<t_1<t_2$. The other cases are treated similarly.
We have
\begin{equation}\label{eq:est2}
|g(t_2)-g(t_1)|\ll (t_1-t_0)^\alpha+(t_2-t_0)^\alpha\ll (t_1-t_0)^\alpha+(t_2-t_1)^\alpha.
\end{equation}
When $(t_1-t_0)^{-r} \le (t_2-t_1)^{-1/2}$, the H\"older estimate follows from
(\ref{eq:est1}), and when the opposite inequality holds, the H\"older estimate follows from
(\ref{eq:est2}). This completes the proof.
\end{proof}

We now obtain the following 

\begin{theorem}\label{alg-vol}\label{p:h2}
Let $X$ be a real algebraic variety 
equipped with a regular volume form $\omega$, $\Psi:X\to [1,\infty)$ a proper function,
and
$$
v(t)=\vol(\{x\in X:\, \Psi(x)<t\}).
$$
Then for some $a\ge 1$, a nonzero polynomial $P$, and $\delta>0$, we have 
$$
v(t)=t^{a-1}P(\log t)+O(t^{a-1-\delta})\quad \hbox{as $t\to\infty$}.
$$
\end{theorem}

\begin{proof}

First note that since the function $\theta(t)=\max\{\|x\|: \Psi(x)<t\}$ is semialgebraic,
there exists $M>0$ such that  $\theta(t)\ll t^M$ for sufficiently large $t$.
This implies that for some $N>0$ and $t\gg 0$, we have $v(t)\ll t^M$.

Now let $w(t)=v(t)-v(1)$ for $t\ge 1$ and $w(t)=0$ for $t<1$.
Consider the transform of $w$:
$$
f(s)=\int_0^\infty t^{-s}w(t)dt
$$
which is convergent for $\Re(s)>M+1$ and
$$
f(s) 
=(s-1)^{-1} \int_{X} \Psi(x)^{-s+1}\, d\omega(x) 
$$

Applying Hironaka resolution of singularities, we may assume $X$ is 
semialgebraic subset of a smooth projective variety $Y$, and
there exists an atlas of maps $\phi_i:(-1,1)^d\to Y$, $U_i$ is open in $Y$,
such that $\phi_i$'s are diffeomorphisms on sets of full measure, 
and $\phi_i^{-1}(U_i\cap X)$ is a union of quadrants, and
$$
\Psi(\phi_i(x))^{-1}=x^{\alpha_i}\Psi_i(x),\quad\quad d\omega(\phi_i(x))=x^{\beta_i}\rho_i(x)dx,
$$
where $x^{\alpha_i}$ and $x^{\beta_i}$ denote monomials, and $F_i$ and $\rho_i$ are
smooth functions nonvanishing on $(-1,1)^d$. Let $\{\eta_i\}$ be a partition of unity
subordinate to the cover $\{U_i\}$. Then
\begin{align*}
\int_{X}\Psi(x)^{-s+1}\, d\omega(x)=\sum_i \int_{\phi_i^{-1}(U_i\cap X)} x^{(s-1)\alpha_i+\beta_i}\Psi_i(x)^{s-1}\rho_i(x)\eta_i(\phi_i(x))\, dx.
\end{align*}
Integrating by parts, we deduce that this expresion has meromorphic continuation
and satisfies the conditions of Proposition \ref{p:tschinkel}.
This implies the claim.
\end{proof}

\begin{theorem}\label{t:h} 
Let $X$ be a real algebraic variety equipped with a regular volume form $\omega$
and $\Psi: X\to [1,\infty)$ a nonconstant regular proper function.
Then for some $\beta>0$, the function $g(t)=\int_{\Psi<t}d\omega$ satisfies
$$
g((1+\epsilon)t)-g(t)\ll \epsilon^\beta \max\{1,g(t)\}\quad\hbox{for all $\epsilon\in (0,1)$ and $t\ge 0$.}
$$
\end{theorem}

\begin{proof}
On finite intervals, this is already proved in Theorem \ref{alg-holder}
Since $\Psi$ is regular, it has only finitely many critical points.
Thus, it remains to consider an interval $t\gg 0$ which contains
no critical points.  It follows from Theorem \ref{alg-vol}
that for $t\ge \epsilon^{-\alpha}$ with arbitrary $\alpha>0$, we have
$$
g((1+\epsilon)t)-g(t)\ll \epsilon^\beta\, g(t)
$$
where $\beta>0$ depends on $\alpha$. To prove the estimate for $t< \epsilon^{-\alpha}$,
we use that
$$
g'(t)=\int_{\Psi=t} \|(\grad \Psi)_x\|^{-1}\,d\omega_t 
$$
where $\omega_t$ is the regular volume form on $\{\Psi=t\}$ induced by $\omega$.
Since the function $t\mapsto \max\{\|(\grad\Psi)_x\|^{-1}: \Psi(x)=t\}$ is semialgebraic,
it follows that there exists $M>0$ such that for $\Psi(x)\gg 0$, 
$$
\|\grad (\log\Psi)_x\|^{-1}\ll \Psi(x)^M.
$$
Similarly, for $\Psi(x)\gg 0$, 
$$
\|x\|\ll \Psi(x)^N.
$$
This implies that for some $N>0$ and $t\gg 0$, we have
$$
g'(t)\ll t^{N}.
$$
Then when $0\ll t< \epsilon^{-\alpha}$ with $\alpha<1/N$, we have H\"older estimate
$$
g(t+\epsilon)-g(t)\ll \epsilon^{1-\alpha N}.
$$
Hence, the claim follows.
\end{proof}

Finally, we combine the foregoing arguments to prove H\"older-admissibility of families defined 
by a height function on a product of affine varieties, and in particular on $S$-algebraic group 
(as stated in Theorem \ref{HHA}).

\begin{theorem}\label{holder heights}
Let $X=X_1\cdots X_N$ be a product of affine varieties $X_i$ over local fields
equipped with regular volume forms.
We denote by $\|\cdot\|_i$ either the Euclidean norm for Archemedian places
or the $\max$-norm for non-Archemedian places and set
$$
X_t=\{(x_1,\ldots,x_N):\, \sum_i \log \|x_i\|_i<t\}.
$$
If at least one of the factors of $X$ is Archimedean, then the function $t\mapsto \vol(X_t)$ is
 uniformly H\"older.
\end{theorem}

\begin{proof}
Setting $v_i(t)=\vol(\{x_i: \log\|x_i\|<t\}$, the claim is deduced 
applying Propositions \ref{p:prod3} and \ref{p:prod2} inductively.
Using restriction of scalars, we can assume that all Archimedean factors are real.
For Archimedean $v_i$'s,
the assumption of Proposition \ref{p:prod3} follow from Theorem \ref{t:h} and the assumption of Proposition \ref{p:prod2} follows from Proposition \ref{p:h2}. Hence, it remains to verify
the assumption of Proposition \ref{p:prod2} at non-Archimedean places.
In this case, it follows from \cite{De} that $\int_{X_i} \|x\|_i^s\, d\omega_i(x)$
is a rational function of $q_i^{s}$ and $q_i^{-s}$ where $q_i$ is the order of the residue field.
Hence,
$$
w_i(n):=\vol(\{x\in X_i:\, \|x_i\|_i=q_i^n\})=\sum_j p_{ij}(n)q_i^{a_{ij}n}
$$
for rational polynomials $p_{ij}$ and $a_{ij}\in\mathbb{Z}$.
This implies that $v_i(t)=\sum_{n<t}w_i(n)$ satisfies
$v_i(t+1)\ll v_i(t)$ for sufficiently large $t$.
 
Now the claim follows from Propositions \ref{p:prod3} and \ref{p:prod2}.
\end{proof}

\end{document}